\documentclass[10pt]{jhuthesis}
\usepackage{amsmath,amssymb,amsthm}
\usepackage{graphicx}
\usepackage{local}
\usepackage[numbers,sort&compress]{natbib}
\usepackage{hypernat}
\usepackage{url}
\usepackage{enumerate}

\usepackage{memhfixc}

\begin{document}
                                %
\frontmatter

\author{Nevin Kapur}
\title{Additive functionals on random search trees}
\pubmonth{May}
\pubyear{2003}
\degree{Doctor of Philosophy}
\maketitle

\singlespacing


\chapter{Abstract}

Search trees are fundamental data structures in computer science.  We
study functionals on random
search trees that satisfy recurrence relations of a simple additive
form.  Many important functionals including the space requirement,
internal path length, and the so-called shape functional fall under
this framework.  Our goal is to derive asymptotics of moments and
identify limiting distributions of these functionals under two
commonly studied probability models---the random permutation model
and the uniform model.

For the random permutation model, our approach is based on establishing
\emph{transfer theorems} that link the order of growth of the input into
a particular (deterministic) recurrence to the order of growth of the
output.  For the
uniform model, our approach is based on
the complex-analytic tool of \emph{singularity analysis}.  To
facilitate a systematic analysis of these additive functionals we
extend singularity analysis, a class of methods by which one can
translate on a term-by-term basis an asymptotic expansion of a
functional around its dominant singularity into a corresponding
expansion for the Taylor coefficients of the function.  The most
important extension is the determination of how singularities are
composed under the operation of Hadamard product of analytic power
series.

The transfer theorems derived are used in conjunction with the method
of moments to establish limit laws for~\( m \)-ary search trees under
the random permutation model. For the uniform model on binary search
trees, the extended singularity analysis
toolkit is employed to establish the asymptotic behavior of the
moments of a wide class of functionals.  These asymptotics are used,
again in conjunction with the method of moments, to derive limit laws.

\bigskip\par\noindent
Advisor: James Allen Fill\\
Readers: James Allen Fill and John C.~Wierman



\chapter{Acknowledgments}

It is impossible for me to express satisfactorily in words the
gratitude I owe  my advisor, Jim Fill.  This dissertation is the
result of many, many sessions of joint toil where I learned the art
of research from him.  For that and much more, I am indebted to him.

This work would have been nothing more than a graduate student's
distant dream were it not for the contribution of Philippe Flajolet.
His outline concerning singularity analysis for Hadamard products
blossomed into this dissertation's central
Chapter~\ref{cha:sing-analys-hadam}, and his vision throughout the
project was invaluable.

I thank Professors Goldman, Priebe, Scheinerman, and Wierman for
serving on my dissertation defense committee and offering suggestions
to improve this work.

I am thankful to the faculty, staff, and students of the Department of
Mathematical Sciences at The Johns Hopkins University for their support.

My stay in Baltimore has been most enjoyable thanks to the
excellent friends I found here.  I will always have fond memories of
the times I spent with Carolyn Cooper, helen and Jason DeVinney, and
EJ Valitutto.

My mother taught me to value education and my parents sacrificed much
so that I could pursue my dream.  To them I owe much.

\bigskip
\hfill \textit{Nevin Kapur}


\clearpage
\tableofcontents
\clearpage
\listoffigures
\clearpage

\vspace*{3.15in}
\setlength{\epigraphwidth}{1.25\epigraphwidth}
\epigraphfontsize{\normalsize}
\epigraph{%
  ``The time has come,'' the Walrus said,\\
  ``To talk of many things: \ldots''}{\textit{The Walrus and The
  Carpenter}\\ \textsc{Lewis Carroll}}

\mainmatter
\pagestyle{simple}

\counterwithin{equation}{section}
\counterwithin{figure}{section}
\counterwithin{table}{section}

\part{Preliminaries}
\chapter{Introduction}
\label{cha:introduction}

\section{Background} \label{sec:intro_background}

We start by giving a brief overview of search trees, the basic object
under scrutiny in this dissertation.  Search trees are fundamental
data structures in computer science, their primary applications being
in searching and sorting.  For integer \( m \geq 2 \), the \( m \)-ary
search tree, or multiway tree, generalizes the binary search tree.
The quantity \( m \) is called the \emph{branching factor}.
According to~\cite{MR93f:68045}, search trees of branching factors
higher than 2 were first suggested by Muntz and
Uzgalis~\cite{muntz71:_dynam} ``to solve internal memory problems with
large quantities of data.''  For more background we refer the reader
to~ \cite{knuth97,knuth98} and~\cite{MR93f:68045}.

\label{def:marytree} An \emph{\( m \)-ary tree} is a rooted tree
with at most \( m \) ``children'' for each \emph{node (vertex)}, each
child of a node being
distinguished as one of \( m \) possible types.
Recursively expressed, an \( m \)-ary tree either is empty or consists
of a distinguished node (called the \emph{root}) together with an
ordered \( m \)-tuple of \emph{subtrees}, each of which is an \( m
\)-ary tree.

\label{def:marysearchtree} An \emph{\( m \)-ary search tree} is an
\( m \)-ary tree in which each node has the capacity to contain \(
m-1 \) elements of some linearly ordered set, called the set of
\emph{keys}.
In typical implementations of \(m\)-ary search trees, the keys at each
node are stored in increasing order and at each node one has \( m \)
pointers to the subtrees. By spreading the input data in \( m \)
directions instead of only 2, as is the case for a binary search tree,
one seeks to have shorter path lengths and thus quicker searches.

\begin{figure}[htbp]
\centering
\includegraphics{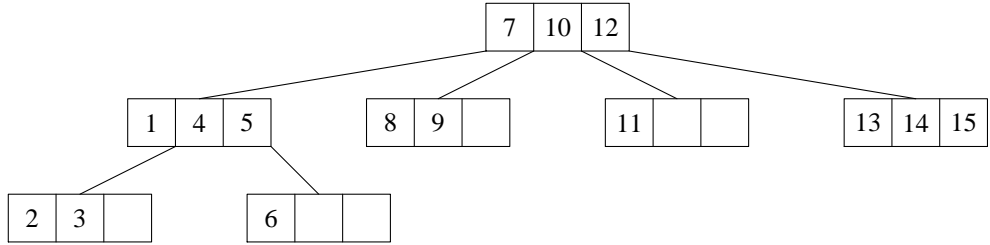}
\caption[A quaternary (i.e., ``$4$-ary'') tree.]{The quaternary (i.e.,
  ``$4$-ary'') tree associated with the sequence
  $s$.}\label{fig:examplemAry}
\end{figure}

We consider the space of \( m \)-ary search trees on \( n \) keys, and
assume that the keys are linearly ordered.  Hence, without loss
of generality, we can take the set of keys to be \( [n] :=
\{1,2,\ldots,n\} \). We construct an \( m \)-ary search tree from a
sequence $s$ of \( n \) distinct keys in the following way:
\begin{enumerate}
\item If \( n < m \), then all the keys are stored in the root node in
  increasing order.
\item If \( n \geq m \), then the first \( m-1 \) keys in the sequence
      are
  stored in the root in increasing order, and the remaining \( n-(m-1)
  \) keys are stored in the subtrees subject to the condition that if
  \( \sigma_1 < \sigma_2 < \cdots < \sigma_{m-1} \) denotes the
  ordered sequence of keys in the root, then the keys in the \( j \)th
  subtree are those that lie between \( \sigma_{j-1} \) and \(
  \sigma_{j} \), where \( \sigma_0 := 0 \) and \( \sigma_{m} := n+1
  \), sequenced as in $s$.
\item All the subtrees are \( m \)-ary search trees that satisfy
  conditions~1, 2, and~3.
\end{enumerate}
Figure~\ref{fig:examplemAry} shows the quaternary tree associated with
the sequence
\[
s := ( 10, 7, 12, 4, 1, 8, 5, 6, 9, 14, 11, 2, 15, 13, 3).
\]

\begin{figure}[htbp]
\centering
\includegraphics[width=\textwidth]{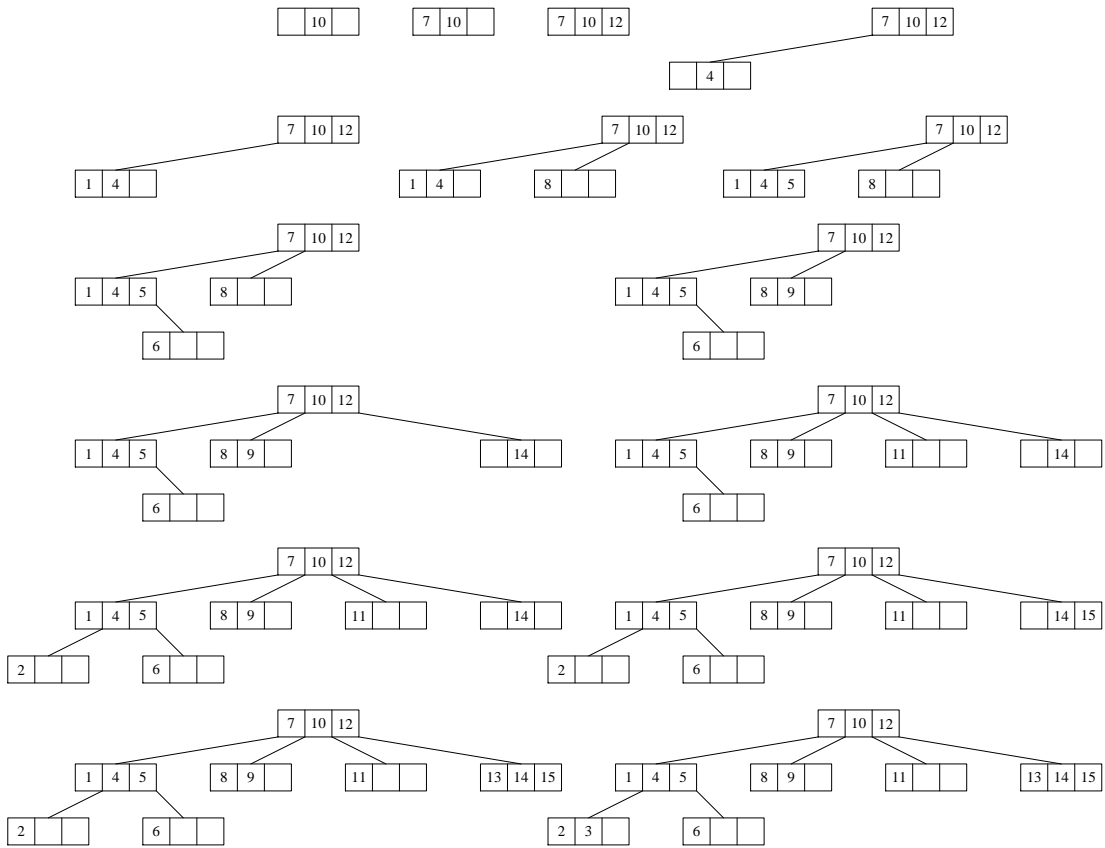}%
\caption[Sequential construction of a quaternary search
tree.]{Sequential construction of the quaternary search tree
  from the sequence $s$.}%
\label{fig:buildmAry}
\end{figure}

We may also build an \( m \)-ary search tree by inserting the successive
elements of a sequence into an initially empty tree. This process is
illustrated for the sequence \( s \) in Figure~\ref{fig:buildmAry}.

It is our goal to study functionals defined on \( m \)-ary search
trees under two natural probability models.  The functionals we
consider are of additive type, each induced by a ``toll function.''
Informally speaking, the value of an additive functional on
an~\(m\)-ary search tree is recursively defined as the sum of the
functional's value on the subtrees rooted at the children of the root of the
tree and a term (the \emph{toll}) that is a function of the size of
the tree.  Additive functionals defined in this manner represent the
cost of divide-and-conquer algorithms, where the inherent recursive
nature of the algorithms lends itself naturally to such a formulation.

\section{Probability models}
\label{sec:probability-models}

The two most common probability models on the space of \( m \)-ary
search trees are the \emph{uniform} model  and the
\emph{random permutation}, or \emph{random insertion}, model.  In the
case $m=2$, the number of search trees on $n$ keys is enumerated by
the $n$th Catalan number.  Thus the uniform model on binary search
trees is also referred to as the \emph{Catalan} model.

Under the uniform model, the probability of obtaining any particular
tree is just the reciprocal of the number of \( m \)-ary search trees
on \( n \) keys. Fill and Dobrow \cite{MR99f:68034} have studied
problems associated with determining this number.
Fill~\cite{MR97f:68021} argued that the functional corresponding to
the toll sequence \( (\ln{n})_{n \geq 1} \) serves as a crude measure
of the ``shape'' of a binary search tree, and explained how this
functional arises in connection with the move-to-root self-organizing
scheme for dynamic maintenance of binary search trees.  He obtained
asymptotic information about the mean and variance under the Catalan
model.  (The latter results were rederived in the
extension~\cite{MR99j:05171} from binary search trees to simply
generated rooted trees.)

In~\cite[Prop.~2]{MR88h:68033} Flajolet and Steyaert gave order-of-growth
information about the mean of functionals induced by tolls of the form \(
n^\alpha \).  Tak{\'a}cs established the limiting (Airy) distribution of
path length in Catalan trees~\cite{MR92m:60057,MR93e:60175,MR92k:60164},
which is the additive functional for the toll~$n-1$.  
The additive
functional for the toll $n^2$ arises in the study of the Wiener index of the
tree and has been analyzed by Janson~\cite{janson:_wiener}.

The uniform model on binary search trees has been also been used
recently by Janson~\cite{janson02:_ideal}
in the analysis of an
algorithm of Koda and Ruskey~\cite{MR94f:94015} for listing ideals in a
forest poset.

\begin{figure}[htbp]
\centering
\includegraphics{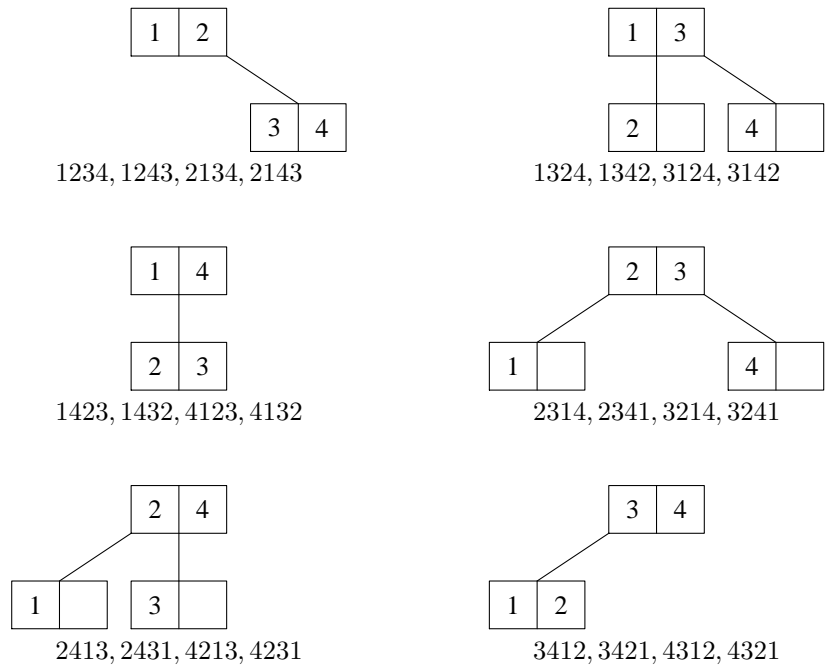}%
\caption{The distinct
  ternary trees and the permutations that generate them, for $n=4$.}%
\label{fig:distribution}
\end{figure}

Now we turn our attention to the random permutation model.  Let \( \pi
\) be a random permutation of \( [n] \):\ each of the \( n! \)
permutations is equally likely to be the value of $\pi$.  The
distribution of trees under the random permutation model is the
distribution induced by constructing an \( m \)-ary search tree from
$\pi$ in the manner we have described in
Section~\ref{sec:intro_background}.  For example, for \( m=3 \) and \(
n=4 \), there are 6 distinct trees each occurring with probability \(
1/6 \), as depicted in Figure~\ref{fig:distribution}. We note that it
is not the case in general that all trees are equiprobable.  For
example, if $m=2$ and $n=3$, there are 5 distinct trees; 4 of these
have probability $1/6$ each, and 1 (corresponding to permutations
$213$ and $231$) has probability $2/6$.

Unlike under the uniform model, functionals of search trees under the random
permutation model have been studied a great deal, the most attention
being paid to binary search trees.  Indeed, much of the motivation for
studying additive functionals comes from the fact that the number of
key comparisons in Quicksort is such a functional and under the random
permutation model it corresponds to the cost of Quicksort on uniformly
random input.  There is an extensive literature on the analysis of
Quicksort; see, for
example,~\cite{MR1932675,FillJansonECP,FillJansonbook,MR2002k:68042,MR2000b:68051,MR97e:68018,MR90k:68131,hennequin91:_analy,MR90k:68132,MR56:7309,MR51:7356,MR92f:68028,01808498,MR96j:68049,janson153}.
Many researchers have studied various other functionals under the
random permutation model;
see~\cite{MR1896364,MR1910527,MR1910529,MR2002h:68035,MR1887300,MR2002f:68023,MR2002k:60088,MR2000c:68037,MR98m:68041,MR99e:68026,MR98m:68042,MR98m:68047,MR97c:68041,MR97f:68021,MR96d:60100,MR96d:60101,MR96k:68033,MR87i:68009,MR2002j:62020}
for treatment of functionals defined on binary search trees.  For~\(
m\)-ary search trees see~\cite{MR95i:68030,MR90a:68012,MR1871558,MR99f:68034}.
Chapter~3 of~\cite{MR93f:68045} is a rich source of results on the
properties of \( m \)-ary search trees.

\section{Roadmap}
\label{sec:roadmap}

In Chapter~\ref{cha:additive-type-tree} we will formally define the
problem of additive tree functionals, introducing the concept of
a toll function.  The recurrences obtained will be translated to
relations among generating functions.  The resulting relations will
serve as a motivation for our main theoretical results---\emph{transfer
theorems} and the extension of the \emph{singularity analysis} toolkit.  In
Chapter~\ref{cha:transfer-theorems} we consider  the
random permutation model and prove general results that relate the
asymptotics of the toll function to that of the corresponding tree
functional.  We will then proceed to extend singularity analysis in
Chapter~\ref{cha:sing-analys-hadam} so that it is readily applicable
to recurrences obtained under the uniform model.

Next we will consider various applications.  In
Chapter~\ref{cha:binary-search-trees} we will combine the use of the
extended singularity analysis toolkit and the transfer theorems
derived in Chapter~\ref{cha:transfer-theorems} to derive a complete
asymptotic expansion of the mean of the shape functional of an $m$-ary
search tree under the random permutation model and the variance for
binary search trees.  Then in Chapter~\ref{cha:limit-distr-m} we will
use the transfer theorems of Chapter~\ref{cha:transfer-theorems} to
determine limiting distributions of additive tree functionals
under the random permutation model.  Finally, in
Chapters~\ref{sec:toll-function-nalpha}
and~\ref{sec:catalan_shape-functional}, the extension of singularity
analysis to Hadamard products will be applied to the Catalan model
for binary search trees, resulting in limit laws under two different
classes of toll functions.



\chapter{Additive tree functionals}
\label{cha:additive-type-tree}

We now formally introduce additive functionals on~\(m\)-ary search
trees.  We will provide some motivating examples and derive the basic
recurrence that governs the behavior of these functions under each of
the random permutation model and the uniform model.

\section{Notation and definitions}
\label{sec:notation}
\begin{figure}[htbp]
\begin{center}
  \centering
  \includegraphics{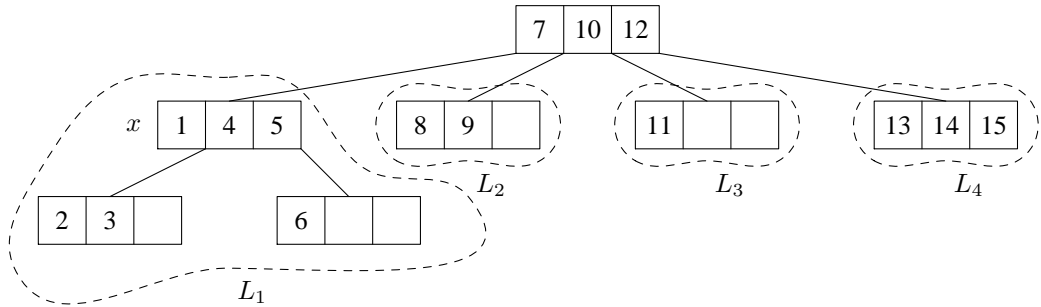}%
  \caption[Example of notation for~$m$-ary search trees.]%
  {Example of notation for the quaternary tree~\(T\) of
    Figures~\ref{fig:examplemAry} and~\ref{fig:buildmAry}. In this
    example,~$|T|=15$ and $T_x = L_1$.}
  \label{fig:terms}
\end{center}
\end{figure}

We first establish some notation.  Let~\(T\) be an~\(m\)-ary search
tree.  We use~\(|T|\) to denote the number of keys in~\(T\).  Call a
node \emph{full} if it contains~\(m-1\) keys.  For~\(1 \leq j \leq
m\), let~\(L_j(T)\) denote the~\(j\)th subtree pendant from the root
of~\(T\).  [The numbering of $L_1(T), \ldots, L_m(T)$ will correspond
to left-to-right ordering in our figures.] When it is
clear to which tree~\(T\) we are referring we will denote the subtrees
simply by~\( L_1,\ldots,L_m \).  For~\( x \) a node in~\( T \),
write~\( T_x \) for the subtree of~\( T \) consisting of $x$ and its
descendants, with $x$ as root.  This notation is illustrated in
Figure~\ref{fig:terms}.

\begin{definition}
  \label{def:additive-functional} Fix~\( m \geq 2 \).  We will call a
  functional~\( f \) on~\( m \)-ary search trees  an
  \emph{additive tree functional} if it satisfies the recurrence
  \begin{equation}
    \label{eq:2.1}
    f(T) = \sum_{i=1}^m f( L_i(T) ) + t_{|T|},
  \end{equation}
  for any tree~\( T \) with~\( |T| \geq m -1 \). Here~\( (t_n)_{n \geq
  m-1} \) is a given sequence, henceforth called the \emph{toll
  function}.
\end{definition}
Note that the recurrence~\eqref{eq:2.1} does not make any reference
to~\( t_n \) for~\( 0 \leq n \leq m-2 \) nor specify~\( f(T) \) for~\(
0 \leq |T| \leq m-2 \).

\section{Examples}
\label{sec:additive-examples}
Several interesting examples can be cast as additive functionals.
\begin{example}
  \label{example:space-requirement} If we specify~\( f(T) \)
  arbitrarily for~\( 0 \leq |T| \leq m-2 \) and take~\( t_n \equiv c
  \) for~\( n \geq m-1 \), we obtain the ``additive functional''
  framework of~\cite[\S3.1]{MR93f:68045}.  (Our definition of an
  additive functional generalizes this notion.)  In particular if we
  define~\( f(\emptyset) := 0 \) and~\( f(T) := 1 \) for the unique \(
  m \)-ary search tree~\( T \) on~\( n \) keys for \( 1 \leq n \leq
  m-2 \) and let \( t_n \equiv 1 \) for \( n \geq m-1 \), then \( f(T)
  \) counts the number of nodes in $T$ and thus gives the \emph{space
  requirement} functional discussed in~\citep[\S3.4]{MR93f:68045}.
\end{example}

\begin{example}
  \label{example:path-length}
  If we define \( f(T) := 0 \) when \( 0 \leq |T| \leq m-2 \) and
  \( t_n := n-(m-1) \) for \( n \geq m-1 \) then \( f \) is the
  \emph{internal path length} functional discussed
  in~\citep[\S3.5]{MR93f:68045}:\ $f(T)$ is the sum of all root-to-key
  distances in $T$.
\end{example}

\begin{example}
  \label{example:shape-functional}
   As noted in Chapter~\ref{cha:introduction}, we can associate with
  each permutation of \([n]\) an \( m \)-ary search tree.  Let \( Q \)
  denote the probability mass function of the distribution of trees
  induced by associating an \( m \)-ary search tree with each of the
  \(n!\) (equally likely)  permutations of \( [n] \).  Dobrow and
  Fill~\cite{MR97k:68038} noted that
  \begin{equation}
    \label{eq:2.2}
    Q(T) = \frac{1}{\prod_{x} \binom{|T_x|}{m-1}},
  \end{equation} where the product in~\eqref{eq:2.2} is over all full
  nodes in \( T \).  This functional is sometimes called the
  \emph{shape functional} as it serves as a crude measure of the
  ``shape'' of the tree, with ``full'' trees (like the complete tree)
  achieving larger values of $Q$.  For further discussions along these
  lines, consult~\cite{MR97k:68038} and~\cite{MR97f:68021}.  If we
  define \( f(T) := 0\) for \( 0 \leq |T| \leq m-2 \) and \( t_n :=
  \ln{\binom{n}{m-1}} \) for \( n \geq m-1 \), then \( f(T) =
  -\ln{Q(T)} \).  In this work we will refer to $-\ln{Q}$ (rather
  than $Q$) as the shape functional.
\end{example}

\section{A common framework for all moments}
\label{sec:additive-common-framework}

Our aim is to study the distribution of the functional \( f(T) \)
satisfying~\eqref{eq:2.1} under the two aforementioned probability
models.  Under the uniform model \( T \) is distributed uniformly over
all \( m \)-ary search trees with $n$ keys, while under the random
permutation model \( T \) has the distribution \(Q\) defined
at~\eqref{eq:2.2}.  In either case we will use \( X_n \) to denote a
random variable having the distribution of $f(T)$ and define \(
\mu_n(k) := \E{X_n^k} \) for \( k \geq 1\), the choice of model being
clear from context.  To address the problem recursively we examine the
joint distribution of the subtree sizes $|L_1|, \ldots, |L_m|$.

\subsection*{Random permutation model}
\label{sec:additive-rand-perm-model}
Under the random permutation model the joint distribution of the
subtree sizes $|L_1|, \ldots, |L_m|$ is uniform over all \(
\binom{n}{m-1} \) \( m \)-tuples  of
nonnegative integers that sum to \( n - (m-1)
\):\ see~\citep[Exercise~3.8]{MR93f:68045}.  We now apply the law of
total expectation by conditioning on $(|L_1|,\ldots,|L_m|)$.
Let \( \sum_{\mathbf{j}} \) denote the sum over
all \( m \)-tuples \(( j_1,\ldots,j_m) \) that sum to \( n-(m-1) \)
and \( \sum_{\mathbf{k}} \) the sum over all \( (m+1) \)-tuples \( (
k_1,\ldots,k_{m+1} ) \) of nonnegative integers that sum to \( k \).
Then, letting \( \oplus \) denote sums of mutually independent random
variables, for \( n \geq m-1 \) we have
\begin{align*}
  \mu_n(k) &= \E \E( X_n^k\; \bigl|\; |L_1|,\ldots,|L_m| ) =
  \frac{1}{\binom{n}{m-1}} \sum_{\mathbf{j}} \E{(
  X_{j_1} \oplus \cdots \oplus X_{j_m} + t_n )^k} \\
  &= \frac{1}{\binom{n}{m-1}} \sum_{\mathbf{j}}{\sum_{\mathbf{k}}}
  \binom{k}{k_1,\ldots,k_m,k_{m+1}} \mu_{j_1}(k_1)\cdots
  \mu_{j_m}(k_m) t_n^{k_{m+1}}.
\end{align*}
We can rewrite this equation as
\begin{equation}
  \label{eq:2.3}
  \mu_n(k) = \frac{m}{\binom{n}{m-1}} \sum_{j=0}^{n-(m-1)} \binom{n-1-j}{m-2}
  \mu_j(k) + r_n(k),
\end{equation}
where
\begin{equation}
  \label{eq:2.4}
  r_n(k) := \sumstark \binom{k}{k_1,\ldots,k_m,k_{m+1}}
  t_n^{k_{m+1}}  \frac{1}{\binom{n}{m-1}} \sum_{\mathbf{j}}
  \mu_{j_1}(k_1) \cdots \mu_{j_m}(k_m),
\end{equation}
with \( \sum_{\mathbf{k}}^{*} \) denoting the same sum as
\(\sum_{\mathbf{k}} \) with the additional restriction that \( k_i < k
\) for \( i=1,\ldots,m \).  We have thus established that the moments
$\mu_n(k)$ each satisfy the same basic recurrence in $n$, differing
as $k$ varies only in the non-homogeneous term \( r_n(k) \).  Observe
that \( r_n(1) = t_n \), the toll function.  We record this important
fact as
\begin{proposition}
  \label{thm:rand-perm-model}
  Under the random permutation model all moments of an additive
  functional satisfy the basic recurrence
  \begin{equation}
    \label{eq:2.7}
    a_n = b_n + \frac{m}{\binom{n}{m-1}} \sum_{j=0}^{n-(m-1)}
    \binom{n-1-j}{m-2} a_j, \qquad n \geq m-1,
  \end{equation} with specified initial conditions (say) \( a_j := b_j
  \), \( 0 \leq j \leq m-2 \).  [Recall the statement following
  Definition~\ref{def:additive-functional} about the initial
  conditions for the recurrence~\eqref{eq:2.1}.]
\end{proposition}

In order to analyze the recurrence relation~\eqref{eq:2.7} we introduce
generating functions
\begin{equation*}
  A(z) := \sum_{n=0}^\infty a_n z^n \quad \text{ and } \quad
  B(z) := \sum_{n=0}^\infty b_n z^n.
\end{equation*}
Furthermore, let
\begin{equation*}
  \rising{x}{r} := \prod_{k=0}^{r-1} (x+k)
\end{equation*}
denote the $r$th rising factorial power of $x$ and
\begin{equation*}
  \falling{x}{r} := \prod_{k=0}^{r-1} (x-k)
\end{equation*}
the $r$th falling factorial power of $x$.  Multiplying~\eqref{eq:2.7} by \(
\falling{n}{m-1} z^{n-(m-1)} \) and summing over \( n \geq m-1 \) we
get the differential equation
\begin{equation}
  \label{eq:2.8}
  A^{(m-1)}(z) = B^{(m-1)}(z) + m! (1-z)^{-(m-1)} A(z).
\end{equation}
Equations of the form~\eqref{eq:2.8} are members of a class known as
\emph{Euler differential equations}.  Using the method of variation of
parameters and combinatorial identities one can obtain the general
solution to this equation.  We only state the results here, referring
the interested reader to Appendix~\ref{appendix:diffeq}.
\begin{theorem}
  \label{thm:diffeq-soln}
  Let \( A \) and \( B \) denote the
  respective ordinary generating functions of the sequences \( (a_n)
  \) and \( (b_n) \) in Proposition~\ref{thm:rand-perm-model}.  Then
  \begin{equation}
    \label{eq:2.11}
    A(z) = \sum_{j=1}^{m-1} c_j (1-z)^{-\lambda_j} + \sum_{j=1}^{m-1}
    \frac{(1-z)^{-\lambda_j}}{\psi'(\lambda_j)} \int_{0}^z
    B^{(m-1)}(\zeta) (1-\zeta)^{\lambda_j+m-2} \; d\zeta,
  \end{equation}
  where \( \psi \) is the indicial polynomial
  \begin{equation}
    \label{eq:2.10}
    \psi(\lambda) := \rising{\lambda}{m-1} - m! = \lambda(\lambda+1)
    \cdots (\lambda+m-2) - m!
  \end{equation} with roots \( 2 =: \lambda_1,\lambda_2, \ldots,
  \lambda_{m-1} \) in nonincreasing order of real parts.
  In~\eqref{eq:2.11}, the
  coefficients  \( c_1, c_2, \ldots, c_{m-1} \)  can
  be written explicitly in terms of the initial conditions \( b_0,
  \ldots, b_{m-2} \) as
  \begin{equation}
    \label{eq:2.9}
    c_j = \frac{m!}{\psi'(\lambda_j)} \sum_{k=0}^{m-2} b_k
    \frac{k!}{\rising{\lambda_j}{k+1}}, \qquad j=1,\ldots,m-1.
  \end{equation}
\end{theorem}
\begin{proof}
  See Appendix~\ref{appendix:diffeq}, in particular
  Corollary~\ref{corollary:A.5.2} and Proposition~\ref{lemma:A.2.1}.
\end{proof}
The indicial polynomial~\eqref{eq:2.10} is well studied;
see~\cite{MR93f:68045,MR1871558,MR96j:68042} and
Appendix~\ref{appendix:indicial}.  We will exploit the
expression~\eqref{eq:2.11} for \( A(z) \) to relate the asymptotic
properties of the sequence \( (b_n) \) to those of \( (a_n) \) in
Chapter~\ref{cha:transfer-theorems} and then, in
Chapter~\ref{cha:limit-distr-m}, use the transfer theorems of
Chapter~\ref{cha:transfer-theorems} to derive limiting distributions
of additive functionals.

\subsection*{Uniform model for binary trees}
\label{sec:additive-uniform-model}

We now turn our attention to the uniform model for binary search
trees.  The techniques we develop in this dissertation can be extended to
handle the uniform model on $m$-ary search trees as well, although
much of that 
analysis will involve singularity analysis of implicitly defined
functions.  Our techniques can also be extended to handle
any simply generated family of
trees~\cite{MR80k:05043,MR86h:68070}.

It is well known that the number of binary trees on \( n
\) nodes is counted by the \( n \)th Catalan
number (see, e.g.,~\cite[Exercise~6.1.19c]{MR2000k:05026})
\begin{equation*}
  \beta_n := \frac{1}{n+1} \binom{2n}{n},
\end{equation*}
with generating function
\begin{equation}
  \label{eq:CAT}
  \CAT(z) := \sum_{n=0}^\infty \beta_n z^n = \frac{1}{2z} (1 -
  \sqrt{1-4z}).
\end{equation}
Let \( \mu_n(k) \) be defined as the \( k \)th moment of the random
variable \( X_n \) and \( \bar{\mu}_n(k) := \beta_n \mu_n(k)/4^n \).
By conditioning on the key stored at the root, noting that the
probability is $\beta_{l-1}\beta_{n-l}/\beta_n$ that the key stored at
the root is $l$, $1 \leq l \leq n$, and multiplying throughout by \(
\beta_n/4^n \) we get, for \( n \geq 1 \),
\begin{equation}
  \label{eq:2.5}
  \bar{\mu}_{n}(k) = \frac12 \sum_{j=1}^{n} \frac{\beta_{n-j}}{4^{n-j}}
  \bar{\mu}_{j-1}(k) +  \bar{r}_n(k),
\end{equation}
where
\begin{equation}
  \label{eq:2.6}
  \bar{r}_n(k) := \frac14 \sum_{\substack{k_1+k_2+k_3=k\\k_1,k_2<k}}
  \binom{k}{k_1,k_2,k_3} t_n^{k_3} \sum_{j=1}^{n} \bar{\mu}_{j-1}(k_1)
  \bar{\mu}_{n-j}(k_2).
\end{equation}
Again we observe that moments of all order satisfy the same basic
recurrence, with
\[
\bar{r}_n(1) = \frac{t_n}{4^n} \sum_{j=1}^n{\beta_{j-1}\beta_{n-j}} =
t_n \frac{\beta_n}{4^n}.
\]
We define generating
functions
\begin{equation}
  \label{eq:2.20}
  \overline{M}_{k}(z) := \sum_{n=0}^\infty \bar{\mu}_{n}(k) z^n, \qquad
  \overline{R}_{k}(z) := \sum_{n=0}^\infty \bar{r}_{n}(k) z^n, \qquad
  T(z) := \sum_{n=0}^\infty t_n z^n.
\end{equation}
Then, using~\eqref{eq:2.5} and~\eqref{eq:CAT}, we get
\begin{equation}
  \label{eq:2.21}
  \overline{M}_{k}(z) = \frac{ \overline{R}_k(z)}{ \sqrt{1-z} },
\end{equation}
so to analyze \( \overline{M}_{k}(z)  \) it is sufficient
to analyze \( \overline{R}_k(z) \).  With this in mind we introduce the
Hadamard product of two power series.
\begin{definition}
  \label{definition:hadamard}
  The \emph{Hadamard product} of two power series \( F \) and \( G
  \), denoted by \( F(z) \odot G(z) \) or \( (F \odot G)(z) \), is the
  power series defined by
  \[
  F(z) \odot G(z) := \sum_n f_n g_n z^n,
  \]
  where
  \[ F(z) = \sum_n f_n z^n \text{ and } G(z) = \sum_{n} g_n z^n. \]
\end{definition}
In light of this definition, the defining relation for \( \overline{R}_k(z)
\) at~\eqref{eq:2.6} leads to the following relation in terms of
generating functions:
\begin{equation}
  \label{eq:2.22}
  \overline{R}_k(z) =
  \sum_{\substack{k_1+k_2+k_3=k\\k_1,k_2<k}}  \binom{k}{k_1,k_2,k_3}
  ( T(z) )^{\odot k_3} \odot [ \frac{z}{4} \overline{M}_{k_1}(z)
  \overline{M}_{k_2}(z)  ],
\end{equation}
where, for \( k \) a nonnegative integer,
\begin{equation*}
  T(z)^{\odot k} := \underbrace{ T(z) \odot \cdots \odot T(z) }_{k}.
\end{equation*}
Our approach for getting asymptotic expansions for moments of all
orders relies on \emph{singularity analysis} of generating
functions~\cite{MR90m:05012}, a tool that, under suitable conditions,
allows the translation of asymptotics of generating functions to
asymptotics of the underlying sequences.  In order to analyze
$\overline{M}_k(z)$ and \( \overline{R}_k(z) \) inductively (in $k$),
we would like to be able to predict how singularities are composed
under Hadamard products.  This is precisely the problem studied in
Chapter~\ref{cha:sing-analys-hadam}.  We will, in
Chapters~\ref{sec:toll-function-nalpha}
and~\ref{sec:catalan_shape-functional}, use the results derived in
Chapter~\ref{cha:sing-analys-hadam}
to get limiting distributions induced by two types of toll functions.


\part{Theory}
\chapter{Random permutation model: transfer theorems}
\label{cha:transfer-theorems}

In Proposition~\ref{thm:rand-perm-model} we have seen that under the
random permutation model,
moments of each specified order for any additive functional on
\(m\)-ary search trees satisfy a recurrence relation of the form
\begin{equation}
  \label{eq:3.1.1}
  a_n = b_n + \frac{m}{\binom{n}{m-1}} \sum_{j=0}^{n-(m-1)}
  \binom{n-1-j}{m-2} a_j, \qquad n \geq m-1,
\end{equation}
with specified initial conditions \( a_j := b_j \) for \(
j=0,\ldots,m-2 \).  By Theorem~\ref{thm:diffeq-soln} the
output generating function
\begin{equation*}
  A(z) := \sum_{n=0}^{\infty} a_n z^n
\end{equation*}
is given in terms of the input generating function
\begin{equation*}
  B(z) := \sum_{n=0}^\infty b_n z^n
\end{equation*}
by
\begin{equation*}
  A(z) = \sum_{j=1}^{m-1} c_j (1-z)^{-\lambda_j} + \sum_{j=1}^{m-1}
  \frac{(1-z)^{-\lambda_j}}{\psi'(\lambda_j)} \int_{0}^z
  B^{(m-1)}(\zeta) (1-\zeta)^{\lambda_j+m-2} \; d\zeta,
\end{equation*}
where \( \psi \) is the indicial polynomial
\begin{equation*}
  \psi(\lambda) := \rising{\lambda}{m-1} - m!
\end{equation*}
with roots \( 2 =: \lambda_1,\lambda_2,\ldots,\lambda_{m-1} \) listed
in nonincreasing order of real parts.  We refer the reader to
Appendix~\ref{appendix:indicial} for relevant facts about
this polynomial, and recall from there the notation
\begin{equation*}
  \alpha := \max_{2 \leq j \leq m-1} \Re{(\lambda_j)}.
\end{equation*}
The coefficients \( c_1, c_2, \ldots, c_{m-1} \) can
be written explicitly in terms of the initial conditions \( b_0,
\ldots, b_{m-2} \).  In particular, by Proposition~\ref{lemma:A.2.1}
with $j=1$ and $g_{\text{c}}^{(k)}(0) = k! b_k$ there, we get [also
using~\eqref{eq:psp2}]
\begin{equation}
  \label{eq:3.1.c_1}
  c_1 = \frac{1}{H_m-1} \sum_{j=0}^{m-2} \frac{b_j}{(j+1)(j+2)}.
\end{equation}
Here \( H_m \) is the $m$th Harmonic number:\ $H_m = \sum_{j=1}^m
\tfrac1j$. 

In this chapter we will obtain results of both exact and asymptotic
natures that link the behavior of the output sequence $(a_n)$ to that
of the input sequence $(b_n)$ in~\eqref{eq:3.1.1}.

\emph{Note.}  In the sequel whenever multi-valued functions are
encountered, we will assume that we follow the principal branch.

\section{The Exact Transfer Theorem}
\label{sec:exact-transf-theor}

Our main result repeats~\eqref{eq:2.11} and then gives an alternative
exact transfer from \( B \) to \( A \):
\begin{theorem}[Exact Transfer Theorem (ETT)]
  \label{theorem:ETT}
  Let
  \begin{equation}
    \label{eq:3.1.1.4}
    \widehat{B}(z) := B(z) - \sum_{j=0}^{m-2} b_j z^j =
    \sum_{n=m-1}^\infty b_n z^n.
  \end{equation}
  Then
  \begin{align}
    \label{eq:3.1.1.5}
    A(z) &= \sum_{j=1}^{m-1} c_j (1-z)^{-\lambda_j} + \sum_{j=1}^{m-1}
    \frac{(1-z)^{-\lambda_j}}{\psi'(\lambda_j)} \int_{0}^z
    B^{(m-1)}(\zeta) (1-\zeta)^{\lambda_j+m-2} \; d\zeta\\
    \label{eq:3.1.1.6}
    &= \sum_{j=1}^{m-1} c_j (1-z)^{-\lambda_j} + \widehat{B}(z) + m!
    \sum_{j=1}^{m-1} \frac{(1-z)^{-\lambda_j}}{\psi'(\lambda_j)}
    \int_{0}^z \widehat{B}(\zeta)(1-\zeta)^{\lambda_j-1}\,d\zeta.
  \end{align}
\end{theorem}
\begin{proof}
  We begin with~\eqref{eq:3.1.1.5}, repeated from
  Theorem~\ref{thm:diffeq-soln}, and note that \( B \) there can be
  replaced by \( \widehat{B} \).  We then use repeated integration by
  parts and invoke Identity~\ref{identity:B.2.2}.  Denoting
  \begin{equation*}
    \widehat{A}(z) := A(z) - \sum_{j=1}^{m-1} c_j (1-z)^{-\lambda_j},
  \end{equation*}
  after \( m-2 \) integrations by parts we find
  \begin{equation*}
    \widehat{A}(z) = \sum_{j=1}^{m-1}
    \frac{(1-z)^{-\lambda_j}}{\psi'(\lambda_j)} (\lambda_j+m-2) \cdots
    (\lambda_j+1) \int_{0}^z \widehat{B}'(\zeta)
    (1-\zeta)^{\lambda_j} \, d\zeta.
  \end{equation*}
  But 
  \begin{equation*}
    (\lambda_j+m-2) \cdots  (\lambda_j+1) =
    \frac{\rising{\lambda_j}{m-1}}{\lambda_j} = \frac{\psi(\lambda_j)
      + m!}{\lambda_j} = \frac{m!}{\lambda_j},
  \end{equation*}
  so
  \begin{equation*}
    \widehat{A}(z) = m! \sum_{j=1}^{m-1}
    \frac{(1-z)^{-\lambda_j}}{\lambda_j \psi'(\lambda_j)}
    \int_{0}^z \widehat{B}'(\zeta) (1-\zeta)^{\lambda_j}\,d\zeta.
  \end{equation*}
  We obtain~\eqref{eq:3.1.1.6} by performing one more integration by
  parts and utilizing Identity~\ref{identity:B.1} with \( \lambda = 0
  \).
\end{proof}
\begin{remark}
  \label{remark:ettterm}
For computations, \eqref{eq:3.1.1.5} might be easier to use when it is
no bother to compute derivatives of \( B \); otherwise,
\eqref{eq:3.1.1.6} is easier.  Equation~\eqref{eq:3.1.1.6} will be
used in establishing part~(a) of the Asymptotic Transfer Theorem in
Section~\ref{sec:asympt-transf-theor}; the proof of part~(b) will
use~\eqref{eq:3.1.1.5}.
\end{remark}

\section{The Asymptotic Transfer Theorem}
\label{sec:asympt-transf-theor}

It is quite easy to transfer asymptotics for \( B \) to asymptotics
for \( A \) using the ETT.  We give three examples important for
applications to moments of additive functionals in the next
theorem, proved independently in~\cite{MR1871558} using a quite
different approach.
\begin{theorem}[Asymptotic Transfer Theorem (ATT)]
  \label{theorem:att}
  Suppose that \( (a_n) \) and \( (b_n) \) satisfy the
  recurrence~\eqref{eq:3.1.1}.
  \begin{enumerate}[(a)]
  \item If
    \begin{equation}
      \label{eq:3.2.1.7}
      b_n = o(n) \quad \text{ and } \quad \sum_{n=0}^\infty
      \frac{b_n}{(n+1)(n+2)} \text{ converges (conditionally),}
    \end{equation}
    then
    \begin{equation}
      \label{eq:3.2.1.8}
      a_n = \frac{K_1}{H_m-1} n + o(n), \quad\text{ where }\quad
      K_1 := \sum_{j=0}^\infty \frac{b_j}{(j+1)(j+2)}.
    \end{equation}
  \item If \( b_n \equiv K_2(n+1) + h_n \) where \( h_n \)
    satisfies~\eqref{eq:3.2.1.7} [with \( (b_n) \) replaced by \( (h_n)
    \)], then
    \begin{equation}
      \label{eq:3.2.1.9}
      a_n = \frac{K_2}{H_m-1} n H_n + \frac{K_3}{H_m-1} n + o(n),
    \end{equation}
    where
    \begin{equation}
      \label{eq:3.2.1.10}
      K_3 := \sum_{j=0}^\infty \frac{h_j}{(j+1)(j+2)} + K_2 \left[
        \frac{H_m-1}{2} - 1 + \frac{H_m^{(2)} - 1}{2(H_m-1)} \right].
    \end{equation}
  \item If \( b_n = K_4 n^v + o( n^{\Re(v)} ) \) with \( \Re(v) > 1 \),
    then
    \begin{equation}
      \label{eq:3.2.1.11}
      a_n = \frac{K_4}{1- \frac{m! \Gamma(v+1)}{\Gamma(v+m)}} n^v + o(
      n^{\Re(v)} ).
    \end{equation}
  \end{enumerate}
\end{theorem}
Before giving the proof we establish a few  lemmas.  The
first result is found in the proof of Lemma~6 of~\cite{MR1871558}.  For
completeness we include a proof.  Here, as throughout, \( [z^n] F(z)
\) is the coefficient of \( z^n \) in the series expansion of \(F\).
\begin{lemma}
  \label{lemma:3.2.4.1}
  Let \( Y(z) \) denote the ordinary generating function of the
  sequence \( (y_n) \) with \( y_0 = 0\).  For any \( \lambda \in
  \mathbb{C} \),
  \begin{equation}
    \label{eq:3.2.4.1}
    [z^n] \left( (1-z)^{-\lambda} \int_{0}^z
    (1-\zeta)^{\lambda-1} Y(\zeta)\,d\zeta \right) = \sum_{k=0}^{n-1}
    \frac{y_k}{k+1} \prod_{j=k+2}^n \left( 1 + \frac{\lambda-1}{j}
    \right), \quad n \geq 0.
  \end{equation}
\end{lemma}
\begin{proof}
  The function
  \begin{equation*}
    W(z) := (1-z)^{-\lambda} \int_{0}^z (1-\zeta)^{\lambda-1}
    Y(\zeta) \,d\zeta
  \end{equation*}
  is the unique solution with \( W(0) = 0 \) to the differential
  equation
  \begin{equation*}
    W'(z) = \lambda(1-z)^{-1}W(z) + (1-z)^{-1} Y(z):
  \end{equation*}
  that is, \( w_n := [z^n] W(z) \) for \( n \geq 0 \) satisfies \(
  w_0=0 \) and
  \begin{equation}
    \label{eq:3.2.4.4}
    w_n = \frac{\lambda}{n} \sum_{k=0}^{n-1} w_k + \frac{1}{n}
    \sum_{k=0}^{n-1} y_k, \qquad n \geq 1.
  \end{equation}
  To solve the recurrence~\eqref{eq:3.2.4.4} we compute \( n w_n -
  (n-1)w_{n-1} = \lambda w_{n-1} + y_{n-1} \) so that
  \begin{equation*}
    w_n  = \left(1+\frac{\lambda-1}{n}\right) w_{n-1} + \frac1n
    y_{n-1}, \qquad n \geq 1.
  \end{equation*}
  Now we iterate to get the result.
\end{proof}
The product in~\eqref{eq:3.2.4.1} may be written (when \( \lambda \in
\mathbb{C} \setminus \{0,-1,-2,\ldots\} \)) as
\begin{equation*}
  \frac{\Gamma(\lambda+n)\Gamma(2+k)}{\Gamma(1+n)\Gamma(\lambda+k+1)},
\end{equation*}
which by Stirling's formula equals
\begin{equation}
  \label{eq:3.2.4.2}
  \frac{ n^{\lambda-1} [1 +
  O(n^{-1})]}{(k+1)^{\lambda-1}[1+O((k+1)^{-1})]}
\end{equation}
for \( n \geq 1 \) and \( k \geq 1 \).  [The product
in~\eqref{eq:3.2.4.1} equals~\eqref{eq:3.2.4.2} even if \( \lambda \in
\{0,-1,-2,\ldots\} \). ]  Also, of special interest is the case \(
\lambda=2 \), in which case~\eqref{eq:3.2.4.1} reduces to
\begin{equation}
  \label{eq:3.2.4.3}
  [z^n] \left( (1-z)^{-2} \int_{0}^z (1-\zeta) Y(\zeta)\,d\zeta
  \right) = (n+1) \sum_{k=0}^{n-1} \frac{y_k}{(k+1)(k+2)}, \qquad n
  \geq 0.
\end{equation}
For part~(a) of the \textsc{ATT} we use the following estimates
from~\cite{MR1871558}.  These follow readily from
Lemma~\ref{lemma:3.2.4.1}, \eqref{eq:3.2.4.2}, and~\eqref{eq:3.2.4.3}.
\begin{lemma}
  \label{lemma:3.2.4.2}
  \begin{enumerate}[(a)]
  \item  If \( \Re(\lambda) < 2 \) and
    \[
    Y(z) = \sum_{n=0}^\infty  y_n z^n
    \]
    satisfies \( y_0 = 0 \) and \( y_n = o(n) \), then
    \begin{equation*}
      [z^n] \left( (1-z)^{-\lambda} \int_{0}^z Y(\zeta)
        (1-\zeta)^{\lambda-1} \,d\zeta \right) = o(n).
    \end{equation*}
  \item With \( \widehat{B} \) defined at~\eqref{eq:3.1.1.4},
      if~\eqref{eq:3.2.1.7} holds, then
      \begin{equation*}
        [z^n] \left( (1-z)^{-2} \int_{0}^z \widehat{B}(\zeta)
        (1-\zeta) \, d\zeta \right) = n \sum_{j=m-1}^\infty
        \frac{b_j}{(j+1)(j+2)} + o(n).
      \end{equation*}
  \end{enumerate}
\end{lemma}

For part~(c) of the ATT, we will need the following simple comparison
lemma.
\begin{lemma}
  \label{lemma:3.2.4.3}
  If \( (b_n) \) and \( (b_n') \) are two toll sequences such that
  \begin{equation*}
    |b_n| \leq b_n' \quad\text{ for all \( n \geq 0 \)},
  \end{equation*}
  then the corresponding output sequences \( (a_n) \) and \( (a_n') \)
  in~\eqref{eq:3.1.1} (with the initial conditions stated there) satisfy
  \begin{equation*}
    |a_n| \leq a_n' \quad\text{ for all \( n \geq 0 \)}.
  \end{equation*}
\end{lemma}
\begin{proof}
  This follows immediately by induction.
\end{proof}

\begin{proof}[Proof of Theorem~\ref{theorem:att}]
  Part~(a) follows immediately from~\eqref{eq:3.1.1.6},
  \eqref{eq:3.2.1.7}, Lemma~\ref{lemma:3.2.4.2},
  \eqref{eq:3.1.c_1}, and~\eqref{eq:psp2}.

  For part~(b) suppose first that \( b_n \equiv n+1 \).  Then \( B(z)
  = (1-z)^{-2} \), so
  \begin{equation*}
    B^{(m-1)}(z) = m! (1-z)^{-(m+1)}.
  \end{equation*}
  Plugging this into~\eqref{eq:3.1.1.5} we find
  \begin{equation*}
    a_n = (n+1)\left[ c_1 + m! \sum_{j=2}^{m-1}
    \frac1{(2-\lambda_j)\psi'(\lambda_j)} \right] +
    \frac{m!}{\psi'(2)} [z^n] [ (1-z)^{-2} \ln{((1-z)^{-1})} ] + o(n).
  \end{equation*} Recall Identity~\ref{identity:B.5.22}; and since \(
  b_n \equiv n+1 \), we have in this case by~\eqref{eq:3.1.c_1}, that
  \( c_1 = 1 \).  Therefore
  \begin{align*}
    a_n &= (n+1)\left( 1 + \frac12 \left[ \frac{H_m^{(2)}-1}{(H_m-1)^2}
        - 1 \right] \right) + \frac1{H_m-1}[ (n+1) H_n - n ] + o(n)\\
      &= \frac1{H_m-1} n H_n + \left[ \frac12 - \frac1{H_m-1} +
    \frac{H_m^{(2)}-1}{2(H_m-1)^2} \right]n + o(n).
  \end{align*}
  This completes the proof of part~(b) for our special case, and the
  general case follows from this and part~(a) using the superposition
  principle.

  Finally for part~(c) suppose first that \( b_n \equiv
  \rising{(v+1)}{n}/n! \sim n^v/\Gamma(v+1) \), so that \( B(z) =
  (1-z)^{-(v+1)} \) and \( B^{(m-1)}(z) = \rising{(v+1)}{m-1}
  (1-z)^{-(v+m)} \).  Plugging this into~\eqref{eq:3.1.1.5} and
  utilizing Identity~\ref{identity:B.1} with \( \lambda=v+1 \) and the
  calculation
  \begin{equation*}
    \frac{\rising{(v+1)}{m-1}}{\rising{(v+1)}{m-1}-m!} = \left[ 1 -
    \frac{m! \Gamma(v+1)}{\Gamma(v+m)} \right]^{-1},
  \end{equation*}
  we find
  \begin{equation*}
    A(z) = \left[ 1 - \frac{m! \Gamma(v+1)}{\Gamma(v+m)} \right]^{-1}
    (1-z)^{-(v+1)} + O(|1-z|^{-2}).
  \end{equation*}
  By singularity analysis (see Theorem~\ref{theorem:transfer}), this
  completes the proof of part~(c) for our special case.

  To complete the proof (by superposition) in the general case, we
  need only show that if \( b_n = o(n^v) \) for real \( v > 1 \), then
  \( a_n = o(n^v) \).  Indeed, fix \( \epsilon > 0 \); then there
  exists a sequence \( (b_n') \) such that \( |b_n| \leq b_n' \) for
  all \( n \) and \( b_n' = \epsilon \rising{(v+1)}{n}/n! \) for all
  large \( n \).  This toll sequence is but a slight modification of
  our special-case toll sequence, and we see that its output
  generating function has coefficients
  \begin{equation*}
    a_n' = \epsilon'n^v + o(n^v), \quad \text{ where } \quad \epsilon'
    := \frac{\epsilon}{\Gamma(v+1)}\left[ 1 -
      \frac{m!\Gamma(v+1)}{\Gamma(v+m)} \right]^{-1}.
  \end{equation*}
  Now Lemma~\ref{lemma:3.2.4.3} implies that
  \begin{equation*}
    {\lim\sup}_n |a_n| n^{-v} \leq \epsilon';
  \end{equation*}
  since \( \epsilon \) (and hence \( \epsilon' \)) can be made
  arbitrarily small, this completes the proof.
\end{proof}

The conditions~\eqref{eq:3.2.1.7} on \( b_n \) are not only sufficient
but also necessary for asymptotic linearity of \( a_n \).  Indeed,
here is a converse:
\begin{proposition}
  \label{proposition:3.2.4.4}
  If \( a_n = Kn + o(n) \) for some constant \( K \),
  then~\eqref{eq:3.2.1.7} holds.
\end{proposition}
\begin{proof}
  Proposition~7 of~\cite{MR1871558} provides the simple proof that \(
  b_n = o(n) \).  Moreover, from~\eqref{eq:3.1.1.6}, the value
  of \( c_1 \) at~\eqref{eq:3.1.c_1}, and~\eqref{eq:3.2.4.3} with \( Y
  \) taken to be \( \widehat{B} \), we find that
  \begin{equation*}
    \sum_{j=0}^{n-1} \frac{b_j}{(j+1)(j+2)} = K(H_m-1) + o(1),
  \end{equation*}
  i.e., the series
  \begin{equation*}
    \sum_{n=0}^\infty \frac{b_n}{(n+1)(n+2)}
  \end{equation*}
   converges [conditionally, to \( K(H_m-1) \)].
 \end{proof}

 The following additional asymptotic transfer results are established by
 calculations similar to those in the proof of the ATT.
 \begin{definition}
   \label{definition:slowly-varying}
   A nonnegative function \( L(n) \) defined for $n \geq n_0 \geq 0$
   and not identically zero is called \emph{slowly varying} if
   \begin{equation*}
     \frac{L(n)}{L(\lfloor tn \rfloor)} \to 1 \quad\text{ as $n \to \infty$}
   \end{equation*}
   for every $t > 0$.
   If \( n_0 > 0 \), we define \( L(n) = 0 \) for \( 0 \leq n < n_0 \).
 \end{definition}
Examples of  slowly varying functions include any powers of $\ln{n}$
or \( \ln\ln{n} \).  The key property of slowly varying
functions that we employ is the regular variation
 fact~\cite[1.5.9a]{MR90i:26003}, that
 \begin{equation}
   \label{eq:3.2.4.10}
   L(n) = o\left( \sum_{k \leq n} \frac{L(k)}{k} \right).
 \end{equation}
 \begin{theorem}
   \label{theorem:3.2.4.5}
   Consider the initial value problem~\eqref{eq:2.7}.
     \begin{enumerate}[(a)]
     \item If \( 2 \leq m \leq 26 \) and \( b_n = o(\sqrt{n}) \), then
       we can refine~\eqref{eq:3.2.1.8} to
       \begin{equation}
         \label{eq:3.2.4.5}
         a_n = \frac{K_1}{H_m-1} n + o(\sqrt{n}).
       \end{equation}
     \item If \( \Re(\lambda_2) =: \alpha < 1 + \beta\) and \( b_n \sim
       n^\beta L(n) \)
       with \( L \) slowly varying, then we can
       refine~\eqref{eq:3.2.1.8} to
       \begin{equation}
         \label{eq:3.2.4.6}
         a_n = \frac{K_1}{H_m-1}n =
         \frac{\rising{(1+\beta)}{m-1}}{m!-\rising{(1+\beta)}{m-1}}
         {n}^\beta L(n) + o( {n}^\beta L(n) ).
       \end{equation}
     \item If \( b_n \sim n L(n) \) with \( L \) slowly varying, then,
       with \(K_1\) defined at~\eqref{eq:3.2.1.8},
       \begin{equation}
         \label{eq:3.2.4.7}
         a_n \sim
         \begin{cases}
           \displaystyle \frac{K_1}{H_m-1}n, & \text{ if }
           \displaystyle\sum^\infty \frac{L(k)}{k} < \infty,\\
           \displaystyle\frac1{H_m-1} n \displaystyle \sum_{k \leq n}
           \frac{L(k)}k, & \text{ if } \displaystyle\sum^\infty
           \frac{L(k)}{k} = \infty.
         \end{cases}
       \end{equation}
     \item Part~(c) of the ATT can be extended as follows.  If \( b_n
       = K_4 n^v L(n) + o(n^{\Re(v)L(n)}) \) with \( \Re(v) > 1 \) and
       \( L \) slowly varying, then
       \begin{equation}
         \label{eq:3.2.4.8}
         a_n = \frac{K_4}{1-\frac{m!\Gamma(v+1)}{\Gamma(v+m)}} n^v L(n) +
         o(n^{\Re(v)} L(n)).
       \end{equation}
     \end{enumerate}
 \end{theorem}
 \begin{proof}[Proof sketch]
   Whenever the conditions~\eqref{eq:3.2.1.7} are met we have by the
   ETT and~\eqref{eq:3.2.4.3}
   \begin{multline}
     \label{eq:3.2.4.9}
     a_n - \frac{K_1}{H_m-1}(n+1) = O(n^{\alpha-1}) + b_n -
     \frac1{H_m-1} (n+1) \sum_{k=n}^\infty
     \frac{\hat{b}_k}{(k+1)(k+2)}\\
     + m! \sum_{j=2}^{m-1}
     \frac{1}{\psi'(\lambda_j)} [z^n] \left( (1-z)^{-\lambda_j}
       \int_{0}^z \widehat{B}(\zeta)
       (1-\zeta)^{\lambda_j-1}\,d\zeta\right),
   \end{multline}
   where \( \alpha := \Re(\lambda_2) \) is the real
   part of the root of $\psi$ with second largest real part, which is
   smaller than \( 3/2 \) when \( m \leq 26 \) (see
   Appendix~\ref{appendix:indicial}).  Simple estimates, including the
   use of~\eqref{eq:3.2.4.2}, give cases~(a) and~(b); for~(b), the
   coefficient of \( {n}^\beta L(n) \) in~\eqref{eq:3.2.4.6} indeed is,
   using Identity~\ref{identity:B.1} and~\eqref{eq:psp2},
   \begin{equation*}
     1 - \frac1{(1-\beta)H_m-1} + m! \sum_{j=2}^{m-1}
     \frac{1}{((1+\beta)-\lambda_j)\psi'(\lambda_j)} =
     -\frac{\rising{(1+\beta)}{m-1}}{m! - \rising{(1+\beta)}{m-1}}.
   \end{equation*}

   In case~(c), from the ETT result~\eqref{eq:3.1.1.6} and simple
   estimates we find
   \begin{equation*}
     a_n = (n+1)\frac{1}{H_m-1} \sum_{k=0}^{m-1}
     \frac{b_k}{(k+1)(k+2)} + O( nL(n) ).
   \end{equation*}
   To finish, we use the regular variation fact~\eqref{eq:3.2.4.10}.

   In case~(d), again from~\eqref{eq:3.1.1.6} and simple estimates we
   find
   \begin{equation*}
     a_n = O(n) + b_n + K_4 n^v L(n) m! \sum_{j=1}^{m-1}
     \frac{1}{\psi'(\lambda_j)(v+1-\lambda_j)} + o( n^v L(n) ).
   \end{equation*}
   The proof is completed by using Identity~\ref{identity:B.1}.
 \end{proof}


\chapter{Singularity analysis of Hadamard products}
\label{cha:sing-analys-hadam}

\section{Background}
\label{sec:sa-background}
\begin{figure}[htbp]
\centering
\includegraphics{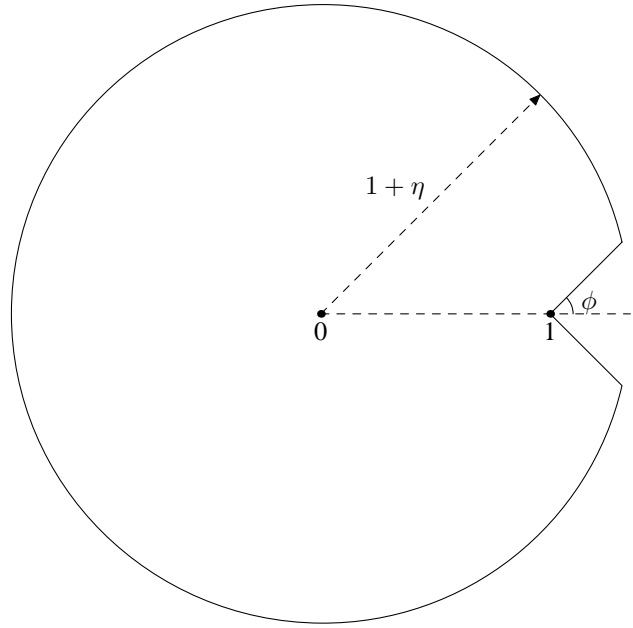}
\caption{The region $\Delta(\phi,\eta)$ in
  Definition~\ref{definition:delta-regular}.}
\label{fig:delta}
\end{figure}

Singularity analysis is ``a class of methods by which one can
translate, on a term-by-term basis, an asymptotic expansion of a
function around a dominant singularity into a corresponding asymptotic
expansion for the Taylor coefficients of the
function''~\cite{MR90m:05012}.  In general, it provides sufficient
conditions for the validity of the implications
\begin{align*}
  F(z) = O(G(z)) & \implies f_n = O(g_n)\\
  F(z) = o(G(z)) & \implies f_n = o(g_n)\\
  F(z) \sim G(z) & \implies f_n \sim g_n,
\end{align*}
where \( F \) and \( G \) are the ordinary generating functions of the
sequences \( (f_n) \) and \( (g_n) \) respectively.  The sufficiency 
conditions are expressed in terms of the following analytic
continuation condition.
\begin{definition}
  \label{definition:delta-regular}
  A function defined by a Taylor series with radius of convergence
  equal to~1 is said to be \emph{\(\Delta\)-regular} if, with the sole
  exception of \( z=1 \), it can be analytically continued in a domain
  \begin{equation}
    \label{eq:4.1}
    \Delta(\phi,\eta) := \{ z: |z| \leq 1+\eta,\; |\arg{(z-1)}| \geq
    \phi \},
  \end{equation}
  for some \( \eta > 0 \) and \( 0 < \phi < \pi/2 \).

  A function \( f \) is said to admit a \emph{\( \Delta \)-singular
    expansion} at \( z=1 \) if it is \( \Delta \)-regular and
  \begin{equation}
    \label{eq:4.2}
    f(z) = \sum_{n=0}^N c_n (1-z)^{\alpha_n} + O(|1-z|^A)
  \end{equation}
  as \( z \to 1 \) in \( \Delta(\phi,\eta) \), for a sequence of
  complex numbers \( (c_n)_{n=0}^N \) and a sequence of real numbers
  \( (\alpha_n)_{n=0}^N \) satisfying \( \alpha_n < A \).
\end{definition}
Note that under these assumptions the \( O(\cdot) \) error term
in~\eqref{eq:4.2} holds uniformly in \( z \in \Delta(\phi,\eta) \).
We also allow in the usual way infinite asymptotic expansions
\begin{equation*}
  f(z) \sim \sum_{n \geq 0} c_n (1-z)^{\alpha_n}
\end{equation*}
to represent an infinite collection of expansions of
type~\eqref{eq:4.2}.  When it is convenient, we will also allow powers
of \(\ln{((1-z)^{-1})}\) in singular expansions.  The region of
analytic continuation in the complex plane is depicted in
Figure~\ref{fig:delta}.  We will sometimes refer to such a region as
an \emph{indented crown}.

For completeness we state relevant results
from~\cite{MR90m:05012} and~\cite{MR2000a:05015}.  These will be used
frequently from here on.
\begin{theorem}[Theorem~1 in~\cite{MR90m:05012}]
  \label{theorem:transfer}
  Let \( f(z) \) be a \( \Delta \)-regular function in \(
  \Delta(\phi,\eta) \) and assume that, as \( z \to 1 \) in
  \(\Delta(\phi,\eta)\),
  \begin{equation*}
    f(z) = O(|1-z|^{-\alpha})
  \end{equation*}
  for some real number \( \alpha \).  Then the \(n\)th Taylor
  coefficient of \( f(z) \) satisfies
  \begin{equation*}
    f_n = [z^n] f(z) = O( n^{\alpha-1} ).
  \end{equation*}
\end{theorem}
Flajolet also showed~\cite{MR2000a:05015} that singularity analysis
is likewise applicable to the \emph{generalized polylogarithm}
$\Li_{\alpha,r}(z)$ defined for a
nonnegative integer \(r\) and  an arbitrary complex
number $\alpha$ and $|z| < 1$ as
\begin{equation}
  \label{eq:4.1.50}
\Li_{\alpha,r}(z) := \sum_{n=1}^\infty \frac{(\ln{n})^r}{n^\alpha} z^n.
\end{equation}
\begin{theorem}[Theorem~1 in~\cite{MR2000a:05015}]
  \label{theorem:Li-sing}
  The function \( \Li_{\alpha,0}(z) \) is analytic in the slit plane
  \( \mathbb{C} \setminus [1, +\infty) \).  For \( \alpha \notin
  \{ 1,2,\ldots \} \) the function \( \Li_{\alpha,0}(z) \) satisfies the
  singular expansion
  \begin{equation*}
    \Li_{\alpha,0}(z) \sim \Gamma(1-\alpha) t^{\alpha-1} + \sum_{j
      \geq 0} \frac{(-1)^j}{j!}\zeta(\alpha-j) t^j,
  \end{equation*}
  with
  \begin{equation*}
    t = -\ln{z} = \sum_{l=1}^\infty \frac{(1-z)^{l}}{l}
  \end{equation*}
  as \( z \to 1 \) in the sector \( -\pi + \epsilon < \arg{(1-z)} <
  \pi-\epsilon \) for any \( \epsilon > 0 \).

  For \( r = 1,2,\ldots \), the singular expansion of \(
  \Li_{\alpha,r}(z) \) is   obtained by
  \begin{equation*}
    \Li_{\alpha,r}(z) = (-1)^r \frac{\partial^r}{\partial{\alpha^r}}
    \Li_{\alpha,0}(z).
  \end{equation*}
\end{theorem}
The case when \( \alpha = m \in \{ 1, 2, \ldots \} \) is also handled
in~\cite[p.~380]{MR2000a:05015}:
\begin{equation}
  \label{eq:4.1.posintalpha}
  \Li_{m,r}(z) = \res{((-1)^r
  \zeta^{(r)}(s+m)\Gamma(s)t^{-s})_{s=1-m}} + \sum_{j \geq 0, j \ne
  m-1} \frac{(-1)^{j+r}}{j!} \zeta^{(r)}(m-j) t^j
\end{equation}
with \(t = -\ln{z}\), where \( \res(\cdot) \) denotes a residue.

Singularity analysis is a very useful tool in the analysis of
combinatorial structures mainly because generating functions of
structures that are constructed combinatorially readily admit the
technical conditions needed for its application.  In such cases
precise asymptotic expansions of the coefficients can be derived.  In
the next two sections we will extend singularity analysis so that it is
applicable to the right-hand side of~\eqref{eq:2.22}.

\section{Singular expansions and differentiation}
\label{sec:sing-expans-diff}

In preparation for the treatment of Hadamard products in
Section~\ref{sec:hadamard-products}, we need a
theorem that enables us to differentiate local expansions of analytic
functions around a singularity. Such a result cannot of course be
unconditionally true; see, for example,~\eqref{eq:4.2.56}. However, it
turns out that the class of functions amenable to singularity analysis
satisfy this property.

\begin{theorem}
  \label{theorem:2A}
  If $f(z)$ is $\Delta$-regular and admits a $\Delta$-expansion near
  its singularity in the sense of
  Definition~\ref{definition:delta-regular}, then for each integer $k
  > 0$, $\tfrac{d^k}{dz^k}f(z)$ is also $\Delta$-regular and admits a
  $\Delta$-expansion obtained through term-by-term
  differentiation:
  \begin{equation*}
    \label{eq:4.2.29}
    \frac{d^k}{dz^k} f(z) = \sum_{n \geq 0} c_n (-1)^k
    \frac{\Gamma(\alpha_n+1)}{\Gamma(\alpha_n+1-k)} 
    (1-z)^{\alpha_n-k} + O( |1-z|^{A-k}).
  \end{equation*}
\end{theorem}
\begin{proof}
  Clearly, all that is required is to establish
  \begin{equation*}
    \label{eq:4.2.30}
    \frac{d}{dz} O( |1-z|^A ) = O( |1-z|^{A-1} ).
  \end{equation*}
  Let us define $g(z) := f(z)/(1-z)^A$ where $f(z) = O( |1-z|^A)$. By
  differentiation, we find
  \begin{equation*}
    \label{eq:4.2.31}
    g'(z) = \frac{f'(z)}{(1-z)^A} + A \frac{f(z)}{(1-z)^{A+1}}.
  \end{equation*}
  This can be rewritten as
  \begin{equation*}
    \label{eq:4.2.32}
    f'(z) = g'(z)(1-z)^A - A\frac{f(z)}{1-z}.
  \end{equation*}
  Now $g(z)$ is by assumption $O(1)$ at $z=1$. If we assume that
  $g'(z) = O( |1-z|^{-1})$ at $z=1$, we obtain the sought result:
  $f'(z) = O( |1-z|^{A-1})$.
  Thus we need only prove the implication
  \begin{equation}
    \label{eq:4.2.33}
    g(z) = O(1) \implies g'(z) = O( |1-z|^{-1}),
  \end{equation}
  or, in more symbolic terms,
  \begin{equation*}
    \frac{d}{dz} O(1) = O\left(\frac{1}{1-z}\right).
  \end{equation*}
  
  We are going to prove~\eqref{eq:4.2.33} in a domain
  $\Delta(\phi+\epsilon,\eta)$, where $\epsilon > 0$ should be fixed
  but can be taken as small as we please.  The starting point is
  Cauchy integral formula
  \begin{equation*}
    \label{eq:4.2.34}
    g(z) = \frac{1}{2\pi i} \int_\gamma g(w) \frac{dw}{w-z},
  \end{equation*}
  which, through differentiation, provides
  \begin{equation*}
    \label{eq:4.2.35}
    g'(z) = \frac{1}{2\pi i} \int_\gamma g(w) \frac{dw}{(w-z)^2}.
  \end{equation*}
  \begin{figure}[htbp]
    \centering
    \includegraphics{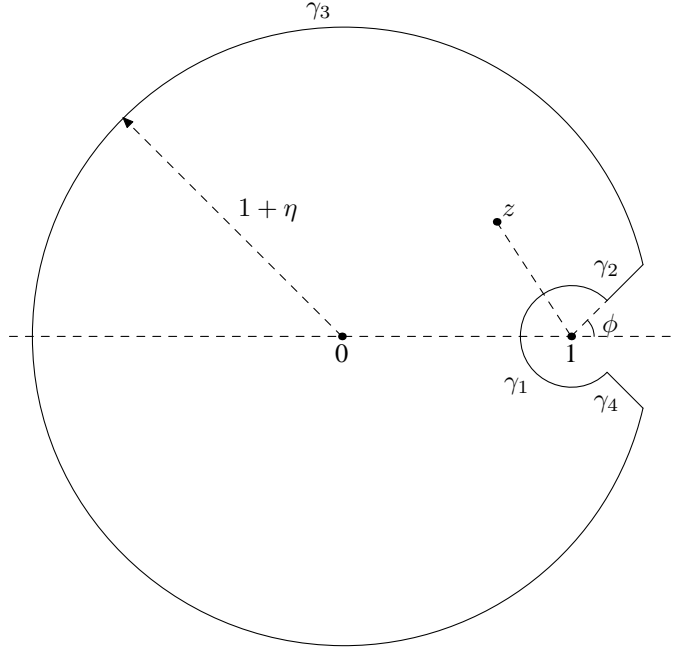}%
    \caption{The contour $\gamma$ used in the proof of the differentiation
      theorem.}%
    \label{fig:diff}%
  \end{figure}

  We select a small $\epsilon > 0$ and restrict $z$ to the domain
  $\Delta(\phi+\epsilon,\eta)$.
  The integral is estimated using the 
  contour $\gamma$
  depicted in Figure~\ref{fig:diff}: $\gamma := \gamma_1 \cup \gamma_2
  \cup \gamma_3 \cup
  \gamma_4$, with
  \begin{align*}
    \gamma_1 &:= \{ w\!: |w-1| = \tfrac{1}{2}|z-1|,\, |\arg(w-1)| \geq
    \phi\},\\
    \gamma_2 &:= \{ w\!: \tfrac{1}{2}|z-1| \leq |w-1|,\, |w| \leq 1+\eta,
    \arg(w-1) = \phi\},\\
    \gamma_3 &:= \{ w\!: |w| = 1 + \eta,\, |\arg(w-1)| \geq \phi\},\\
    \gamma_4 &:= \{ w\!: \tfrac{1}{2}|z-1| \leq |w-1|,\, |w| \leq 1+\eta,
    \arg(w-1) = -\phi \}.
  \end{align*}

  We then decompose the integral giving $g'(z)$ as $\int_\gamma =
  \int_{\gamma_1} +
  \int_{\gamma_2} + \int_{\gamma_3} + \int_{\gamma_4}$. Along 
  $\gamma_1$ we apply
  trivial
  upper bounds:
  \begin{equation}
    \label{eq:4.2.36}
    \left|\int_{\gamma_1}\right| = O(1) \cdot O\left(\frac{1}{ 
        |1-z|^2}\right) \cdot
    O(|1-z| ) = O( |1-z|^{-1}).
  \end{equation}
  The integral along the external circle $\gamma_3$ satisfies trivially
  \begin{equation}
    \label{eq:4.2.37}
    \left|\int_{\gamma_3}\right| = O(1) = O( | 1-z |^{-1}).
  \end{equation}
  If we set $|1-z|=\delta$, we find that the integrals along $\gamma_2$ and
  $\gamma_4$ are of no bigger order than
  \begin{equation*}
    \int_{\delta / 2}^\infty = \frac{du}{\delta^2 + u^2} =
    O\left(\frac{1}{\delta}\right) = O( |1-z|^{-1}).
  \end{equation*}
  Thus we have
  \begin{equation}
    \label{eq:4.2.38}
    \left|\int_{\gamma_2}\right|, \; \left|\int_{\gamma_4}\right| = O\left(
      |1 - z|^{-1} \right).
  \end{equation}
  Combining the estimates \eqref{eq:4.2.36}--\eqref{eq:4.2.38} yields the
  asymptotic form of $g'(z)$ as $O( |1-z|^{-1})$, as required at~\eqref{eq:4.2.33}.
\end{proof}

Thus, for instance, taking
\[
g(z) = \cos{\ln{\left(\frac{1}{1-z}\right)}} \quad \text{and} \quad
g'(z) = - \frac{1}{1-z}\sin{\ln{\left(\frac{1}{1-z}\right)}},
\]
we correctly predict that $g(z) = O(1) \Rightarrow g'(z) = O(
|1-z|^{-1})$. On the other hand, the apparent paradox given by the
pair
\begin{equation}
  \label{eq:4.2.56}
  g(z) = \cos{\exp{\left(\frac{1}{1-z}\right)}} \quad \text{and} \quad
  g'(z) = - \frac{1}{(1-z)^2}\exp{\left(\frac{1}{1-z}\right)}
  \sin{\exp{\left(\frac{1}{1-z}\right)}} ,
\end{equation}
is resolved by observing that in no nondegenerate sector around $z=1$
do we have $g(z) = O(1)$.

It is also well known that integration of asymptotic expansions is
usually easier than differentiation, and we can indeed prove
the following theorem.
\begin{theorem}
  \label{theorem:2B}
  If $f(z)$ has a $\Delta$-expansion that reduces to 
  $
  f(z) =
  O(|1-z|^A|\log{(1-z)}|^B)
  $,
  then we have
  \begin{equation*}
    \int_0^z f(t) \,dt =
    \begin{cases}
      O( |1-z|^{A+1}|\log{(1-z)}|^B ) &
      \text{if } A < -1\\
      O(|\log{(1-z)}|^{B+1}) & \text{if } A = -1,\, B \neq -1\\
      O(|\log{\log{(1-z)^{-1}}}|) & \text{if } A=-1,\, B=-1\\
      \displaystyle\int_0^1 f(t) \,dt +
      O(|1-z|^{A+1}|\log{(1-z)}|^{B}) & \text{if } A
      > -1.
    \end{cases}
  \end{equation*}
\end{theorem}
\begin{proof}
  The $\Delta$-expansion condition tells us that $|f(z)| \leq
  K|1-z|^A|\ln{(1-z)}|^B$ for some $0 \leq K < \infty$.
  
  Suppose $A < -1$.  By Cauchy integral formula, we can choose any
  path of integration that stays within the region of analyticity of
  $f$. We choose the contour $\gamma := \gamma_1 \cup \gamma_2$, shown in
  Figure~\ref{fig:int}.
  \begin{figure}[htbp]
    \centering
    \includegraphics{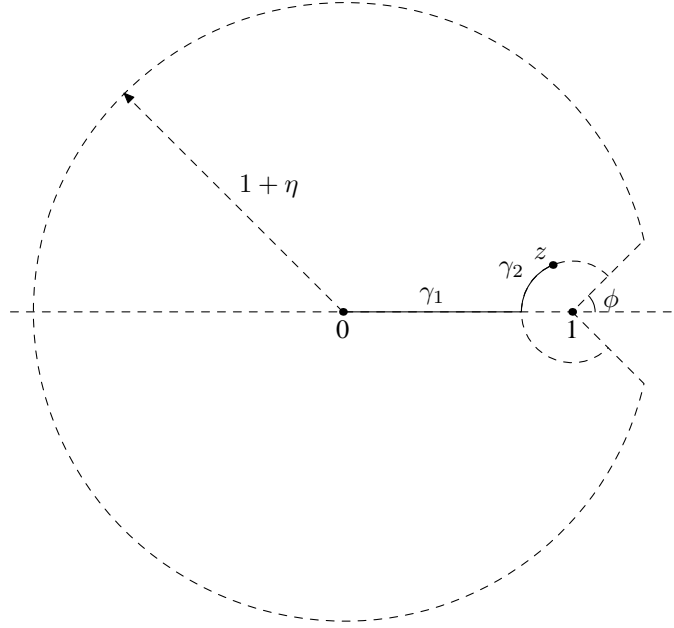}%
    \caption{The contour used in the proof of the integration theorem.}%
    \label{fig:int}%
  \end{figure}
  
  Then
  \begin{align*}
    \left|\int_{\gamma} f(t) \,dt \right| &= \left| \int_{\gamma_1} f(t) \,dt +
      \int_{\gamma_2} f(t) \,dt \right|
    \leq  \left| \int_{\gamma_1} f(t) \,dt \right| + \left|
      \int_{\gamma_2} f(t) \,dt \right|\\
    &\leq  K \int_{\gamma_1} |1-t|^A|\ln{(1-t)}|^B \,dt + K
    \int_{\gamma_2} |1-t|^A|\ln{(1-t)}|^B \,|dt|.
  \end{align*}
  The last integral along $\gamma_2$ is
  \[
  |1-z|^A |\ln{(1-z)}|^B \int_{\gamma_2} |dt| =
  O(|1-z|^{A+1}|\log{(1-z)}|^{B}).
  \]
  The last integral along $\gamma_1$ is a real integral.  Define $x := 1 -
  |1-z|$, so that the last integral along $\gamma_1$ equals
  \begin{align*}
    & \int_{0}^x (1-t)^{A} (\ln{(1-t)^{-1}})^B \,dt\\
    &= \frac{1}{|A|-1} \left[(1-x)^{A+1}(\ln{(1-x)^{-1}})^B
      - B\int_{0}^x (1-t)^A(\ln{(1-t)^{-1}})^{B-1} \,dt\right]\\
    &= \frac{1}{|A|-1} (1-x)^{A+1} (\ln{(1-x)^{-1}})^B +
    o\left(\int_{0}^x (1-t)^{A} (\log{(1-t)^{-1}})^B \,dt\right).
  \end{align*}%
  We have shown that
  \begin{align*}
    \int_{\gamma_1} |1-t|^A|\ln{(1-t)}|^B\,dt &= (1 + o(1)) \frac{1}{|A|-1}
    (1-x)^{A+1}(\ln{(1-x)^{-1}})^B\\
    &= (1 + o(1)) \frac{1}{|A|-1} |1-z|^{A+1} (\ln{|1-z|^{-1}})^B\\
    &= (1 + o(1)) \frac{1}{|A|-1} |1-z|^{A+1} |\ln{(1-z)^{-1}}|^B\\
    &= \Theta( |1-z|^{A+1} | \log{(1-z)^{-1}} |^B )\\
    &= O( |1-z|^{A+1} | \log{(1-z)^{-1}} |^B ).
  \end{align*}
  Hence
  \[ \int_{\gamma} f(t) \,dt = O(|1-z|^{A+1} |\log{(1-z)}|^B).\]
  
  The case $A=-1$ is similar.  As above,
  \[ \left|\int_{\gamma} f(t) \,dt\right| \leq K \int_{\gamma_1} |1-t|^{-1}
  |\ln{(1-t)}|^B \,dt + K \int_{\gamma_2} |1-t|^{-1} |\ln{(1-t)}|^B \,|dt|.
  \]
  The integral here along $\gamma_2$ is
  \[
  |1-z|^{-1}|\ln{(1-z)}|^B \int_{\gamma_2} \,|dt| = O(|\log{(1-z)}|^B) = o(
  |\log{(1-z)}|^{B+1}).
  \]
  The integral along $\gamma_1$ is a real integral.  If $B \neq -1$, then
  setting $x := 1-|1-z|$, we integrate to get
  \begin{equation*}
    \int_{0}^x(1-t)^{-1}(\ln{(1-t)^{-1}})^B\,dt
    = \frac{1}{B+1}(\ln{(1-x)^{-1}})^{B+1} = O(|\log{(1-z)}|^{B+1}).
  \end{equation*}
  If $B=-1$,
  \begin{equation*}
    \left|\int_{\gamma_1} f(t)\,dt\right|
    \leq O(1) + K \int_{1/2}^x
    (1-t)^{-1}(\ln{(1-t)^{-1}})^{-1}\,dt
    = O( |\log{\log{(1-z)^{-1}}}| ).
  \end{equation*}
  Combining the two integrals gives us the result.
  
  The case $A > -1$ is obtained by
  observing that the asymptotic condition guarantees the existence of
  $\int_0^1$ and by decomposing $\int_0^z = \int_0^1 + \int_1^z$. The contour
  of integration is the line segment joining 1 to $z$.  It is enough
  to show that
  \[
  \left| \int_z^1 f(t)\,dt \right| = O(|1-z|^{A+1}|\log{(1-z)}|^{B}).
  \]
  Indeed, performing the change of variable $t = uz + (1-u)$, we can
  bound the absolute value of this integral by $K$ times
  \begin{align*}
    &\int_0^1 |1-z|^A u^A |\ln{(u(1-z))}|^B |1-z|\,du\\
    &=
    |1-z|^{A+1}|\ln{(1-z)}|^{B} \int_0^1 u^A \left(
      \frac{|\ln{(u(1-z))}|}{|\ln{(1-z)}|} \right)^B\, du
    =
    O(|1-z|^{A+1}|\log{(1-z)}|^{B}),
  \end{align*}
  since the last integral is bounded.
\end{proof}

\section{Hadamard products}
\label{sec:hadamard-products}
In this section we  examine the way singular expansions are
composed under Hadamard products. The Hadamard product is a bilinear
form. So if we have a set of functions 
admitting known singular expansions,
we need to establish their composition law, and this will give
composition rules for finite terminating expansions. In order to
extend this to asymptotic expansions, either finite with $O( |1-z|^A
)$ error terms or infinite, we need to establish a theorem providing
the form of
\[
O( |1-z|^A ) \odot O( |1-z|^B ).
\]

The composition rule for polylogarithms~[see~\eqref{eq:4.1.50}] is
trivial, since
\[
\Li_{\alpha,r}(z) \odot \Li_{\beta,s}(z) = \Li_{\alpha+\beta,r+s}(z).
\]
However, polylogarithms do not have a simple composition rule with
respect to ordinary products.

We thus turn to the composition rule for the basis formed by functions
of the form $(1-z)^{-\alpha}$, where $\alpha$ may be any
complex
number.  From the expansion
\begin{equation*}
  \label{eq:4.3.39}
  (1-z)^{-\alpha} = 1 + \frac{\alpha}{1} z + 
  \frac{\alpha(\alpha+1)}{2!} z^2 + \cdots
\end{equation*}
around the origin,
we get through term-by-term multiplication,
\begin{equation}
  \label{eq:4.3.40}
  (1-z)^{-\alpha} \odot (1-z)^{-\beta} = 1 + \frac{\alpha \beta}{(1!)^2} z +
  \frac{\alpha(\alpha+1)\beta(\beta+1)}{(2!)^2} z^2 + \cdots.
\end{equation}
The classical \emph{hypergeometric function} of Gauss is defined for
$\gamma$ not a nonpositive integer by
\begin{equation*}
  \label{eq:4.3.43}
  \FF[\alpha,\beta;\gamma;z] = 1 + \frac{\alpha\beta}{\gamma} \frac{z}{1!} +
  \frac{\alpha(\alpha+1)\beta(\beta+1)}{\gamma(\gamma+1)} 
  \frac{z^2}{2!} + \cdots,
\end{equation*}
so that~\eqref{eq:4.3.40} reduces to
\begin{equation}
  \label{eq:4.3.41}
  (1-z)^{-\alpha} \odot (1-z)^{-\beta} = \FF[\alpha,\beta;1;z].
\end{equation}
From the transformation theory of hypergeometrics---see,
e.g.,~\cite[p.~163]{MR56:12235}---we know that, in general,
hypergeometrics can be expanded in the vicinity of $z=1$ by means of
the $z \mapsto 1 - z$ transformation. If we instantiate this
transformation with $\gamma=1$, we get
\begin{multline}
  \label{eq:4.3.42}
  \FF[\alpha,\beta;1;z] = 
  \frac{\Gamma(1-\alpha-\beta)}{\Gamma(1-\alpha)\Gamma(1-\beta)}
  \,\FF[\alpha,\beta;\alpha+\beta;1-z]\\
  + \frac{\Gamma(\alpha+\beta-1)}{\Gamma(\alpha)\Gamma(\beta)} 
  (1-z)^{-\alpha-\beta+1}
  \,\FF[1-\alpha,1-\beta;2-\alpha-\beta;1-z].
\end{multline}
In other words, we can state the following theorem.
\begin{proposition}
  \label{theorem:3}
  When \( \alpha\), \(\beta\), and $\alpha+\beta$ are not integers,
  the Hadamard product
  \[ (1-z)^{-\alpha} \odot (1-z)^{-\beta} \]
  has an infinite $\Delta$-expansion with exponent scale
  \(
  \{0,1,2,\ldots\} \cup \{-\alpha-\beta+1, -\alpha-\beta+2,\ldots\}
  \),
  namely,
  \[
  (1-z)^{-\alpha} \odot (1-z)^{-\beta} \sim \sum_{k \geq 0} 
  \lambda_k^{(\alpha,\beta)}
  \frac{(1-z)^k}{k!} + \sum_{k \geq 0} \mu_k^{(\alpha,\beta)}
  \frac{(1-z)^{-\alpha-\beta+1+k}}{k!},
  \]
  where the
  coefficients
  $\lambda$ and $\mu$ are given by
  \begin{align*}
    \lambda_k^{(\alpha,\beta)} &=
    \frac{\Gamma(1-\alpha-\beta)}
    {\Gamma(1-\alpha)\Gamma(1-\beta)}
    \frac{\alpha^{\bar{k}}\beta^{\bar{k}}}{(\alpha+\beta)^{\bar{k}}},\\
    \mu_k^{(\alpha,\beta)} &= 
    \frac{\Gamma(\alpha+\beta-1)}{\Gamma(\alpha)\Gamma(\beta)}
    \frac{(1-\alpha)^{\bar{k}}
    (1-\beta)^{\bar{k}}}{(2-\alpha-\beta)^{\bar{k}}}.
  \end{align*}
\end{proposition}

We are thus facing a situation where polylogarithms are simple for
Hadamard products and complicated for ordinary products, with the reverse
situation for power functions.

We now examine how $O(\cdot)$ terms get composed. The task is easier when
the resulting function becomes large at its singularity. Fortunately,
thanks to the results of
Section~\ref{sec:sing-expans-diff}, all cases can be reduced (see
Proposition~\ref{theorem:4}) to the
following.
\begin{proposition}
  \label{lem:3}
  Assume that $f(z)$ and $g(z)$ are $\Delta$-regular in
  $\Delta(\psi_0, \eta)$.  Additionally, suppose that
  \[
  f(z) = O(|1-z|^a) \text{ and } g(z) = O(|1-z|^b) \text{ as  } z \to
  1,
  \]
  where $a$ and $b$
  satisfy $a+b+1 < 0$. Then the Hadamard product $(f \odot g)(z)$ is
  $\Delta$-regular and admits the expansion
  \begin{equation*}
    \label{eq:4.3.44}
    (f \odot g)(z) = O( |1-z|^{a+b+1} ).
  \end{equation*}
\end{proposition}
\begin{proof}
  We first give the proof in the simpler case when%
  $f(z)$ and $g(z)$ are analytic in the whole of the complex plane
  slit along $z \geq 1$.

  Assume that $s$ is a complex number satisfying $|s| < 1$. Then
  consider the integral
  \begin{equation}
    \label{eq:4.3.45}
    I = \frac{1}{2\pi i} \int_{\gamma_0} f(w)g\left(\frac{s}{w}\right)
    \frac{dw}{w},
  \end{equation}
  taken (counterclockwise) along a contour~$\gamma_0$ (see
  Figure~\ref{fig:5-1a})
  which is simply a circle of
  radius $R$ centered at the origin and such that $|s| < R < 1$.
  Setting
  \[
  f(z) = \sum_{n \geq 0} f_nz^n \qquad \text{and} \qquad g(z) = \sum_{n \geq
    0} g_nz^n,
  \]
  and computing the integral~(\ref{eq:4.3.45}) by expanding the functions,
  we find
  \[
  I = \sum_{n \geq 0} f_ng_n s^n.
  \]
  This is the classical formula of Hadamard for Hadamard products.
  \begin{figure}[htbp]%
    \centering
    \includegraphics{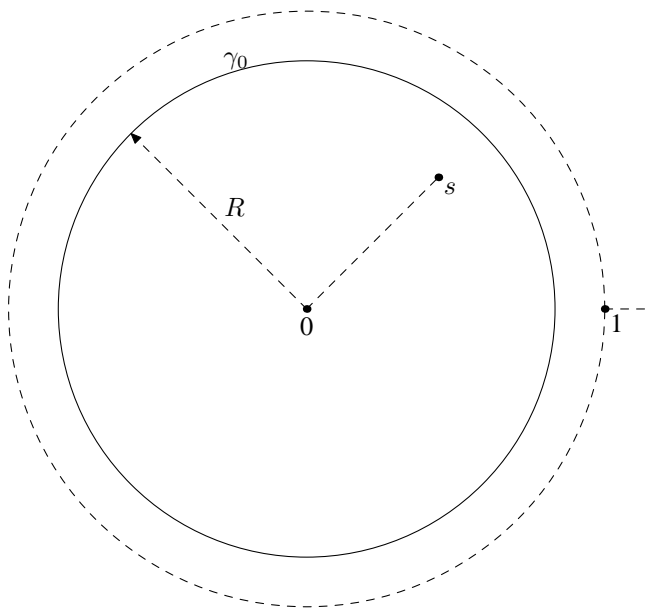}%
    \caption{The contour $\gamma_0$.}%
    \label{fig:5-1a}%
  \end{figure}

  We propose
  to continue $f \odot g$ to a point $z$ such that $|z| > 1$ where
  $z/s$ is real and positive. Because of the analytic continuation
  properties of $f$ and $g$, we can take as the contour of integration
  in~(\ref{eq:4.3.45}) the contour $\gamma$ depicted in
  Figure~\ref{fig:5-1b}.  The quantity $f \odot g$ evaluated at
    $s$ is still equal to the integral form along $\gamma$ since
  there is a continuous deformation from $\gamma_0$ to $\gamma$ such
  that, along the intermediate contours, both $w$ and $s/w$ stay
  within the domains of analyticity of $f$ and $g$.
  \begin{figure}[ht]%
    \centering
    \includegraphics{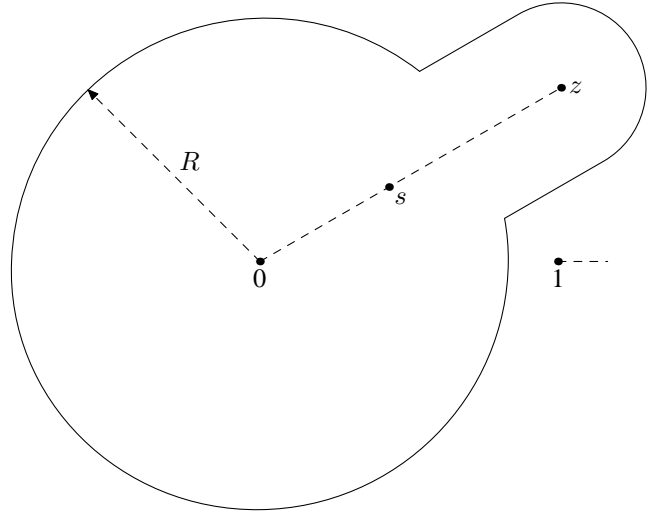}%
    \caption{The contour $\gamma$.}%
    \label{fig:5-1b}%
  \end{figure}

  We now consider the integral in~(\ref{eq:4.3.45}) as a function of the
  parameter $s$, and let $s$ vary continuously from $z$ along the ray
  from the origin to $z$. The integrand remains an analytic function of
  both
  $s$ and $z$ so that the integral~(\ref{eq:4.3.45}) provides the
  analytic continuation of $f \odot g$. Using this contour, we thus capture
  the singular behaviors of both $f$ and $g$.

  In order to obtain the asymptotics for $f \odot g$, we fully specify
  in Figure~\ref{fig:5-1c} the geometry of the contour $\gamma$.  Without
  loss of generality, we assume that $z$ is in the upper half-plane.
  \begin{figure}[htbp]%
    \centering
    \includegraphics{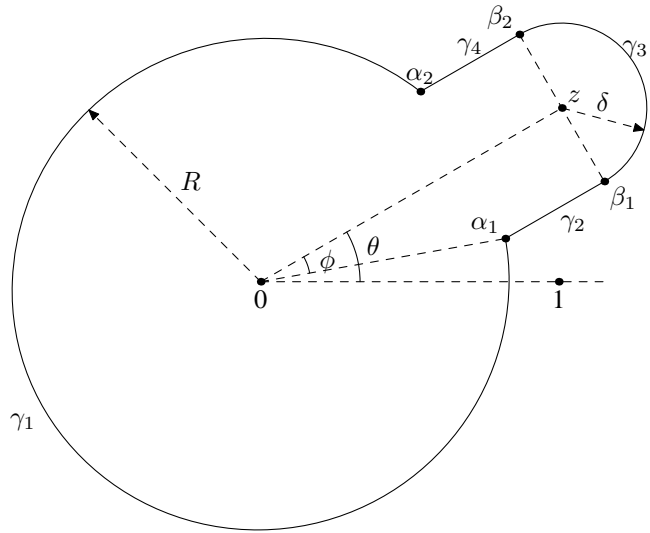}%
    \caption{The geometry of the contour $\gamma$.}%
    \label{fig:5-1c}%
  \end{figure}

  Define
  \begin{align*}
    \theta &:= \arg{z}, & \phi &:= \arcsin{\frac{c\delta}{R}},\\
    \alpha_1 & := Re^{i(\theta - \phi)}, & \beta_1 &:= z - c \delta
    e^{i(\theta +
      \frac{\pi}{2})},\\
    \alpha_2 & := Re^{i(\theta + \phi)}, &\beta_2 &:= z + c \delta
    e^{i(\theta + \frac{\pi}{2})},
  \end{align*}
  where we can (it is not hard to see) and do choose $0 <c \leq
  \tfrac14\sin\psi_0$ small enough that $\gamma$ (and its interior) stays
  entirely within the interior of $\Delta(\psi_0,\eta)$. (The reason for
  the restriction $c \leq \tfrac14\sin\psi_0$ will become clear later.)


  Then, the contour $\gamma$ is $\gamma_1 \cup \gamma_2 \cup \gamma_3
  \cup \gamma_4$, where
  \begin{align*}
    \gamma_1 &= \{ w\!: |w| = R,\; \arg{w} \in [0,\theta - \phi] \cup [\theta
    + \phi, 2\pi]\},\\
    \gamma_2 &= \{ w\!: w = (1-t)\alpha_1 + t\beta_1,\; 0 \leq t \leq 1\},\\
    \gamma_3 &= \{ w\!: |w-z| = \delta/2,\; \arg(w - z) \in [0,\pi]\},\\
    \gamma_4 &= \{ w\!: w = t\alpha_2 + (1-t)\beta_2,\; 0 \leq t \leq 1 \}.
  \end{align*}
  By the asymptotic conditions on $f$ and $g$, we know that there are
  positive constants $K_1$ and $K_2$ such that $f(w) \leq K_1 |1-w|^a$
  and $g(w) \leq K_2 |1-w|^b$ for all $w$ on $\gamma$.  We also note that $R
  \leq |w| \leq 1+\eta$, so we can find a constant $K$ such that
  \begin{align}
    \label{eq:4.3.49}
    |(f \odot g)(z)| &= \left|\frac{1}{2\pi i} \int_{\gamma}\,f(w)
      g\left(\frac{z}{w}\right) \frac{dw}{w}\right| \notag
     \leq \frac{K}{2\pi} \int_{\gamma}\,|1 - w|^a |w - z|^b\,|dw| \notag\\
    & \leq \frac{K}{2\pi} \left( \int_{\gamma_1} + \int_{\gamma_2} +
      \int_{\gamma_3} + \int_{\gamma_4} \right) |1 - w|^a |w - z|^b\,|dw|.
  \end{align}
  We now show that each of the four integrals in~\eqref{eq:4.3.49} is
  $O(\delta^{a + b + 1}) = O( |1-z|^{a+b+1})$, completing the proof 
  of~\eqref{eq:4.3.44}.
  \begin{enumerate}
  \item Uniformly along $\gamma_1$, the integrand is $\Theta(1)$. 
    Hence the integral is
    $\Theta(1)$, which is $O(\delta^{a+b+1})$ since $a+b+1 < 0$.
  \item Along $\gamma_2$, we can express the variable of integration $w$ as%
    \[
    w \equiv w(u) = (\alpha_1 - \beta_1)u + \beta_1, \quad 0 \leq u \leq 1,
    \]
    so that
    \[ dw = (\alpha_1 - \beta_1)du.\]
    Note that $|\beta_1 - \alpha_1| = (1 + o(1))(1-R)$.
    The integral of interest is then bounded by a constant times
    \[
    \int_{0}^\infty |1 - w(u)|^a|w(u) - z|^b\,du .
    \]

    We split this integral into two pieces, one where $u \leq M\delta$, say
    $\gamma_{21}$, and the other where $u > M\delta$, say
    $\gamma_{22}$, where $M$ is a 
    suitably chosen large constant. We will show that each of these
    integrals is $O(\delta^{a+b+1})$.  To this end, observe that
    \begin{align}
      1 - w &= 1 - \beta_1 + (\beta_1 - \alpha_1)u
      = 1 - z + c\delta e^{i(\theta + \tfrac\pi2)} + u |\beta_1 -
      \alpha_1|e^{i\theta} \notag\\
      &= e^{i\theta}( e^{-i\theta} - |z| + i c\delta + u|\beta_1 -
      \alpha_1| ), \label{eq:4.3.57}
    \end{align}
    We defer to Lemma~\ref{lemma:3.5} a calculation that establishes
    $|1-w| = \Omega(\delta)$ uniformly along all of $\gamma_2$. To get an
    upper bound, using the triangle inequality,
    \begin{equation*}
      |1-w| \leq |1-\beta_1| + |\beta_1-w|
      \leq \delta + c\delta + |\beta_1-w| = (c+1)\delta + |\beta_1-w|.
    \end{equation*}
    Restricting attention to $\gamma_{21}$ gives $|\beta_1-w| \leq
    |\beta_1-\alpha_1|M\delta$.  Hence, uniformly
    along $\gamma_{21}$, $|1-w| = \Theta(\delta)$.
    
    For $|w-z|$, a trivial lower bound is $|w-z| \geq |\beta_1-z| =
    c\delta$. Along $\gamma_{21}$, an upper bound is established by
    the triangle inequality as
    \[
    |w-z| \leq |\beta_1-z| + |w-\beta_1| \leq c\delta +
     |\beta_1-\alpha_1|M\delta.
    \]
    Hence $|w-z| = \Theta(\delta)$ uniformly along $\gamma_{21}$, so
    that the integral along $\gamma_{21}$ is $\Theta(\delta^{a+b+1})$
    since the integrand is uniformly $\Theta(\delta^{a+b})$ and the
    path of integration is of length $\Theta(\delta)$.
    
    For the integral along $\gamma_{22}$, we will establish that
    $|1-w| = \Theta(u)$ and $|w-z| = \Theta(u)$ for suitably large
    $M$.  The desired bound on the integral will then follow, since
    $a+b+1 < 0$.
    
    We use~\eqref{eq:4.3.57} and the fact that the real part of a
    complex number is a lower bound on its absolute value to get
    \[
    |1-w| \geq \cos{\theta} - |z| + |\beta_1-\alpha_1|u.\ 
    \]
    Writing $z=1+ \delta e^{i\psi}$, we note that
    \[
    \cos\theta = \frac{1+\delta\cos\psi}{|z|} \geq
    \frac{1+\delta\cos\psi}{1+\delta} \geq 
    \frac{1-\delta}{1+\delta} \geq 1 - 2\delta.
    \]
    Hence
    \begin{equation*}
      |1-w| \geq 1 - 2\delta - (1+\delta) + (1 + o(1))(1-R)u
      = (1 + o(1))(1-R)u -3\delta.
    \end{equation*}
    Since $u \geq M\delta$, if we choose $M$ large enough, we get $|1-w| =
    \Omega(u)$.

    For the upper bound, 
    \begin{equation*}
      |1-w|  \leq |w-\beta_1| + |1-\beta_1|
      = |\alpha_1 - \beta_1|u + |1-\beta_1|
      \leq  (1 + o(1))(1-R)u + (c+1)\delta
      = O(u).
    \end{equation*}
    Hence $|1-w| = \Theta(u).$

    The estimate for $|w-z|$ is similar.  For the upper bound, observe that
    \begin{align*}
      w-z &= \beta_1 - z + (\alpha_1 - \beta_1)u, \intertext{whence}
      |w-z| &\leq |\beta_1-z| + |\alpha_1-\beta_1|u
      c\delta + (1 + o(1))(1-R)u = O(u).
    \end{align*}
    For the lower bound, we express
    \[ w-z = -c\delta e^{i(\theta+\tfrac{\pi}{2})} +
    u|\alpha_1-\beta_1|e^{i\theta} =
    e^{i\theta}(u|\alpha_1-\beta_1| - ic\delta), \]
    so that
    \begin{equation*}
      |w-z| = |u|\alpha_1-\beta_1| - ic\delta|
      \geq u|\alpha_1 - \beta_1|
      = (1 + o(1))(1-R)u.
    \end{equation*}
    This establishes  $|w-z| = \Omega(u)$ and hence $|w-z| = \Theta(u)$.
    
  \item Along $\gamma_3$, $|w-z| = c\delta$, and $(1-c)\delta = \delta
    - c\delta \leq |1-w| \leq \delta + c\delta= (1+c)\delta$ using the
    triangle inequality.  The contour of integration has length $\pi
    c\delta$. Hence the integral is $\Theta(\delta^{a+b+1})
    = O(|1-z|^{a+b+1} )$.
  \item Along $\gamma_4$, the argument is similar to (and somewhat simpler
    than) that for the integral
    along $\gamma_2$.
  \end{enumerate}
\end{proof}
We now present the proof of the lower-bound on $|1-w|$ along
$\gamma_2$ used in the proof of Proposition~\ref{lem:3}.
\begin{lemma}
  \label{lemma:3.5}
  In the notation of the proof of Proposition~\ref{lem:3}, $|1-w| =
  \Omega(\delta)$ uniformly along all of $\gamma_2$.
\end{lemma}
\begin{proof}
  If $\sin\theta > 2c\delta$, we lower-bound $|1-w|$
  using~\eqref{eq:4.3.57} and the fact that the imaginary part of a
  complex number is a lower bound on the its absolute value,
  \[
  |1-w| \geq -\Im(e^{-i\theta} - |z| + ic\delta + u|\beta_1-\alpha_1|)
  = \sin\theta - c\delta 
  > c\delta.
  \]
  Otherwise,  we use the real part as a lower bound:
  \[
  |1-w| \geq \Re(e^{-i\theta} - |z| + ic\delta + u|\beta_1-\alpha_1|)
  \geq \cos\theta - |z| + u|\beta_1-\alpha_1|.
  \]
  We choose $c$ small enough so that if $\Re(z) = 1 - |z|\cos\theta <
  0$, then $\sin\theta > 2c\delta$.  Indeed, if $\Re(z) > 1$, then $z
  = 1 + \delta e^{i\psi}$ where $\psi_0 < \psi < \pi/2$, and so
  \[ 
  \sin\theta = \frac{\delta\sin\psi}{|z|} > \frac{\delta\sin\psi_0}{2}
  \]
  and any $c \leq \frac{1}{4}\sin\psi_0$ will suffice.  So we may
  assume $$\sin\theta \leq 2c\delta \text{ and } 1 - |z|\cos\theta
  \geq 0.$$
  
  Now $\cos^2\theta \geq 1 - 4c^2\delta^2$, whence $\cos\theta = 1 -
  o(\delta)$, and we need only prove that $1 - |z| \geq 0$ and $1 -
  |z| \geq \Omega(\delta)$.  To this end, observe that $\delta^2 =
  |1-z|^2 = (1-|z|\cos\theta)^2 + (|z|\sin\theta)^2$, so that
  \begin{align*}
    1 = |z|\cos\theta + (1 - |z|\cos\theta)
    &\geq |z|\sqrt{1- 4c^2\delta^2} + \delta\sqrt{1-4c^2|z|^2}\\
    &\geq |z|(1 - 4c^2\delta^2) + \delta(1 - 4c^2|z|^2),
  \end{align*}
  and hence, for all sufficiently small $\delta$,
  \begin{align*}
    1 - |z| &\geq \delta(1 - 4c^2|z|^2) - 4|z|c^2\delta^2\\
    &\geq \delta(1-16c^2) - 8c^2\delta^2 \quad \text{(since
      $|z| \leq 2$ for $\delta \leq 1$)}\\
    &\geq  p\delta,
  \end{align*}
  where $p$ is a positive constant depending only on $\psi_0$, since
  $c \leq \tfrac14\sin\psi_0 < \tfrac14$. 
\end{proof}

Next we show how to modify Proposition~\ref{lem:3} whenever $a + b >
-2$.
\begin{proposition}
  \label{theorem:4}
  Assume that $f(z)$ and $g(z)$ are $\Delta$-regular and that
  \[
  f(z)  = O(|1-z|^a) \text{ and }g(z) = O(|1-z|^b) \text{ as  } z \to 1.
  \]
  \begin{enumerate}[(a)]
  \item If $k < a + b +1 < k+1$ for some integer $-1 \leq k < 
    \infty$, then
    \begin{equation*}
      \label{eq:4.3.46}
      (f \odot g)(z) = \sum_{j=0}^k \frac{(-1)^j}{j!}
      (f\odot{}g)^{(j)}(1)(1-z)^j + O(|1-z|^{a+b+1}).
    \end{equation*}
  \item If $a+b+1$ is a nonnegative integer then
    \begin{equation*}
      \label{eq:4.3.47}
      (f \odot g)(z) = \sum_{j=0}^{a+b} \frac{(-1)^j}{j!}
      (f\odot{}g)^{(j)}(1)(1-z)^j 
      + O( |1-z|^{a+b+1}|\log((1 - z)^{-1})|).
    \end{equation*}
  \end{enumerate}
\end{proposition}
\begin{proof}
  Let $D$ denote the operator $\tfrac{d}{dz}$ and let $\theta$ denote
  the operator $zD$.
  \begin{enumerate}[(a)]
    \item  We first show that
    \begin{equation}
      \label{eq:4.3.62}
      (f\odot{}g)(z) = \sum_{j=0}^k \frac{(-1)^j}{j!}
      \left[ \theta^j(f\odot{}g) \right](1)(-\ln{z})^j + O(|1-z|^{a+b+1}).
    \end{equation}
    The proof is by induction on $k$. The basis case of
    $k=-1$ is a special case of Proposition~\ref{lem:3}. Observe that
    \[
    \theta(f\odot{}g) = (\theta{}f)\odot{}g
    \]
    and by Theorem~\ref{theorem:2A} that
    \[
    (\theta{}f)(z) = zf'(z) = O(|1-z|^{a-1}).
    \]
    If $k \geq 0$, then by the induction hypothesis
    \[
    \left[ \theta(f\odot{}g) \right](\zeta) = (\theta{}f\odot{}g)(\zeta) 
    = \sum_{j=0}^{k-1} \frac{(-1)^j}{j!}
    \left[ \theta^j(\theta{}f\odot{}g) \right](1) (-\ln\zeta)^j + O(
    |1-\zeta|^{a+b}).
    \]
    We divide by $\zeta$ and then integrate both sides
    with respect to $\zeta$ from $z$ to $1$ to get
    \begin{align*}
      (f \odot g)(z) &= (f \odot g)(1) + \sum_{j=0}^{k-1}
      \frac{(-1)^{j+1}}{(j+1)!} \left[ \theta^{j+1}(f\odot{}g) \right](1)
      (-\ln z)^{j+1}
      + O(|1-z|^{a+b+1})\notag\\
      &= \sum_{j=0}^k \frac{(-1)^j}{j!}
      \left[ \theta^j(f\odot{}g) \right](1)(-\ln z)^j +
      O(|1-z|^{a+b+1}),\label{eq:4.3.61}
    \end{align*}
    as required.

    To obtain the expansion in the asserted form, we use the
    following two well-known facts: 
    \begin{enumerate}
    \item For fixed nonnegative integers $j$ and $k$, as $z \to 1$ we
      have
      \[
      (-\ln z)^j = j! \sum_{r=0}^k\stirlingone{r}{j}
      \frac{(1-z)^r}{r!} + O(|1-z|^{k+1}),
      \]
      where $\stirlingone{r}{j}$ denotes a Stirling number of the
      first kind (see~1.2.9-(26) in~\cite{knuth97}).
    \item The operators $D$ and $\theta$ satisfy
      \[
      z^rD^r = \sum_{j=1}^r(-1)^{r-j}
      \stirlingone{r}{j}\theta^j, \quad r \geq 1.
      \]
    \end{enumerate}
    Using these facts in~\eqref{eq:4.3.62} and recalling that $a+b+1 <
    k+1$, we get
    \begin{align*}
      (f\odot{}g)(z) &= \sum_{j=0}^k (-1)^j\left[ \theta^j(f\odot{}g)
      \right](1)
      \sum_{r=0}^k \stirlingone{r}{j} \frac{(1-z)^r}{r!} +
      O(|1-z|^{a+b+1})\\
      &= (f\odot{}g)(1) + \sum_{r=1}^k (-1)^r \frac{(1-z)^r}{r!}
      \sum_{j=1}^k (-1)^{r-j}
      \stirlingone{r}{j}\left[ \theta^j(f\odot{}g) \right](1)
      + O(|1-z|^{a+b+1})\\
      &= \sum_{r=0}^k (-1)^r
      \frac{(1-z)^r}{r!}\left.\bigr[ w^r(f\odot{}g)^{(r)}(w)
      \bigr]\right\rvert_{w=1}(1)
      + O(|1-z|^{a+b+1})\\
      &= \sum_{j=0}^k \frac{(-1)^j}{j!}
      (f\odot{}g)^{(j)}(1)(1-z)^j + O(|1-z|^{a+b+1}),
    \end{align*}
    as desired.
    \item
    The structure of the proof is the same as that of
    part~(a). One
    first proves a result similar to~\eqref{eq:4.3.62}.
    Again, the proof is by induction on $k := a+b+1$.  When
    $k=0$, after applying Theorem~\ref{theorem:2A} and
    Proposition~\ref{lem:3} to obtain
    $\theta(f \odot g)(z) = O(|1-z|^{-1})$ we can
    apply  Theorem~\ref{theorem:2B} to get
    $(f \odot g)(z) = O(|\log((1 - z)^{-1})|)$. 
    The induction step is similar to that in (a) above. Finally,
    using the same argument as for part~(a), we obtain the expansion in
    the asserted form.
  \end{enumerate}
\end{proof}

For convenience, we now combine Propositions~\ref{lem:3}
and~\ref{theorem:4} to get a result covering all possible values of
$a$ and $b$.
\begin{theorem}
  \label{thm:hadamard}
  If $f$ and $g$ are amenable to singularity analysis and
  \begin{equation*}
    f(z) = O(|1-z|^a) \quad\text{and}\quad g(z) = O(|1-z|^b) \quad
    \text{as $z \to 1$},
  \end{equation*}
  then  $f \odot g$ is also amenable to singularity
  analysis.  Furthermore,
  \begin{enumerate}[(a)]
  \item If $a+b+1 < 0$ then
    \begin{equation*}
      f(z) \odot g(z) = O(|1-z|^{a+b+1}).
    \end{equation*}\label{item:hadamard1}
  \item If $k < a+b+1 < k+1$ for some integer $-1 \leq k < \infty$,
    then
    \begin{equation*}
      f(z) \odot g(z) = \sum_{j=0}^k \frac{(-1)^j}{j!} (f \odot
      g)^{(j)}(1) (1-z)^j + O(|1-z|^{a+b+1}).
    \end{equation*}
    \label{item:hadamard2}
  \item   If $a+b+1$ is a nonnegative integer then
    \begin{equation*}
      f(z) \odot g(z) = \sum_{j=0}^{a+b} \frac{(-1)^j}{j!} (f \odot
      g)^{(j)}(1) (1-z)^j + O(|1-z|^{a+b+1}|\log{((1-z)^{-1})}|).
    \end{equation*}    \label{item:hadamard3}
  \end{enumerate}
\end{theorem}
We will make extensive use of Theorem~\ref{thm:hadamard} in the
subsequent chapters.  We mention in passing that the results of this
chapter are also applicable to the analysis of the Cayley tree
recurrence~\cite{MR89i:05024,MR81a:68049}, which models the evolution
of a graph from being totally disconnected to being tree-like, when
successive edges are added at random; see~\cite{FFK} for details.


\part{Applications}
\chapter{Random permutation model: shape functional}
\label{cha:binary-search-trees}
As our first application of the theoretical results of
Chapters~\ref{cha:transfer-theorems} and~\ref{cha:sing-analys-hadam},
we demonstrate how to get a complete asymptotic expansion for the mean
of the shape functional (see Example~\ref{example:shape-functional})
for \( m \)-ary search trees and we determine the asymptotics of its
variance for binary search trees.  Our results generalize those
obtained in~\cite{MR97f:68021}.

\section{Complete asymptotic expansion of the mean}
\label{sec:bst_mean-shape-funct}
As noted in Example~\ref{example:shape-functional}, the shape
functional corresponds to the toll function \( t_n := \ln{
  \binom{n}{m-1} } \).  In order to obtain a complete asymptotic
expansion of its mean using the ETT (Theorem~\ref{theorem:ETT}), our
first task is to obtain a complete asymptotic expansion for its
generating function
\begin{equation*}
  T(z) := \sum_{n=m-1}^\infty t_n z^n.
\end{equation*}

\subsection*{The generating function of the toll sequence}
\label{sec:gener-funct-toll}

By definition
\begin{align*}
  T(z) &= \sum_{n=m-1}^\infty t_nz^n = \sum_{n=m-1}^\infty
  \left[\ln{\binom{n}{m-1}}\right] z^n\\ &= \sum_{n=m-1}^\infty \left[
  \sum_{l=0}^{m-2} \ln{(n-l)}\right] z^n -
  \left[\ln{(m-1)!}\right]\sum_{n=m-1}^\infty z^n = P(z) - Q(z),
\end{align*}
where
\begin{equation*}
  P(z) := \sum_{n=m-1}^\infty \left[ \sum_{l=0}^{m-2} \ln{(n-l)}
  \right] z^n
\end{equation*}
and
\begin{equation*}
  Q(z) := \left[\ln{(m-1)!}\right]
  \sum_{n=m-1}^\infty z^n = \left[\ln{(m-1)!}\right]\left[(1-z)^{-1} -
  \sum_{n=0}^{m-2} z^n\right].
\end{equation*}
But
\begin{align}
  P(z) = \sum_{l=0}^{m-2} \sum_{n=m-1}^\infty
  \left[\ln{(n-l)}\right]z^n
  &= \sum_{l=0}^{m-2} z^l \sum_{n=l+1}^\infty \left[ \ln{(n-l)}
  \right]z^{n-l} - 
  \sum_{l=0}^{m-2} \sum_{n=l+1}^{m-2} \left[ \ln{(n-l)} \right] z^{n} \notag\\
  \label{eq:5.1.4}
  &= \frac{1-z^{m-1}}{1-z} \Li_{0,1}(z) -
  \sum_{n=1}^{m-2}z^n \sum_{l=0}^{n-1} \ln{(n-l)},
\end{align}
where \(\Li\) denotes the generalized polylogarithm defined
at~\eqref{eq:4.1.50}.
The term subtracted in~\eqref{eq:5.1.4} equals
\begin{multline}
  \label{eq:5.1.1}
  \sum_{n=1}^{m-2} z^n \ln{n!}
  = \sum_{n=1}^{m-2} \sum_{k=0}^n (-1)^k \binom{n}{k}
  (1-z)^k \ln{n!}\\
  = \sum_{k=0}^{m-2} (-1)^k (1-z)^k \sum_{n=\max{\{k,1\}}}^{m-2}
  \binom{n}{k} \ln{n!} = \sum_{k=0}^{m-2} (-1)^k (1-z)^k \sum_{n=k}^{m-2}
  \binom{n}{k} \ln{n!},
\end{multline}
while the first term in~\eqref{eq:5.1.4} is the product of
\((1-z)^{-1} \Li_{0,1}(z)\) and
\begin{equation*}
  1-z^{m-1} 
  = 1 - \sum_{k=0}^{m-1}(-1)^k
  \binom{m-1}{k} (1-z)^k
  = \sum_{k=1}^{m-1} (-1)^{k+1} \binom{m-1}{k} (1-z)^k,
\end{equation*}
i.e., it is \(\Li_{0,1}(z)\) times
\begin{equation}
  \label{eq:5.1.2} \frac{1-z^{m-1}}{1-z} = \sum_{k=0}^{m-2} (-1)^k
  \binom{m-1}{k+1} (1-z)^k.
\end{equation}
Using Theorem~\ref{theorem:Li-sing} we know that
\begin{equation}
  \label{eq:5.1.3} \Li_{0,1}(z) \sim t^{-1}\ln{t^{-1}} + \Gamma'(1)t^{-1}
  - \sum_{j \geq 0} \frac{(-1)^j}{j!}\zeta'(-j)t^j
\end{equation}
as \(z \to 1\) in a suitable indented crown, with \(t :=
-\ln{z}\).  We note that \(\Gamma'(1) = -\gamma\), where \(\gamma
\doteq 0.5772\) is Euler's constant.

In order to carry out the integration required in the ETT, we need to
express \(T(z)\) in terms of powers of \((1-z)\) and
\(\ln{((1-z)^{-1})}\).  This is our first goal.  [For notational
convenience we will define $L(z) := \ln{((1-z)^{-1})}$.]  To this end,
expression \eqref{eq:5.1.1} can be written as
\begin{equation}
  \label{eq:5.1.5} \sum_{k=0}^{m-2} S_1(k) (1-z)^k
  \quad\text{ where }\quad
  S_1(k) := (-1)^k \sum_{n=k}^{m-2} \binom{n}{k}
  \ln{n!}.
\end{equation}
Expression \eqref{eq:5.1.2} can be written as
\begin{equation}
  \label{eq:5.1.2.3} \sum_{k=0}^{m-2} S_2(k) (1-z)^k,
  \quad\text{ where }\quad
  S_2(k) := (-1)^k \binom{m-1}{k+1}.
\end{equation}
For~\eqref{eq:5.1.3}, we will need the following two facts: First,
\begin{equation}
  \label{eq:5.1.6.5}
t^{j} = \sum_{r=j}^\infty T_1(r,j)(1-z)^r, \quad j \geq -1
\quad\text{ where }\quad
T_1(r,i) := 
\begin{cases}
  [z^r]L^i(z) & i \neq 0\\ \delta_{r,0} & i = 0.
\end{cases}
\end{equation}
Next,
\begin{equation}
  \label{eq:5.1.6.75}
\ln{t^{-1}} = L(z) + \sum_{r=1}^\infty T_2(r)(1-z)^r,
\quad\text{ where }\quad
T_2(r) := \sum_{i=1}^r \frac{(-1)^i}{i}
T_3(r+i,i),
\end{equation}
with
\begin{equation}
  \label{eq:5.1.6.875}
T_3(r,i) := [w^r]\left[L(w) - w \right]^i.
\end{equation}
The final ingredient in our expression for \(T(z)\) is
\begin{align}
  \label{eq:5.1.6} \sum_{n=0}^{m-2} z^n &= \sum_{n=0}^{m-2}
  (1-(1-z))^n
  = \sum_{n=0}^{m-2} \sum_{k=0}^n \binom{n}{k}
  (-1)^k (1-z)^k\notag\\
  &= \sum_{k=0}^{m-2} (1-z)^k (-1)^k\sum_{n=k}^{m-2} \binom{n}{k}
  = \sum_{k=0}^{m-2} S_3(k)(1-z)^k,
\end{align}
where
\begin{equation}
  \label{eq:5.1.3.5} S_3(k) := (-1)^k \sum_{n=k}^{m-2} \binom{n}{k}.
\end{equation}

We now express \(\Li_{0,1}(z)\) in terms of \((1-z)\):
\begin{align*}
  \Li_{0,1}(z) &\sim \left(\sum_{r=-1}^\infty
    T_1(r,-1)(1-z)^r\right)\left(L(z) + \sum_{r=1}^\infty
    T_2(r)(1-z)^r\right)\notag\\
  & \quad + \Gamma'(1) \sum_{r=-1}^\infty T_1(r,-1)(1-z)^r
   - \sum_{j \geq 0} \frac{(-1)^j}{j!} \zeta'(-j)
  \sum_{r=j}^\infty T_1(r,j)(1-z)^r\notag\\
  &= \sum_{r=-1}^\infty
  T_1(r,-1)(1-z)^r L(z) + \sum_{r=-1}^\infty
  T_1(r,-1)(1-z)^r \sum_{j=1}^\infty T_2(j) (1-z)^j\notag\\ & \quad
  + \Gamma'(1) \sum_{r=-1}^\infty T_1(r,-1)(1-z)^r
  + \sum_{r=0}^\infty (1-z)^r \sum_{j=0}^r \frac{(-1)^j}{j!} \zeta'(-j)
  T_1(r,j)\\
  &= \sum_{r=-1}^\infty
  T_1(r,-1)(1-z)^r L(z) + \sum_{r=-1}^\infty
  \sum_{j=1}^\infty T_1(r,-1) T_2(j)(1-z)^{r+j}\\
  & \quad +  \Gamma'(1) \sum_{r=-1}^\infty T_1(r,-1)(1-z)^r +
  \sum_{r=0}^\infty S_4(r)(1-z)^r
\end{align*}
where
\begin{equation}
  \label{eq:5.1.6.3}
  S_4(r) :=  \sum_{j=0}^r \frac{(-1)^j}{j!} \zeta'(-j) T_1(r,j).
\end{equation}
Hence
\begin{align*}
  \Li_{0,1}(z) &\sim \sum_{r=-1}^\infty
  T_1(r,-1)(1-z)^r L(z) +
  \sum_{r=-1}^\infty\sum_{j=r+1}^\infty T_1(r,-1)T_2(j-r)
  (1-z)^j\notag\\
  & \quad + \Gamma'(1) \sum_{r=-1}^\infty
  T_1(r,-1)(1-z)^r + \sum_{r=0}^\infty S_4(r)(1-z)^r\notag\\
  &= \sum_{r=-1}^\infty T_1(r,-1)(1-z)^r L(z) +
  \sum_{j=0}^\infty (1-z)^j \sum_{r=-1}^{j-1}
  T_1(r,-1)T_2(j-r) \\
  & \quad+  \Gamma'(1)\sum_{r=-1}^\infty
  T_1(r,-1)(1-z)^r + \sum_{r=0}^\infty S_4(r)(1-z)^r\notag\\ &=
  \sum_{r=-1}^\infty T_1(r,-1)(1-z)^r L(z) +
  \sum_{r=-1}^\infty S_5(r)(1-z)^r,
\end{align*}
where
\begin{equation}
  \label{eq:5.1.6.4}
  S_5(r) := 
  \begin{cases}
    \Gamma'(1)T_1(-1,1) & r=-1\\ \Gamma'(1)T_1(r,-1) + S_4(r) +
    \displaystyle\sum_{j=-1}^{r-1} T_1(j,-1)T_2(r-j) & r \geq 0.
  \end{cases}
\end{equation}

Now we multiply~\eqref{eq:5.1.2} and~\eqref{eq:5.1.3},
using~\eqref{eq:5.1.2.3} and~\eqref{eq:5.1.6.3}, to get
\begin{align*}
  & \left(\sum_{k=0}^{m-2} S_2(k)(1-z)^k\right)
  \left(\sum_{r=-1}^\infty
    T_1(r,-1)(1-z)^r L(z) + \sum_{r=-1}^\infty
    S_5(r)(1-z)^r\right)\notag\\
  &= \sum_{k=0}^{m-2} \sum_{r=-1}^\infty
  S_2(k)T_1(r,-1)(1-z)^{k+r} L(z) +
  \sum_{k=0}^{m-2} \sum_{r=-1}^\infty S_2(k) S_5(r)
  (1-z)^{k+r}\notag\\
  &= \sum_{k=0}^{m-2} \sum_{r=k-1}^\infty
  S_2(k)T_1(r-k,-1)(1-z)^r L(z) +
  \sum_{k=0}^{m-2}\sum_{r=k-1}^\infty S_2(k)S_5(r-k)(1-z)^r\notag\\
  &= \sum_{r=-1}^\infty (1-z)^r L(z) \sum_{k=0}^{\min{\{r+1,m-2\}}}
  S_2(k)T_1(r-k,-1)\\
  & \quad + \sum_{r=-1}^\infty
  (1-z)^r \sum_{k=0}^{\min{\{r+1,m-2\}}}S_2(k)S_5(r-k)\notag\\
  \label{eq:5.1.7.5}
  &= \sum_{r=-1}^\infty S_6(r)(1-z)^r L(z) + \sum_{r=-1}^\infty
  S_7(r)(1-z)^r
\end{align*}
where
\begin{equation}
  \label{eq:5.1.7.6}
  S_6(r) := \sum_{k=0}^{\min{\{r+1,m-2\}}} S_2(k)T_1(r-k,-1)
  \quad\text{ and }\quad
  S_7(r) := \sum_{k=0}^{\min{\{r+1,m-2\}}} S_2(k)S_5(r-k).
\end{equation}
We are now ready to plug~\eqref{eq:5.1.7.6}, \eqref{eq:5.1.5},
and~\eqref{eq:5.1.6} into the expression for \(T(z)\). This gives us
\begin{align}
  \label{eq:5.1.8.2}
  T(z) &\sim \sum_{r=-1}^\infty
  S_6(r)(1-z)^r L(z) + \sum_{r=-1}^\infty
  S_7(r)(1-z)^r \notag\\
  \quad &- \sum_{k=0}^{m-2} S_1(k)(1-z)^k -
  [\ln{(m-1)!}]\left[(1-z)^{-1} - \sum_{k=0}^{m-2}
  S_3(k)(1-z)^k\right]\notag\\
  &= \sum_{r=-1}^\infty
  S_6(r)(1-z)^r L(z) + \sum_{r=-1}^\infty S_8(r)(1-z)^r,
\end{align}
where
\begin{equation}
  \label{eq:5.1.8.3} S_8(r) :=
  \begin{cases}
    S_7(-1) - \ln{(m-1)!} & r=-1\\ S_7(r) - S_1(r) +
    [\ln{(m-1)!}]S_3(r) & 0 \leq r \leq m-2\\ S_7(r) & r \geq m-1.
  \end{cases}
\end{equation}

\subsection*{Application of the ETT}
\label{sec:application-ett}
We now apply the ETT (Theorem~\ref{theorem:ETT}).  For our toll
function,
\begin{equation*}
  \widehat{T}(z) := T(z) - \sum_{j=0}^{m-1} t_n z^j = T(z).
\end{equation*}
Thus, the generating function of the mean of the shape functional,
call it \(G_1\), is given by
\begin{equation}
  \label{eq:5.1.8.4}
  G_1(z) = T(z) + m! \sum_{j=1}^{m-1}
  \frac{(1-z)^{-\lambda_j}}{\psi'(\lambda_j)} \int_{\zeta=0}^z
  T(\zeta)(1-\zeta)^{\lambda_j-1} \, d\zeta.
\end{equation}
We rewrite~\eqref{eq:5.1.8.2} as
\begin{equation}
  \label{eq:5.1.9.4} T(z) \sim \sum_{r=-2}^\infty S_9(r)
  (1-z)^{r+1}L(z) + \sum_{r=-2}^\infty
  S_{10}(r)(1-z)^{r+1},
\end{equation}
where
\begin{equation}
  \label{eq:5.1.15.5}
  S_9(r) := S_6(r+1) \quad \text{and} \quad S_{10}(r) := S_8(r+1).
\end{equation}
For a fixed \(\lambda \in \mathbb{C}\), the terms
\((1-z)^{\lambda+r}L(z)\) and \((1-z)^{\lambda+r}\) are
both integrable (resp.,\ neither integrable) at \(z=1\) iff \(r \geq
-\lceil\Re(\lambda)\rceil\) (resp.,\ iff \(r <
-\lceil\Re(\lambda)\rceil\)).
    
We therefore write
\begin{subequations}
  \begin{align}
    &\int_{\zeta=0}^z T(\zeta)(1-\zeta)^{\lambda-1}d\zeta\notag\\
    \label{eq:5.1.10.1a}
    &= \sum_{-2 \leq r < -\lceil\Re(\lambda)\rceil} \int_{\zeta=0}^z
    [S_9(r)L(\zeta) +
    S_{10}(r)](1-\zeta)^{\lambda+r} \,d\zeta\\
    \label{eq:5.1.10.1b}
    & \quad + \int_{\zeta=0}^z \biggl[T(\zeta)(1-\zeta)^{\lambda-1}
     - \sum_{-2 \leq r <
      -\lceil\Re(\lambda)\rceil} [S_9(r) L(\zeta) +
    S_{10}(r)](1-\zeta)^{\lambda+r}\biggr] \,d\zeta.
  \end{align}
\end{subequations}
The integrand in~\eqref{eq:5.1.10.1b} is integrable at \(\zeta=1\) and is
asymptotically equivalent to
\[ \sum_{r \geq -\lceil\Re(\lambda)\rceil} [S_9(r) L(\zeta) +
S_{10}(r)](1-\zeta)^{\lambda+r},\] as \(\zeta \to 1\) in the sector
$-\pi + \epsilon < \arg{(1-\zeta)} < \pi - \epsilon$ for any $\epsilon
> 0$.

For~\eqref{eq:5.1.10.1a}, we have the following formulas, valid for any
\(\alpha\) with \(\Re(\alpha) > -1\):
\begin{align*}
  \int_{0}^z (1-\zeta)^\alpha \,d\zeta &=
  -\frac{1}{\alpha+1}[(1-z)^{\alpha+1} - 1],\\
  \intertext{and}
  \int_0^z (1-\zeta)^\alpha L(\zeta) \,d\zeta &=
  -\frac{1}{\alpha+1}(1-z)^{\alpha+1} L(z) -
  \frac{1}{(\alpha+1)^2}[(1-z)^{\alpha+1} - 1].\\
\end{align*}
In the special case that \(m\) is odd and \(\lambda=\lambda_{m-1}=-m\)
and \(r=m-1\), we also need the formulas
\begin{equation*}
  \int_{\zeta=0}^z (1-\zeta)^{-1} \,d\zeta = L(z)
\end{equation*}
and
\begin{equation*}
\int_{\zeta=0}^z (1-\zeta)^{-1} L(\zeta) \,d\zeta =
  \frac{1}{2} L^2(z).
\end{equation*}
We now conclude that~\eqref{eq:5.1.10.1a} equals, except in the special
case described above,
\begin{multline}
  \label{eq:5.1.11.4} \sum_{-2 \leq r < -\lceil\Re(\lambda)\rceil} \Biggl[
  -\frac{S_9(r)}{\lambda+r+1}(1-z)^{\lambda+r+1} L(z)
  - \frac{S_9(r)}{(\lambda+r+1)^2}(1-z)^{\lambda+r+1}\\
  + \frac{S_9(r)}{(\lambda+r+1)^2}
  + \frac{S_{10}(r)}{\lambda+r+1} - \frac{S_{10}(r)}{\lambda+r+1}
  (1-z)^{\lambda+r+1} \Biggr].
\end{multline}
In the special case (\(m\) is odd, \(\lambda=-m\), \(r=m-1\)),
\eqref{eq:5.1.10.1a} equals
\begin{multline}
  \label{eq:5.1.11.5}
  \sum_{r=-2}^{m-2} \Biggl[ -
  \frac{S_9(r)}{\lambda+r+1} (1-z)^{\lambda+r+1} L(z) -
  \frac{S_9(r)}{(\lambda+r+1)^2}(1-z)^{\lambda+r+1}
  + \frac{S_9(r)}{(\lambda+r+1)^2}\\
  + \frac{S_{10}(r)}{\lambda+r+1} -
  \frac{S_{10}(r)}{\lambda+r+1}(1-z)^{\lambda+r+1} \Biggr]
  + \frac{1}{2} S_9(m-1) L^2(z) +
  S_{10}(m-1) L(z).
\end{multline}
  
In the case \(m\) is even,
\begin{align}
  \label{eq:5.1.12.3} &
  \sum_{j=1}^{m-1}
  \frac{(1-z)^{-\lambda_j}}{\psi'(\lambda_j)} \times \eqref{eq:5.1.10.1a}
  \Biggl\rvert_{\lambda=\lambda_j} \notag\\ &=
  \sum_{j=1}^{m-1}\frac{1}{\psi'(\lambda_j)} \sum_{-2 \leq r <
  -\lceil\Re(\lambda_j)\rceil} \Biggl[ -
  \frac{S_9(r)}{\lambda_j+r+1}(1-z)^{r+1} L(z)\notag\\
  & \quad - \left(\frac{S_{10}(r)}{\lambda_j+r+1} +
  \frac{S_9(r)}{(\lambda_j+r+1)^2}\right)(1-z)^{r+1} 
  + \frac{S_9(r)}{(\lambda_j+r+1)^2} +
  \frac{S_{10}(r)}{\lambda_j+r+1}\Biggr]\notag\\ &= \sum_{r=-1}^{m-1}
  [S_{11}(r)(1-z)^r L(z) + S_{12}(r)(1-z)^r] + S_{13},
\end{align}
where
\begin{equation*}
  S_{11}(r) := -\sum_{j : \Re(\lambda_j) \leq -r}
  \frac{S_9(r-1)}{\psi'(\lambda_j)(\lambda_j+r)},
\end{equation*}
\begin{equation*}
S_{12}(r) := -
  \sum_{j: \Re(\lambda_j) \leq -r} \frac{1}{\psi'(\lambda_j)} \left(
    \frac{S_9(r-1)}{(\lambda_j+r)^2} + \frac{S_{10}(r-1)}{\lambda_j +
      r}\right),
\end{equation*}
and
\begin{equation}
  \label{eq:5.1.23.99}
S_{13} := \sum_{r=-1}^{m-1} \sum_{\Re(\lambda_j) \leq -r}
  \frac{1}{\psi'(\lambda_j)} \left[\frac{S_9(r-1)}{(\lambda_j+r)^2} +
    \frac{S_{10}(r-1)}{\lambda_j+r+1}\right].
\end{equation}

If \(m\) is odd, we need to add
\[ \frac{(1-z)^m}{\psi'(-m)} \left(\frac{1}{2}
  S_9(m-1) L^2(z) + S_{10}(m-1) L(z) \right)\] to the right-hand-side
of~\eqref{eq:5.1.12.3} to get the left-hand-size
of~\eqref{eq:5.1.12.3}.

For~\eqref{eq:5.1.10.1b}, we define
\begin{equation}
  \label{eq:5.1.19.5}
J_\lambda(\zeta) := T(\zeta)(1-\zeta)^{\lambda-1} - \sum_{-2 \leq r
  < -\lceil\Re(\lambda)\rceil}
[S_9(r) L(z) + S_{10}(r)](1-\zeta)^{\lambda+r}
\end{equation}
so that~\eqref{eq:5.1.10.1b} equals $\bar{K}(\lambda) -
\displaystyle\int_{\zeta=z}^1 J_\lambda(\zeta) \,d\zeta$, where
\begin{equation}
  \label{eq:5.1.70}
  \bar{K}(\lambda) := \int_{\zeta=0}^1 J_\lambda(\zeta) \,d\zeta  \in
  \mathbb{C}.
\end{equation}
Now we make use of
\begin{equation}
  \label{eq:5.1.13.5} \int_{\zeta=z}^1 J_\lambda(\zeta) \,d\zeta \sim
  \sum_{r \geq -\lceil\Re(\lambda)\rceil} \int_{\zeta=z}^1
  \Bigr[S_9(r) L(\zeta) +
  S_{10}(r)\Bigr](1-\zeta)^{\lambda+r} \,d\zeta.
\end{equation}
Again, exact integration is possible: for all \(\alpha \in
\mathbb{C}\) with \(\Re(\alpha) > -1\),
\begin{equation*}
  \int_{\zeta=z}^1 (1-\zeta)^\alpha \,d\zeta =
  \frac{1}{\alpha+1}(1-z)^{\alpha+1},
\end{equation*}
\begin{equation*}
  \int_{\zeta=z}^1
  (1-\zeta)^\alpha L(\zeta) \,d\zeta =
  \frac{1}{\alpha+1}(1-z)^{\alpha+1} L(z) +
  \frac{1}{(\alpha+1)^2}(1-z)^{\alpha+1}.
\end{equation*}
Therefore
\begin{multline*}
  \int_{\zeta=z}^1 J_\lambda(\zeta) \,d\zeta \sim \sum_{r \geq
    -\lceil\Re(\lambda)\rceil}\Bigl[\frac{S_9(r)}{\lambda+r+1}
  (1-z)^{\lambda+r+1} L(z) \\
  + \frac{S_9(r)}{(\lambda+r+1)^2}
  (1-z)^{\lambda+r+1}
  + \frac{S_{10}(r)}{\lambda+r+1} (1-z)^{\lambda+r+1} \Bigr].
\end{multline*}

We can thus conclude
\begin{align}
  \label{eq:5.1.14.3} &
  \sum_{j=1}^{m-1}
  \frac{(1-z)^{-\lambda_j}}{\psi'(\lambda_j)} \times
  \eqref{eq:5.1.10.1b}\Bigl\rvert_{\lambda=\lambda_j}\notag\\ & \sim
  \sum_{j=1}^{m-1} \frac{\bar{K}(\lambda_j)}{\psi'(\lambda_j)}
  (1-z)^{-\lambda_j}
   - \sum_{j=1}^{m-1}
  \frac{1}{\psi'(\lambda_j)} \sum_{r \geq -\lceil\Re(\lambda)\rceil}
  \Bigl[ \frac{S_9(r)}{\lambda_j+r+1} (1-z)^{r+1} L(z)\notag\\
  &\qquad\qquad\qquad\qquad\qquad\qquad\qquad\qquad
  + \Bigl( \frac{S_9(r)}{(\lambda_j+r+1)^2} +
  \frac{S_{10}(r)}{\lambda_j+r+1}\Bigr)(1-z)^{r+1} \Bigr]\notag\\ &=
  \sum_{j=1}^{m-1} \frac{\bar{K}(\lambda_j)}{\psi'(\lambda_j)}
  (1-z)^{-\lambda_j} + \sum_{r
    \geq -1} [S_{14}(r)(1-z)^r L(z) + S_{15}(r)(1-z)^r]
\end{align}
where
\begin{align*}
  S_{14}(r) &:= - \sum_{j: \Re(\lambda_j) > -r}
  \frac{S_9(r-1)}{(\lambda_j+r)\psi'(\lambda_j)}\\ S_{15}(r) &:= -
  \sum_{j:\Re(\lambda_j) > -r}
  \frac{1}{\psi'(\lambda_j)}\Bigl[\frac{S_9(r-1)}{(\lambda_j+r)^2} +
  \frac{S_{10}(r-1)}{\lambda_j+r} \Bigr].
\end{align*}

\subsection*{Case $m$ even}
\label{sec:case-m-even}

Using~\eqref{eq:5.1.12.3}, \eqref{eq:5.1.14.3},
and~\eqref{eq:5.1.8.2} in~\eqref{eq:5.1.8.4}, we get
\begin{equation}
  \label{eq:5.1.15.1} G_1(z) = \sum_{j=1}^{m-1}
  \frac{K(\lambda_j)}{\psi'(\lambda_j)} (1-z)^{-\lambda_j} + M +
  \sum_{r \geq -1}[ K_r(1-z)^r L(z) + L_r(1-z)^r],
\end{equation}
where 
\begin{equation}
  \label{eq:5.1.15.2} K(\lambda) := m! \bar{K}(\lambda)
  \quad\text{ and }\quad
  M := m!S_{13},
\end{equation}
\begin{align}
  \label{eq:5.1.15.4}
  K_r &:= m![S_{11}(r) + S_{14}(r)] + S_6(r)
  = -m!S_9(r-1)\sum_{j=1}^{m-1} \frac{1}{(\lambda_j+r)\psi'(\lambda_j)} +
  S_6(r),\notag\\
  &= \frac{m!S_9(r-1)}{\rising{(-r)}{m-1} -m!} + S_6(r),
\end{align}
and
\begin{align}
  \label{eq:5.1.26}
  L_r &:= m![S_{12}(r) + S_{15}(r)] + S_8(r)\notag\\
  &= -m!\Biggl[S_9(r-1)\sum_{j=1}^{m-1}
  \frac{1}{(\lambda_j+r)^2\psi'(\lambda_j)} 
  + S_{10}(r-1) \sum_{j=1}^{m-1}
  \frac{1}{(\lambda_j+r)\psi'(\lambda_j)}\Biggr]
  + S_8(r)\notag\\
  &= -m!\Biggl[S_9(r-1)\sum_{j=1}^{m-1}
  \frac{1}{(\lambda_j+r)^2\psi'(\lambda_j)} -
  \frac{S_{10}(r-1)}{\rising{(-r)}{m-1} - m!}\Biggr] + S_8(r).
\end{align}
These constants can be obtained in a closed form using
Identities~\ref{identity:B.1} and~\ref{identity:B.2}.

\subsection*{Case $m$ odd}
\label{case-m-odd}

In this case we have the same result as~\eqref{eq:5.1.15.1} except
that the \(r=m\) term must be replaced by
\begin{multline*}
  m! \biggl\{ \frac{(1-z)^m}{\psi'(-m)} \left[\frac{1}{2}
  S_9(m-1) L^2(z) + S_{10}(m-1) L(z)\right]\\
  + S_{14}(m)(1-z)^m L(z) + S_{15}(m)(1-z)^m\biggr\}
  + S_6(m)(1-z)^m L(z) + S_8(m)(1-z)^m\\
  = N_m(1-z)^m L^2(z) + K_m(1-z)^m L(z) +
  L_m(1-z)^m,
\end{multline*}
where
\begin{equation*}
  N_m  := \frac{m!S_9(m-1)}{2\psi'(-m)} = -
  \frac{S_9(m-1)}{2(H_m-1)},\quad K_m  := m! \Bigl[
  \frac{S_{10}(m)}{\psi'(-m)} + S_{14}(m) \Bigr] + S_6(m),
\end{equation*}
and
\begin{equation*}
L_m  :=  m!S_{15}(m) + S_8(m).
\end{equation*}
(We note that $K_m$ and $L_m$ defined here for $m$ odd are different
[compare~\eqref{eq:5.1.15.4} and~\eqref{eq:5.1.26}] from those defined
for $m$ even.)
The constants can be simplified using Identities~\ref{identity:B.3}
and~\ref{identity:B.4}.  Indeed,
\begin{align*}
  S_{14}(m) &= -\sum_{j:\Re(\lambda_j)>-m}
  \frac{S_9(m-1)}{(\lambda_j+m)\psi'(\lambda_j)} = -S_9(m-1)
  \sum_{j\neq m-1} \frac{1}{(\lambda_j+m)\psi'(\lambda_j)}\\
  & =  -\frac{S_9(m-1)}{2m!} \Bigl[1 - \frac{H_m^{(2)} -
    1}{(H_m-1)^2}\Bigr],
\end{align*}
and
\begin{align*}
  S_{15}(m) &= -S_9(m-1) \sum_{j \neq m-1}
  \frac{1}{(\lambda_j+m)^2\psi'(\lambda_j)}
  -  S_{10}(m-1)\sum_{j \neq m-1} \frac{1}{(\lambda_j+m)
    \psi'(\lambda_j)}\\
  &= -\frac{S_9(m-1)}{m!}\left[ \frac{H_m-1}{12} -
    \frac{1}{3}\frac{H_m^{(3)}-1}{(H_m-1)^2} + \frac{1}{4}
    \frac{(H_m^{(2)} - 1)^2}{(H_m-1)^3}\right]\\
  &\qquad
   - \frac{S_{10}(m-1)}{2 \cdot m!}\left[1 -
    \frac{H_m^{(2)} - 1}{(H_m-1)^2}\right].
\end{align*}

\subsection*{Some observations}
The following facts can be established for the constants.
\begin{remark}
  \label{remark:5.1.1}
  \(K_r = 0\) for \(0 \leq r \leq m-2\).  
\end{remark}
\begin{proof}
  \begin{equation*}
    K_r = \frac{m!S_9(r-1)}{\rising{(-r)}{m-1} - m!} + S_6(r)
    = \frac{m! S_6(r)}{\rising{(-r)}{m-1} - m!} + S_6(r) =
    S_6(r) \frac{\rising{(-r)}{m-1}}{\rising{(-r)}{m-1} -
      m!},
  \end{equation*} which contains the factor 0 if \(0 \leq r \leq m-2.\)
\end{proof}
\begin{remark}
  \label{remark:5.1.2}
  \(K_{-1} = -1\).
\end{remark}
\begin{proof}
  \begin{align*}
    K_{-1}                          & = S_6(-1)\frac{(m-1)!}{(m-1)!-m!}
    = \frac{S_6(-1)}{1-m} = \frac{(m-1)T_1(-1,-1)}{1-m}\\
                                    & = -T_1(-1,-1) = -[z^{-1}]L(z) = -1.
  \end{align*}
\end{proof}
\begin{remark}
  \label{remark:5.1.3}
  \(K_{m-1} = 0\) for all odd values of \(m\).
\end{remark}
\begin{proof}
  When \(m\) is odd,
  \begin{equation*}
    K_{m-1} = \frac{S_6(m-1) (-1)^{m-1} (m-1)!}{(-1)^{m-1} (m-1)! -
      m!} = -\frac{S_6(m-1)}{m-1}
  \end{equation*}
  Now \(S_6(m-1)\) is the coefficient of
  \((1-z)^{m-1}L(z) \) in the product \(\eqref{eq:5.1.2} \times
  \eqref{eq:5.1.3}\), i.e.,
  \begin{align*}
    S_6(m-1)                         & = [(1-z)^{m-1}L(z)]\left(
      \eqref{eq:5.1.2} \times \eqref{eq:5.1.3}\right) \\
    &= [(1-z)^{m-1} L(z)]
    \left(
      \frac{1-z^{m-1}}{1-z} (-\ln z)^{-1} L(z)
    \right)\\
       & = [(1-z)^{m-1}] \left( \frac{1-z^{m-1}}{1-z}
       (-\ln{z})^{-1}\right)       \\ 
     & = [z^{m-1}] \left( \frac{1 - (1-z)^{m-1}}{z L(z)}
         \right)
     = [z^m] \frac{1 - (1-z)^{m-1}}{L(z)} = 0.
  \end{align*}
  To see the last equality, observe that
  \[ \frac{1 - (1-z)^{m-1}}{L(z)} = \int_{x=0}^{m-1} (1-z)^x \,dx \]
  and
  \begin{equation*}
    [z^m] \int_{x=0}^{m-1} (1-z)^x \,dx = [z^m]\int_{x=0}^{m-1}
    \sum_{n=0}^\infty \frac{\rising{(-x)}{n}}{n!} z^n \,dx =
    \frac{1}{m!}\int_{x=0}^{m-1} \rising{(-x)}{m} \,dx.
  \end{equation*}
  But the last integral vanishes because its
  integrand is antisymmetric about the midpoint
  of the interval of integration:
  \[ \rising{-(m-1-x)}{m} = (-1)^m(m-1-x)_m =
  (-1)^m\rising{-x}{m} = -\rising{(-x)}{m},\]
  the last equality following from the assumption that \(m\) is odd.
\end{proof}
\begin{remark}
  \label{remark:5.1.5}
  \[
  L_{-1} = \frac{\ln{(m-1)!}}{m-1} - \frac{m}{m-1}H_{m-1} +
  \gamma.
  \]
\end{remark}
\begin{proof}
  \begin{align*}
    L_{-1} &= m!\Bigl[ S_9(-2)\Bigl(-
    \frac{H_{m-1}}{(m-1)!(m-1)^2}\Bigr) + S_{10}(-2) \Bigl(
    \frac{1}{(m-1)! - m!}\Bigr)\Bigr] + S_8(-1)\\
    &=  m!\Bigl[-S_6(-1) \frac{H_{m-1}}{(m-1)!(m-1)^2} - S_8(-1) \Bigl(
    \frac{1}{(m-1)!(m-1)}\Bigr)\Bigr] + S_8(-1).
  \end{align*}
  Now \(S_6(-1) = m-1\) as noted in the proof of Remark~\ref{remark:5.1.2},
  and \(S_8(-1) = -(m-1)\gamma - \ln{(m-1)!}\). Plugging
  these in the display above, we get the desired result.
\end{proof}

\subsection*{The complete asymptotic expansion}
\label{sec:compl-asympt-expans}

The computations in this section have thus established that the
generating function of the mean of the shape functional for \( m
\)-ary search trees under the random permutation model, denoted by \(
G_1(z)\), has the complete asymptotic expansion
\begin{multline}
  \label{eq:5.1.100}
  G_1(z) \sim \sum_{j=1}^{m-1} \frac{K(\lambda_j)}{\psi'(\lambda_j)}
  (1-z)^{-\lambda_j} + M - (1-z)^{-1} L(z) +
  L_{-1}(1-z)^{-1}\\
  + \sum_{r=0}^{m-2} L_r(1-z)^r
  + \sum_{r \geq m-1} [ K_r(1-z)^r L(z) + L_r (1-z)^{r} ]
\end{multline}
when \( m \) is even and
\begin{multline*}
  G_1(z) \sim \sum_{j=1}^{m-1} \frac{K(\lambda_j)}{\psi'(\lambda_j)}
  (1-z)^{-\lambda_j} + M - (1-z)^{-1} L(z) +
  L_{-1}(1-z)^{-1}
  + \sum_{r=0}^{m-1} L_r(1-z)^r \\ + N_m(1-z)^{m} L^2(z)  +
  K_m(1-z)^m L(z) 
  + L_{m} (1-z)^{m} + \sum_{r \geq m+1} [ K_r(1-z)^r L(z)
  + L_r (1-z)^{r} ]
\end{multline*}
when \( m \) is odd.

Using singularity analysis (Theorem~\ref{theorem:transfer}), the
complete asymptotic expansion can be translated on a term-by-term
basis to a complete asymptotic expansion of the mean itself, denoted
by \( \mu_n \).  Thus
\begin{equation}
  \label{eq:5.1.60}
  \mu_n \sim \sum_{j=1}^{m-1} \frac{K(\lambda_j)}{\psi'(\lambda_j)}
  \frac{\rising{\lambda_j}{n}}{n!} - H_n + L_{-1} + \sum_{r \geq m-1}
  K_r \Delta^r\left((\tfrac{1}{n})\right)
\end{equation}
when \( m \) is even and
\begin{equation*}
  \mu_n \sim \sum_{j=1}^{m-1} \frac{K(\lambda_j)}{\psi'(\lambda_j)}
  \frac{\rising{\lambda_j}{n}}{n!} - H_n + L_{-1} + N_m \Delta^m
  \left( \left( \tfrac{2H_{n-1}}{n} \right) \right) + K_m
  \Delta^m\left((\tfrac{1}{n})\right) + \sum_{r \geq m+1}
  K_r \Delta^r\left((\tfrac{1}{n})\right)
\end{equation*}
when \( m \) is odd, where \( \Delta^r(s) \) denotes the \( r \)th-order
difference of a sequence $s = (s_n)$.  Note that $\Delta^r((1/n)) =
\Theta(n^{-(r+1)})$ and $\Delta^m((2H_{n - 1}/n)) =
\Theta(n^{-(m+1)}\log n)$.

\subsection*{The binary search tree}
\label{sec:binary-search-tree}

In preparation for the computation of the variance of the shape
functional for the binary search trees we state a special case of the
complete asymptotic expansion derived above.  When \( m=2 \), the toll
function $t_n$ equals $\ln{n}$ and the
indicial polynomial $\psi(\lambda)$ equals $\lambda - 2$, with root
$\lambda_1 = 2$.  Thus~\eqref{eq:5.1.60},
Remark~\ref{remark:5.1.5}, and the well-known expansion for the
harmonic numbers~\cite[1.2.7-(3)]{knuth97} immediately yield
\begin{equation}
  \label{eq:5.1.61}
  \mu_n = K(2) \frac{\rising{2}{n}}{n!} - H_n - 2 + \gamma + O(n^{-2})
  = K(2)(n+1) - \ln{n} - 2 - \frac{1}{2n} + O(n^{-2}),
\end{equation}
where, by~\eqref{eq:5.1.15.2},~\eqref{eq:5.1.70},
and~\eqref{eq:5.1.19.5},
\begin{equation*}
  K(2) = 2 \bar{K}(2) = 2\int_{\zeta=0}^1 T(\zeta) (1-\zeta)\,d\zeta =
  2 \int_{\zeta=0}^1 \sum_{k=1}^\infty (\ln{k}) \zeta^k
  (1-\zeta)\,d\zeta = 2 \sum_{k=1}^\infty \frac{\ln{k}}{(k+1)(k+2)}.
\end{equation*}

\section{The variance for binary search trees}
\label{sec:bst-vari-shape-funct}

Let \( X_n \) denote the mean of the shape functional for binary
search trees under the random permutation model.  In the sequel it
will be convenient to work with the approximately centered random
variable \( \widetilde{X}_n \)
defined as
\begin{equation*}
  \widetilde{X}_n := X_n - C_1(n+1) \quad \text{for \( n \geq 0 \)},
\end{equation*}
where $C_1 := K(2)$, and to define
\begin{equation*}
  \tilde{\mu}_n(k) := \E{\widetilde{X}_n^k} \quad \text{for \( k
  \geq 1 \)},
\end{equation*}
with \( \tilde{\mu}_n := 0 \).
Then, following the arguments leading to~\eqref{eq:2.3}
and~\eqref{eq:2.4}, we have
\[
\tilde{\mu}_n(2) = \frac{2}{n} \sum_{j=0}^{n-1} \tilde{\mu}_j(2) +
r_n(2),
\]
where
\begin{align*}
 r_n(2)
  &:=  \sum_{\substack{k_1,k_2,k_3 \geq 0\\k_1+k_2+k_3=2\\k_1 < 2,\ k_2 <
      2}} \binom{2}{k_1,k_2,k_3}(\ln{n})^{k_3} \frac{1}{n}
  \sum_{\substack{j_1+j_2=n-1\\ j_1,j_2 \geq 0}}
    \tilde{\mu}_{j_1}(k_1)\tilde{\mu}_{j_2}(k_2)\\
&= (\ln{n})^2 + 4(\ln{n})\frac{1}{n}\sum_{j=0}^{n-1}
\tilde{\mu}_j(1)
+ \frac{2}{n}\sum_{j=0}^{n-1}\tilde{\mu}_j(1)\tilde{\mu}_{n-1-j}(1)\\
&= (\ln{n})^2 + 4(\ln{n})\frac{1}{n} \sum_{j=0}^{n-1}(\mu_j-C_1(j+1))\\
 &\qquad+ \frac{2}{n}
  \sum_{j=0}^{n-1}(\mu_j-C_1(j+1))(\mu_{n-1-j}-C_1(n-j)),
\end{align*}
and $r_0(2) = \E{\widetilde{X}_0^2} = C_1^2.$  Note that knowing \(
\tilde{\mu}_n(1) \) and \( \tilde{\mu}_n(2) \) is enough to compute
the variance of the random variable \( X_n \), call it \( \sigma_n^2
\), since
\begin{equation}
  \label{eq:5.2.70}
  \sigma_n^2 = \Var{\widetilde{X}_n} = \tilde{\mu}_n(2) - \tilde{\mu}_n(1)^2.
\end{equation}

Our proposal is to use the ETT~(Theorem~\ref{theorem:ETT}) again.  To
that end, let $B(z)$ be the generating function of $r_n(2)$.  Then
using the expression for \( r_n(2) \) above we get
\begin{multline}
  \label{eq:5.2.V3.1}
  B(z) = C_1^2 + \Li_{0,2}(z)+ 4\Li_{0,1}(z)\odot \Bigl[ \int_0^z
  (1-\zeta)^{-1}\bigl(G_1(\zeta) - C_1(1-\zeta)^{-2}\bigr)\,d\zeta\Bigr]\\
  + 2\int_0^z [G_1(\zeta) - C_1(1-\zeta)^{-2}]^2\,d\zeta.
\end{multline}
Our goal is to get each term in this sum to $O(|1-z|^{-\epsilon})$ for
arbitrarily small $\epsilon > 0$.

Using Theorem~\ref{theorem:Li-sing} we get
\[
\Li_{\alpha,0}(z) \sim \Gamma(1-\alpha)t^{\alpha-1} + \sum_{j \geq 0}
\frac{(-1)^j}{j!}\zeta(\alpha-j)t^j,
\]
\[
\Li_{\alpha,1}(z) \sim \Gamma'(1-\alpha)t^{\alpha-1} -
\Gamma(1-\alpha)t^{\alpha-1}\ln{t} - \sum_{j \geq 0}
\frac{(-1)^j}{j!} \zeta'(\alpha-j)t^j,
\]
and
\begin{multline*}
\Li_{\alpha,2}(z) \sim \Gamma''(1-\alpha)t^{\alpha-1} -
2\Gamma'(1-\alpha)t^{\alpha-1}\ln{t} \\
+
\Gamma(1-\alpha)t^{\alpha-1}(\ln{t})^2 + \sum_{j\geq 0}
\frac{(-1)^j}{j!} \zeta''(\alpha-j)t^j,
\end{multline*}
where
\begin{equation*}
  t = -\ln{z} = (1-z) + \frac{(1-z)^2}{2} + O(|1-z|^3).
\end{equation*}
Hence
\[
t^{-1} = (1-z)^{-1} - \frac12 + O(|1-z|)
\]
and [recalling $L(z) := \ln((1-z)^{-1})$]
\[
\ln{t} = -L(z) + \frac12 (1-z) + O(|1-z|^2).
\]
This gives
\begin{equation*}
  (\ln{t})^2  = L^2(z) - (1-z)L(z) + O(|(1-z)^2L(z)|),
\end{equation*}
so that
\[
t^{-1}(\ln{t})^2 = (1-z)^{-1}L^2(z) - \frac{1}{2}L^2(z) - L(z) +
  O(|(1-z)L(z)|)
\]
and
\[
t^{-1}\ln{t} = -(1-z)^{-1}L(z) + \frac{1}{2}L(z) + O(1).
\]
Combining these expansions we get
\begin{equation}
  \label{eq:5.2.V5.4}
\Li_{0,1}(z) = (1-z)^{-1}L(z) + \Gamma'(1)(1-z)^{-1} - \frac{1}{2}L(z)
+ O(1)
\end{equation}
and
\begin{align}
  \label{eq:5.2.V5.7}
  \Li_{0,2}(z)
  &= (1-z)^{-1}L^2(z) - \frac{1}{2}L^2(z) - L(z)\notag\\
  & \quad - 2\Gamma'(1)\bigl[-(1-z)^{-1}L(z)
  + \frac12L(z)\bigr]
  + \Gamma''(1)(1-z)^{-1} + O(1)\notag\\
  &= (1-z)^{-1}L^2(z) + 2\Gamma'(1)(1-z)^{-1}L(z)\notag\\
  & \quad +  \Gamma''(1)(1-z)^{-1}
  - \frac{1}{2}L^2(z) - (\Gamma'(1)+1)L(z) + O(1).
\end{align}

We also need an asymptotic expansion for the generating function \(
G_1(z) - C_1(1-z)^{-2}\) of $\tilde{\mu}_n(1)$.  This is available to
us using~\eqref{eq:5.1.100}.  Indeed, using $M=0$
from~\eqref{eq:5.1.15.2}, $L_0 = 1/4$ from~\eqref{eq:5.1.26},
$L_{-1}=-(2-\gamma)$ from Remark~\ref{remark:5.1.5}, and $K_1 = -1/36$
from~\eqref{eq:5.1.15.4} [and the
definitions~\eqref{eq:5.1.15.5},~\eqref{eq:5.1.7.6},~\eqref{eq:5.1.2.3},~\eqref{eq:5.1.6.5},~\eqref{eq:5.1.23.99}],
we get
\begin{equation*}
G_1(z) - C_1(1-z)^{-2} = -(1-z)^{-1}L(z) +
L_{-1}(1-z)^{-1}
+ \frac{1}{4} - \frac{1}{36}(1-z)L(z) +O(|1-z|).
\end{equation*}
Thus, invoking Theorem~\ref{theorem:2B}, we get
\begin{align*}
  &\int_0^z (1-\zeta)^{-1}[G_1(\zeta) - C_1(1-\zeta)^{-2}]\,d\zeta\notag\\
  &= \int_0^z [-(1-\zeta)^{-2}L(\zeta) + L_{-1}(1-\zeta)^{-2}
  + \frac{1}{4}(1-\zeta)^{-1}
  - \frac{1}{36}L(\zeta) + O(1)]\,d\zeta \notag\\
  &= -[(1-z)^{-1}L(z) - (1-z)^{-1}] + L_{-1}(1-z)^{-1}
  + \frac{1}{4}L(z) + O(1) \notag\\
&= -(1-z)^{-1}L(z) + (L_{-1}+1)(1-z)^{-1} + \frac14L(z) + O(1),
\end{align*}
and
\begin{multline}
  \label{eq:5.2.2}
  [G_1(z) - C_1(1-z)^{-2}]^2
  = (1-z)^{-2}L^2(z) - 2L_{-1}(1-z)^{-2}L(z) + L_{-1}^2(1-z)^{-2}\\
  \qquad -\frac12(1-z)^{-1}L(z) + \frac{1}{2}L_{-1}(1-z)^{-1} + O(|L^2(z)|).
\end{multline}
We now integrate~\eqref{eq:5.2.2} term by term.  Integration by parts
establishes
\begin{multline*}
\int_0^z (1-\zeta)^rL^2(\zeta)\,d\zeta =\\
\begin{cases}
  -\displaystyle\frac{1}{r+1}(1-z)^{r+1}L^2(z) -
   \displaystyle\frac{2}{(r+1)^2} (1-z)^{r+1}L(z) - 
   \displaystyle\frac{2}{(r+1)^3}[(1-z)^{r+1} - 1] & (r \ne -1)\\
   \frac{1}{3}L^3(z) & (r = -1)
\end{cases}
\end{multline*}
and
\[
\int_0^z (1-\zeta)^rL(\zeta)\,d\zeta =
\begin{cases}
  -\displaystyle\frac{1}{r+1}(1-z)^{r+1}L(z) -
   \displaystyle\frac{1}{(r+1)^2}(1-z)^{r+1} + 
   \displaystyle\frac{1}{(r+1)^2} & (r \ne -1)\\
   \displaystyle\frac{1}{2} L^2(z) & (r = -1).
\end{cases}
\]
By Theorem~\ref{theorem:2B} the integral
of the remainder in~\eqref{eq:5.2.2} is $O(1) = O(|1-z|^{-\epsilon})$. Hence
\begin{align}
  &\int_0^z [G_1(\zeta) - C_1(1-\zeta)^{-2}]^2\,d\zeta \notag\\
  &= (1-z)^{-1}L^2(z) - 2(1-z)^{-1}L(z) + 2(1-z)^{-1}\notag\\
   & \quad - 2L_{-1}[(1-z)^{-1}L(z) - (1-z)^{-1}]
  + L_{-1}^2(1-z)^{-1} + O(|1-z|^{-\epsilon})\notag\\
  &= (1-z)^{-1}L^2(z) - 2(L_{-1}+1)(1-z)^{-1}L(z)
  \label{eq:5.2.V9.1}
  + (2 + 2L_{-1} + L_{-1}^2)(1-z)^{-1} + O(|1-z|^{-\epsilon}).
\end{align}

We now have an asymptotic expansion for each term in the
sum~\eqref{eq:5.2.V3.1} but the one involving the Hadamard product.
We make use of the following facts.
\begin{enumerate}
\item $(1-z)^{-1}$ is the identity element for Hadamard products,
  i.e., \( (1-z)^{-1} \odot f(z) = f(z) \).
\item $L(z) = O(|1-z|^{-\epsilon})$ for any $\epsilon > 0$, so that
  $(1-z)^{-1}L(z) = O(|1-z|^{-1-\epsilon})$. Hence, by
  Theorem~\ref{thm:hadamard}, $[(1-z)^{-1}L(z)] \odot O(1) =
  O(|1-z|^{-\epsilon})$.
\item Similarly,
  \begin{equation*}
    [(1-z)^{-1}L(z)] \odot L(z) = O(|1-z|^{-\epsilon}),\qquad
    L(z) \odot L(z) = O(1),
  \end{equation*}
  \begin{equation*}
    L(z) \odot O(1) = O(1), \qquad
    O(1) \odot O(1) = O(1).
  \end{equation*}
\item The following lemma will also be used.
  \begin{lemma}
    \label{lem:V1}
    \(
    [(1-z)^{-1}L(z)] \odot [(1-z)^{-1}L(z)] = (1-z)^{-1}L^2(z) +
    (1-z)^{-1}\Li_{2,0}(z).
    \)
  \end{lemma}
\begin{proof}
  Let $H_n$ denote the $n$th harmonic number.  Then
  \begin{align*}
    &[(1-z)^{-1}L(z)] \odot [(1-z)^{-1}L(z)]
    = \sum_{n \geq 1} H_n^2 z^n
    = \sum_{n \geq 1} \left(\sum_{k=1}^n \frac{1}k \right)^2z^n\\
    &= \sum_{n \geq 1} \left(\sum_{k=1}^n \frac1{k^2}\right)z^n + 2
    \sum_{n \geq 1} \Biggl[ \sum_{k=2}^n \frac1k \sum_{j=1}^{k-1}
    \frac1j\Biggr] z^n
    = \sum_{k \geq 1} \frac{1}{k^2} \Biggl(\sum_{n \geq k} z^n\Biggr) +
    2\sum_{n \geq 1} \Biggl[ \sum_{k=2}^n \frac{H_{k-1}}{k}\Biggr] z^n\\
    &= \sum_{k \geq 1} \frac{z^k}{k^2} (1-z)^{-1} + 2\sum_{n \geq 1}
    \Biggl[\sum_{k=2}^n \frac{H_{k-1}}{k}\Biggr] z^n.
  \end{align*}
  The first sum is $(1-z)^{-1}\Li_{2,0}(z)$.  For the second sum,
  observe that
  \begin{align*}
    \sum_{n \geq 1}\Biggl[ \sum_{k=2}^n \frac{H_{k-1}}{k}\Biggr] z^n 
    &= \sum_{k \geq 2} \frac{H_{k-1}}k \sum_{n \geq k}z^n = (1-z)^{-1}
    \sum_{k \geq 1} \frac{H_{k}}{k+1} z^{k+1}\\
    &= (1-z)^{-1}\int_0^z(1-\zeta)^{-1}L(\zeta)\,d\zeta =
    \frac{1}{2}(1-z)^{-1}L^2(z), 
  \end{align*}
  and the claim follows.
\end{proof}
\end{enumerate}

In light of these facts we can now conclude that
\begin{align}
  \label{eq:5.2.3}
  &\Li_{0,1}(z) \odot \int_0^z(1-\zeta)^{-1}[G_1(\zeta) -
  C_1(1-\zeta)^{-2}]\,d\zeta\\
  &= [(1-z)^{-1}L(z) + \Gamma'(1)(1-z)^{-1} - \frac12L(z) + O(1)]\notag\\
  &\qquad\odot [-(1-z)^{-1}L(z) + (L_{-1}+1)(1-z)^{-1} +\frac14 L(z) +
  O(1)]\notag\\ 
  &=-\bigl[(1-z)^{-1}L^2(z) + (1-z)^{-1}\Li_{2,0}(z)\bigr] +
  (L_{-1}+1)(1-z)^{-1}L(z) \notag\\
  &\qquad - \Gamma'(1)(1-z)^{-1}L(z)
  + \Gamma'(1)(L_{-1}+1)(1-z)^{-1} + O(|1-z|^{-\epsilon}). \notag
\end{align}
To get an asymptotic expansion for \( \Li_{2,0} \) we
use~\eqref{eq:4.1.posintalpha}:
\begin{equation*}
  \Li_{2,0}(z) = t(\ln{t} - H_1) + \frac{\pi^2}{6} + O(|1-z|^2) =
  \frac{\pi^2}6 - (1-z)L(z) + O(|1-z|).
\end{equation*}
Hence we see that the Hadamard product~\eqref{eq:5.2.3} has the asymptotic
expansion
\begin{align}
  &-(1-z)^{-1}L^2(z) - (1-z)^{-1}\Bigl(\frac{\pi^2}{6} - (1-z)L(z) +
  O(|1-z|)\Bigr)\notag\\
  &\quad + (L_{-1}  + 1 - \Gamma'(1))(1-z)^{-1}L(z) +
  \Gamma'(1)(L_{-1} + 1)(1-z)^{-1} + O(|1-z|^{-\epsilon})\notag\\
  &= -(1-z)^{-1}L^2(z) + (L_{-1} + 1-\Gamma'(1))(1-z)^{-1}L(z)\notag\\
  &\quad + \Bigl[\Gamma'(1)(L_{-1} + 1)-\frac{\pi^2}6\Bigr](1-z)^{-1} +
  O(|1-z|^{-\epsilon}).
  \label{V10.5}
\end{align}

Using~\eqref{eq:5.2.V5.7},~\eqref{eq:5.2.V9.1}, and~\eqref{V10.5}
in~\eqref{eq:5.2.V3.1}, we get
\begin{align}
  &B(z)\notag\\ &= C_1^2 + (1-z)^{-1}L^2(z) + 2\Gamma'(1)(1-z)^{-1}L(z) +
  \Gamma''(1)(1-z)^{-1}\notag\\
  &+ 4\Bigl[-(1-z)^{-1}L^2(z) +
  (L_{-1}+1-\Gamma'(1))(1-z)^{-1}L(z)
  +\bigl[\Gamma'(1)(L_{-1}+1)-\frac{\pi^2}{6}\bigr](1-z)^{-1}
  \notag\Bigr]\\ 
  &+ 2\Bigl[(1-z)^{-1}L^2(z) - 2(L_{-1}+1)(1-z)^{-1}L(z)
   + (2+2L_{-1}+ L_{-1}^2)(1-z)^{-1}\Bigr] +
  O(|1-z|^{-\epsilon})\notag\\ 
  \label{eq:5.2.V12.3}
  &= -(1-z)^{-1}L^2(z) + V_1(1-z)^{-1}L(z) + V_2(1-z)^{-1} +
  O(|1-z|^{-\epsilon}), 
\end{align}
where
\begin{equation}
  \label{eq:5.2.defV1}
  V_1 := -2\Gamma'(1) = 2\gamma
\end{equation}
and
\begin{equation}
  \label{eq:5.2.defV2}
  V_2 := \Gamma''(1) + 4\Bigl( \Gamma'(1)(L_{-1}+1) - \frac{\pi^2}6\Bigr) +
  2(2+2L_{-1}  + L_{-1}^2) = 4 - \gamma^2 - \frac{\pi^2}2.
\end{equation}

Let $\widehat{B}(z) := B(z) - r_0(2) = B(z) - C_1^2$.  By direct
application of the ETT, $A(z)$, defined to be the generating function
of $\bigl(\tilde{\mu}_n(2)\bigr)$, is given by
\[
A(z) = c_1(1-z)^{-2} + \widehat{B}(z) + 2(1-z)^{-2}\int_0^z
\widehat{B}(\zeta)(1-\zeta)\,d\zeta, 
\]
where $c_1 = r_0(2) = C_1^2$.  The integral here can be broken up as
\[
\int_0^z \widehat{B}(\zeta)(1-\zeta)\,d\zeta = V - \int_z^1
\widehat{B}(\zeta)(1-\zeta)\,d\zeta,
\]
where the first term on the right is the constant
\begin{equation*}
V := \int_0^1 \widehat{B}(\zeta)(1-\zeta)\,d\zeta = \int_0^1 \left(
  \sum_{k=1}^\infty r_k(2)\zeta^k\right) (1-\zeta)\,d\zeta
= \sum_{k=1}^\infty \frac{r_k(2)}{(k+1)(k+2)}.
\end{equation*}
For the second term we use the following consequences of
integration by parts:
\[
\int_z^1 L^2(\zeta)\,d\zeta = (1-z)L^2(z) + 2(1-z)L(z) + 2(1-z),
\]
\[
\int_z^1 L(\zeta)\,d\zeta = (1-z)L(z) + (1-z)
\]
and the crude remainder estimate
\[
\int_z^1 O(|1-\zeta|^{1-\epsilon})\,d\zeta = O(|1-z|^{2-\epsilon}).
\]
Thus
\begin{align*}
  &(1-z)^{-2}\int_z^1 \widehat{B}(\zeta)(1-\zeta)\,d\zeta\\
  &= (1-z)^{-2}\bigl[ -(1-z)L^2(z) - 2(1-z)L(z) - 2(1-z) + V_1(1-z)L(z)\\
 & \qquad\qquad + V_1(1-z) + V_2(1-z) + O(|1-z|^{2-\epsilon})\bigr]\\
  &= -(1-z)^{-1}L^2(z) + (V_1-2)(1-z)^{-1}L(z)  +
  (V_1+V_2-2)(1-z)^{-1}  +  O(|1-z|^{-\epsilon}), 
\end{align*}
and
\begin{align}
  A(z) &= C_1^2(1-z)^{-2} - (1-z)^{-1}L^2(z) + V_1(1-z)^{-1}L(z) + V_2(1-z)^{-1} +
  O(|1-z|^{-\epsilon}) \notag\\
  &\quad + 2\Bigl[ V(1-z)^{-2} + (1-z)^{-1}L^2(z) -
  (V_1-2)(1-z)^{-1}L(z)\notag \\ 
  &\qquad - (V_1+V_2-2)(1-z)^{-1} + O(|1-z|^{-\epsilon})\notag\Bigr]\\
  &= (C_1^2+2V)(1-z)^{-2} + (1-z)^{-1}L^2(z) + (V_1 -
  2(V_1-2))(1-z)^{-1}L(z) \notag\\
  & \qquad + (V_2 - 2(V_1+V_2-2))(1-z)^{-1} + O(|1-z|^{-\epsilon})\notag\\
  &= (C_1^2+2V)(1-z)^{-2} + (1-z)^{-1}L^2(z) + (4-V_1)(1-z)^{-1}L(z)\notag\\
  \label{eq:5.2.V16.3}
  &\qquad + (4-2V_1-V_2)(1-z)^{-1} + O(|1-z|^{-\epsilon}).
\end{align}

We can now transfer this asymptotic expansion for the generating
function to one for its Taylor coefficients, invoking
Theorem~\ref{theorem:transfer} for the remainder estimate.  We know that
\[
[z^n](1-z)^{-1}L(z) = H_n = \ln{n} + \gamma + O(n^{-1})
\]
and, from the proof of Lemma~\ref{lem:V1},
\begin{align*}
[z^n]\bigl[(1-z)^{-1} L^2(z)\bigr] &= H_n^2 - \sum_{k=1}^n \frac{1}{k^2} =
\Bigl(\ln{n}+\gamma+O(n^{-1})\Bigr)^2 - \sum_{k=1}^\infty
\frac{1}{k^2} + O(n^{-1})\\
&= (\ln{n})^2 + 2\gamma\ln{n} + \biggl(\gamma^2 - \frac{\pi^2}6\biggr) +
O\biggl(\frac{\ln{n}}{n}\biggr).
\end{align*}
So, by singularity analysis (Theorem~\ref{theorem:transfer}),
\begin{multline}
  \label{eq:5.2.V17.3}
  \tilde{\mu}_n(2) = (C_1^2+2V)(n+1) + (\ln{n})^2 + 2\gamma\ln{n} +
  \biggl(\gamma^2 - \frac{\pi^2}6 \biggr) +  
  (4-V_1)(\ln{n} + \gamma)\\  + (4-2V_1 - V_2) + O(n^{-1+\epsilon}).
\end{multline}

In order to obtain the variance using~\eqref{eq:5.2.70} we also need
[recall~\eqref{eq:5.1.61}]
\begin{equation*}
  \tilde{\mu}_n(1)^2 = (\ln{n})^2 + 4\ln{n} + 4 + O(n^{-1+\epsilon}).
\end{equation*}
Hence, using~\eqref{eq:5.2.defV1} and~\eqref{eq:5.2.defV2} we obtain
an expansion for the variance:
\begin{align}
  \sigma_n^2 &= (C_1^2+2V)(n+1) + (2\gamma-V_1)\ln{n}
  + \biggl[\Bigl(\gamma^2 -
  \frac{\pi^2}6\Bigr) + (4-V_1)\gamma  - 2V_1 - V_2 \biggr] +
  O(n^{-1+\epsilon})\notag\\
  \label{eq:5.2.V19.1}
  &= (C_1^2+2V)(n+1) 
  - \Bigl(4 - \frac{\pi^2}3\Bigr) +  O(n^{-1+\epsilon}).
\end{align}
Numerical calculations using \texttt{Mathematica} suggest that the
next term in the expansion is $-\tfrac1{2n}$.  One could possibly
obtain this by carrying out the program of this section, keeping track
of more terms in the asymptotic expansions.


\chapter{Random permutation model: limiting distributions}
\label{cha:limit-distr-m}

In this chapter we will apply the Asymptotic Transfer Theorem
(Theorem~\ref{theorem:att}) of Chapter~\ref{cha:transfer-theorems} to
derive limiting distributions for additive functionals defined
on \(m\)-ary search trees under the random permutation model.  Our
main results are
Theorems~\ref{theorem:6.1.5.1},~\ref{theorem:clt5.4},~\ref{thm:F2},
and~\ref{thm:F3}.

Throughout \(\mathcal{N}(\mu,\sigma^2)\) denotes the normal
distribution with mean \(\mu\) and variance \(\sigma^2\).

\section{Small toll functions}
\label{sec:small-toll-functions-1}

We begin with two central limit theorems
(Theorems~\ref{theorem:6.1.5.1} and~\ref{theorem:clt5.4}) for
``small'' toll functions.  Throughout, we write \( \sum_{\mathbf{j}}
\) as shorthand for the sum over \(m\)-tuples \( (j_1,\ldots,j_m) \)
of nonnegative integers summing to \( n-(m-1) \).
\begin{theorem}
  \label{theorem:6.1.5.1}
  If \( 2 \leq m \leq 26 \) and the real toll sequence \( (t_n) \)
  satisfies
  \begin{equation}
    \label{eq:6.1.5.1}
    \text{(a) } t_n = o(\sqrt{n}) \quad\text{ \textnormal{and} }\quad
    \text{(b) } \sum^{\infty} n^{-1} \max_{n^\delta \leq k \leq n}
    \frac{t_k^2}{k} < \infty \quad \text{\textnormal{for some} \( 0 <
      \delta < 1 \)}, 
  \end{equation}
  then the mean \( \mu_n \) and the variance \( \sigma_n^2 \) of the
  corresponding additive functional \( X_n \) on \(m\)-ary search
  trees with the random permutation model satisfy, respectively,
  \begin{equation}
    \label{eq:6.1.5.2}
    \mu_n = \frac{K_1}{H_m-1}n + o(\sqrt{n}) =: \mu n + o(\sqrt{n}),
  \end{equation}
  with \( K_1 \) defined as
  \begin{equation}
    \label{eq:6K1}
    K_1 := \sum_{j=0}^\infty \frac{t_j}{(j+1)(j+2)},
  \end{equation}
  and
  \begin{equation}
    \label{eq:6.1.5.3}
    \sigma_n^2 = \sigma^2 n + o(n), \qquad\text{ where }\qquad
    \sigma^2 := \frac{1}{H_m-1} \sum_{j=0}^\infty \frac{r_j}{(j+1)(j+2)},
  \end{equation}
  with the sequence \( (r_n) \) defined by \( r_j := 0 \) for \( 0
  \leq j \leq m-2 \) and
  \begin{equation}
    \label{eq:6.1.5.4}
    r_n := \frac{1}{\binom{n}{m-1}} \sum_{\mathbf{j}} [t_n + \mu_{j_1}
    + \cdots + \mu_{j_m} - \mu_n]^2, \qquad n \geq m-1.
  \end{equation}
  Moreover,
  \begin{equation*}
    \frac{X_n-\mu n}{\sqrt{n}} \stackrel{\mathcal{L}}{\to}
    \mathcal{N}(0,\sigma^2),
  \end{equation*}
  and there is convergence of moments of every order.
\end{theorem}
\begin{proof}
  One can check (for \emph{any} \( 2 \leq m < \infty \)) that the
  variance \( \sigma_n^2 \) vanishes for all \( n \geq m-1 \) if and
  only if the toll sequence is chosen as
  \begin{equation*}
    t_n = t \min\{m-1,n\}, \qquad n \geq 0
  \end{equation*}
  for some constant \( t \in \mathbb{R} \).  In that case without
  loss of generality \( t=1 \) and then \( X_n \equiv n \) is just the
  number of keys and we have exact (though degenerate) normality.  So
  we shall assume that \( \sigma^2 > 0 \).
  
  We use the method of moments together with the asymptotic transfer
  results of Section~\ref{sec:asympt-transf-theor}.

  Given  the toll sequence \( (t_n) \) defining the sequence \( (X_n)
  \) of random functionals of interests, the means \( (\mu_n) \)
  satisfy precisely the recurrence~\eqref{eq:2.7} with $(b_n)$
  replaced by $(t_n)$.
  Thus~\eqref{eq:6.1.5.2} simply repeats the asymptotic transfer
  result~\eqref{eq:3.2.4.5}.

  According to the law of total variance, the sequence \(
  (\sigma_n^2) \) also satisfies the recurrence~\eqref{eq:2.7}, but
  with \( (b_n) \) replaced by \( (r_n) \).  Furthermore, according to
  Lemma~\ref{lemma:6.1.5.6} to follow, the sequence \( (r_n) \)
  satisfies the condition~\eqref{eq:3.2.1.7}.  Then~\eqref{eq:6.1.5.3}
  is immediate from part~(a) of the ATT (Theorem~\ref{theorem:att}).

  Let \( \widetilde{X}_n := X_n - \mu(n+1) \) for \( n \geq 0 \).  We will
  complete the proof by showing by induction on \( k \) that
  \begin{equation}
    \label{eq:6.1.5.12}
    \tilde{\mu}_n(k) := \E{\widetilde{X}_n^k}, \quad k \geq 1 \quad
    \text{[with \( \tilde{\mu}_n(0) := 1 \)]}
  \end{equation}
  satisfies
  \begin{equation}
    \label{eq:6.1.5.13}
    \tilde{\mu}_n(2k) \sim \frac{(2k)!}{2^k k!} \sigma^{2k} n^k,
    \qquad\text{ and }\qquad \tilde{\mu}_n(2k-1) = o(n^{k-\tfrac12}),
    \quad k \geq 1.
  \end{equation}
  Observe that~\eqref{eq:6.1.5.2} and~\eqref{eq:6.1.5.3} imply
  that~\eqref{eq:6.1.5.13} holds for \( k=1 \).

  The key to the induction step for~\eqref{eq:6.1.5.13} is to apply
  the law of total expectation to~\eqref{eq:6.1.5.12}, by conditioning
  on the principal subtree sizes \( |L_1|,\ldots,|L_m| \).  Letting \(
  \sum_{\mathbf{k}} \) denote the sum over \( (m+1) \)-tuples \( (k_1,
  \ldots, k_{m+1} ) \) of nonnegative integers summing to \( k \), for \(
  n \geq m-1 \) we have in a manner identical to that for~\eqref{eq:2.3},
  \begin{equation}
    \label{eq:6.1.5.15}
    \tilde\mu_n(k) =  \frac{m}{\binom{n}{m-1}} \sum_{j=0}^{n-(m-1)}
    \binom{n-1-j}{m-2}   \tilde\mu_j(k) + r_n(k),
  \end{equation}
  where
  \begin{equation*}
    r_n(k) := \sumstark \binom{k}{k_1,\ldots,k_m,k_{m+1}}
    t_n^{k_{m+1}}  \frac{1}{\binom{n}{m-1}} \sum_{\mathbf{j}}
    \tilde\mu_{j_1}(k_1) \cdots \tilde\mu_{j_m}(k_m),
  \end{equation*}
  with \( \sum_{\mathbf{k}}^{*} \) denoting the same sum as
  \(\sum_{\mathbf{k}} \) with the additional restriction that \( k_i < k
  \) for \( i=1,\ldots,m \).  Observe that~\eqref{eq:6.1.5.15} is
  again of the same basic form~\eqref{eq:2.7}.  We will apply the ATT
  after evaluating \( r_n(k) \) asymptotically.

  We will treat the induction step in detail only for \(
  \tilde\mu_n(2k) \) at~\eqref{eq:6.1.5.13}, the case for \(
  \tilde\mu_n(2k-1) \) being similar and somewhat easier.  Let \(
  \sum_{\mathbf{k}}^{**} \) denote the sum over \(m\)-tuples \(
  (k_1,\ldots,k_m) \) of nonnegative integers, each strictly less than
  \( k\), summing to \(k\) (i.e., the same sum as \(
  \sum_{\mathbf{k}}^{*} \) with the additional restriction that \(
  k_{m+1}=0 \)).  For \( k \geq 2 \), we clearly have, by induction,
  \begin{align*}
    r_n(2k) &= o(n^k) + \sumstarsk
    \binom{2k}{2k_1,\ldots,{2k_m}}  \frac{1}{\binom{n}{m-1}}
    \sum_{\mathbf{j}} \tilde\mu_{j_1}(2k_1) \cdots
    \tilde\mu_{j_m}(2k_m) \\
    &= o(n^k) + \sumstarsk
    \binom{2k}{2k_1,\ldots,{2k_m}}  \frac{1}{\binom{n}{m-1}}
    \sum_{\mathbf{j}} \frac{(2k)!}{2^{k_1}k_1!} \sigma^{2k_1}
    j_1^{k_1} \cdots \frac{(2k_m)!}{2^{k_m}k_m!} \sigma^{2k_m}
    j_m^{k_m}\\
    &= o(n^k) + \frac{(2k)!}{2^k k!} \sigma^{2k} n^k
    \sumstarsk \binom{k}{k_1,\ldots, k_m} 
    \frac{1}{\binom{n}{m-1}} \sum_{\mathbf{j}}
    \left(\frac{j_1}{n}\right)^{k_1} \cdots
    \left(
      \frac{j_m}{n}
    \right)^{k_m}.
  \end{align*}
  But
  \begin{align*}
    &\frac{1}{\binom{n}{m-1}} \sum_{\mathbf{j}}
    \left(
      \frac{j_1}{n}
    \right)^{k_1} \cdots
    \left(
      \frac{j_m}{n}
    \right)^{k_m}\\
    &\to (m-1)! \int x_1^{k_1} \cdots x_{m-1}^{k_m-1}
    (1-x_1-\cdots-x_{m-1})^{k_m} \, dx_1\cdots dx_{m-1}\\
    &= (m-1)! \frac{\Gamma(k_1+1)\cdots\Gamma(k_m+1)}{\Gamma(k+m)} =
    \frac{1}{\binom{k}{k_1,\ldots,k_m}\binom{k+m-1}{m-1}},
  \end{align*}
  where the above integral is over \( (x_1,\ldots,x_{m-1})  \in
  [0,1]^{m-1} \) with sum not exceeding unity.  Since the number of
  terms in \( \sum_{\mathbf{k}}^{**} \) is \( \binom{k+m-1}{m-1} - m\), we
  therefore have
  \begin{equation*}
    r_n(2k) = \frac{(2k)!}{2^k k!} \sigma^{2k} n^k
    \frac{\binom{k+m-1}{m-1}-m}{\binom{k+m-1}{m-1}} + o(n^k)
    = \frac{(2k)!}{2^k k!} \sigma^{2k} n^k
    \left[
      1 - \frac{m! \Gamma(k+1)}{\Gamma(k+m)}
    \right]  + o(n^k), \qquad k \geq 2.
  \end{equation*}
  Similarly,
  \begin{equation*}
    r_n(2k-1) = o(n^{k-\tfrac12}), \qquad k \geq 2.
  \end{equation*}
  Now part~(c) of the ATT (Theorem~\ref{theorem:att})
  implies~\eqref{eq:6.1.5.13}.
\end{proof}

The following lemma lies at the heart of the proof
of Theorem~\ref{theorem:6.1.5.1}.
\begin{lemma}
  \label{lemma:6.1.5.6}
  In the context of Theorem~\ref{theorem:6.1.5.1}, the sequence \(
  (r_n) \) defined at~\eqref{eq:6.1.5.4} satisfies the
  conditions~\eqref{eq:3.2.1.7}.
\end{lemma}
\begin{proof}
  Clearly
  \begin{equation}
    \label{eq:6.1.5.16}
    r_n = \frac{1}{\binom{n}{m-1}} \sum_{\mathbf{j}} [t_n + \tilde\mu_{j_1}
    + \cdots + \tilde\mu_{j_m} - \tilde{\mu}_n]^2, \qquad n \geq m-1,
  \end{equation}
  with
  \begin{equation}
    \label{eq:6.1.5.17}
    \tilde\mu_n := \mu_n - \mu(n+1),
  \end{equation}
  which is \( o(\sqrt{n}) \) by~\eqref{eq:6.1.5.2}.  Recall the
  inequality
  \begin{equation}
    \label{eq:6.2.5.18}
    \left[
      \sum_{i=1}^k \xi_i
    \right]^2 \leq k \sum_{i=1}^k \xi^2
  \end{equation}
  for real numbers \( \xi_1,\ldots,\xi_k \).  Applying this
  to~\eqref{eq:6.1.5.16},
  \begin{equation}
    \label{eq:6.1.5.19}
    \frac{r_n}{m+2} \leq t_n^2 + \tilde\mu_n^2 +
    \frac{m}{\binom{n}{m-1}} \sum_{j=0}^{n-(m-1)} \binom{n-1-j}{m-2}
    \tilde\mu_j^2, 
  \end{equation}
  from which it is clear that \( r_n = o(n) \).

  To establish the summability of \( r_n/n^2 \) we need only establish
  that of \( \tilde\mu_n^2/n^2 \).  Indeed we may then
  use~\eqref{eq:6.1.5.19} again, together with the fact
  that~\eqref{eq:6.1.5.1}(b) implies
  \begin{equation}
    \label{eq:6.1.5.5}
    \sum^\infty \frac{t_n^2}{n^2} < \infty
  \end{equation}
  and the estimate
  \begin{multline*}
    \sum_{n=m-1}^\infty n^{-2} \frac{m}{\binom{n}{m-1}}
    \sum_{j=0}^{n-(m-1)} \binom{n-1-j}{m-2} \tilde\mu_j^2 =
    m(m-1)\sum_{j=0}^\infty \tilde\mu_j^2 \sum_{n=j+m-1}^\infty
    \frac{\falling{(n-1-j)}{m-2}}{n^2\falling{n}{m-1}} \\
    \leq m(m-1) \sum_{j=0}^\infty \tilde\mu_j^2 \sum_{n=j+m-1}^\infty
    n^{-3}
    = O\left( \sum \frac{\tilde\mu_j^2}{j^2} \right)< \infty.
  \end{multline*}

  To establish the summability of \( \tilde\mu_n^2/n^2 \), we recall
  from~\eqref{eq:3.2.4.9} and~\eqref{eq:3.2.4.2} that
  \begin{equation}
    \label{eq:6.1.5.20}
    \tilde\mu_n = O( n^{\alpha-1} ) + t_n - \frac{1}{H_m-1}(n+1)
    \sum_{k=n}^\infty \frac{\hat{t}_k}{(k+1)(k+2)} + \sum_{j=2}^{m-1}
    O
    \left(
      n^{\alpha_j-1} \sum_{k=0}^{n-1} \frac{|\hat{t}_k|}{(k+1)^{\alpha_j}}
    \right),
  \end{equation}
  writing \( \alpha_j := \Re(\lambda_j) \) (with \( \alpha = \alpha_2
  < 3/2 \), since \( m \leq 26 \)).  Using~\eqref{eq:6.2.5.18}, we
  need only establish the summability of \( n^{-2} \) times the square
  of each of the four terms on the right in~\eqref{eq:6.1.5.20}.  The
  first of these verifications is trivial, and the second was carried
  out at~\eqref{eq:6.1.5.5}.  For the third we apply the
  Cauchy--Schwarz inequality to give
  \begin{align*}
    \left[
      \sum_{k=n}^\infty \frac{\hat{t}_k}{(k+1)(k+2)}
    \right]^2 &\leq \frac{1}{n}
    \left[
      \sum_{k=n}^\infty \sqrt{n} \frac{\sqrt{k}}{(k+1)(k+2)}
      \frac{|t_k|}{\sqrt{k}}
    \right]^2\\
    &= O
    \left(
      \frac{1}{n} \sum_{k=n}^\infty \sqrt{n} k^{-3/2}
      \frac{t_k^2}{k}
    \right)
    = O
    \left(
      n^{-1/2} \sum_{k=n}^\infty k^{-5/2} t_k^2
    \right),
  \end{align*}
  whence
  \begin{multline*}
    \sum_n
    \left[
      \sum_{k=n}^\infty \frac{\hat{t}_k}{(k+1)(k+2)}
    \right]^2 = O
    \left(
      \sum_n n^{-1/2} \sum_{k=n}^\infty k^{-5/2} t_k^2
    \right)\\
    = O
    \left(
      \sum_k k^{-5/2} t_k^2 k^{1/2}
    \right)
    = O
    \left(
      \sum_k \frac{t_k^2}{k^2}
    \right) < \infty
  \end{multline*}
  by~\eqref{eq:6.1.5.5} again.

  We pause  to note that when \( m=2 \) the proof is finished
  here, and that up to now we have used only~\eqref{eq:6.1.5.5}, not
  the stronger assumption~\eqref{eq:6.1.5.1}(b).

  For our fourth and final verification, it suffices [again by
  invoking~\eqref{eq:6.2.5.18}] to establish the summability of
  \begin{equation}
    \label{eq:6.1.5.21}
    n^{2\rho-4}
    \left[
      \sum_{k=1}^{n-1} \frac{|t_k|}{k^\rho}
    \right]^2
  \end{equation}
  for any real \( \rho < 3/2 \).  To do this we break the sum into \(
  \sum_{k < n^\delta} \) and \( \sum_{n^\delta \leq k < n} \) and once
  again invoke~\eqref{eq:6.2.5.18}.  In the range
  \(\sum_{k<n^\delta}\) we simply use \(t_k = O(\sqrt{k})\) and note that
  \begin{equation*}
    n^{2\rho-4}
    \left[
      \sum_{k < n^\delta} O
      \left(
        k^{\tfrac12-\rho}
      \right)
    \right]^2 = O
    \left(
      n^{2\rho-4} (n^\delta)^{3-2\rho}
    \right) = O(n^\tau)
  \end{equation*}
  with \( \tau < -1 \).  In the range \( \sum_{n^\delta \leq k < n} \)
  we use Cauchy--Schwarz again:
  \begin{align*}
    n^{2\rho-4}
    \left[
      \sum_{n^\delta \leq k < n} \frac{|t_k|}{k^\rho}
    \right]^2 &= n^{2\rho-4} n^{3-2\rho}
    \left[
      \sum_{n^\delta \leq k < n}
      \frac{k^{\tfrac12-\rho}}{n^{\tfrac32-\rho}} \frac{|t_k|}{k^{1/2}}
    \right]^2\\
    &= O
    \left(
      n^{-1} \sum_{n^\delta \leq k < n}
      \frac{k^{\tfrac12-\rho}}{n^{\tfrac32-\rho}} \frac{t_k^2}{k}
    \right) = O
    \left(
      n^{-1} \max_{n^\delta \leq k < n} \frac{t_k^2}{k}
    \right),
  \end{align*}
  which is summable by assumption~\eqref{eq:6.1.5.1}(b).
\end{proof}

\begin{remark}
  \begin{enumerate}[(a)]
  \item We have already noted in the proof above that
    condition~\eqref{eq:6.1.5.1}(b) trivially implies
    \begin{equation}
      \label{eq:6.1.5.5a}
      \sum^\infty \frac{t_n^2}{n^2} < \infty,
    \end{equation}
    which in turn implies that~\eqref{eq:3.2.1.7} holds with absolute
    convergence; indeed, since the nonnegative numbers \(
    [(n+1)(n+2)]^{-1} \), \( n \geq 0 \), sum to unity, we have
    \begin{equation}
      \label{eq:6.1.5.6}
      \left[
        \sum_{n=0}^\infty \frac{|t_n|}{(n+1)(n+2)}
      \right]^2 \leq \sum_{n=0}^\infty \frac{t_n^2}{(n+1)(n+2)} < \infty.
    \end{equation}
  \item
    If
    \begin{equation}
      \label{eq:6.1.5.7}
      |t_n| = O(\tilde{t}_n) \quad\text{ with }\quad 0 \leq
      \frac{\tilde{t}_n}{\sqrt{n}} \downarrow \quad\text{ and }\quad
      \sum^\infty \frac{\tilde{t}_n^2}{n^2} < \infty,
    \end{equation}
    then we claim that~\eqref{eq:6.1.5.1} holds, and then as a corollary
    \begin{equation*}
      \frac{t_n}{\sqrt{n}} \downarrow 0 \quad\text{ and }\quad \sum^\infty
      \frac{t_n^2}{n^2} < \infty
    \end{equation*}
    implies~\eqref{eq:6.1.5.1}.  To see this claim, first observe that the
    condition~\eqref{eq:6.1.5.7} implies~\eqref{eq:6.1.5.1}(a); moreover,
    we observe that the series (say, over \(2 \leq n < \infty\))
    in~\eqref{eq:6.1.5.1}(b) is bounded by a constant times
    \begin{equation*}
      \sum_{n=2}^\infty n^{-1} \max_{n^\delta \leq k \leq n}
      \frac{\tilde{t}_k^2}{k} = \sum_{n=2}^\infty n^{-1}
      \frac{\tilde{t}_{\lceil n^\delta \rceil}^2}{\lceil n^\delta \rceil} =
      \sum_{k=2}^\infty \frac{\tilde{t}_k^2}{k} \sum_{(k-1)^{1/\delta} < n
        \leq k^{1/\delta}} n^{-1} = O
      \left(
        \sum_{k=2}^\infty \frac{\tilde{t}_k^2}{k^2}
      \right) < \infty.
    \end{equation*}
  \item
    We have also observed that, when \( m=2 \), the proof of
    Theorem~\ref{theorem:6.1.5.1} requires only~\eqref{eq:6.1.5.1}(a)
    and \eqref{eq:6.1.5.5a}.  In that case we obtain a strengthening of
    ``Case S1'' of Theorem~2 in~\cite{01808498} (for deterministic
    toll sequences); they required \(
    t_n=O(\sqrt{n}/(\ln{n})^{\tfrac12+\epsilon}) \) for some \( \epsilon
    > 0 \).
  \end{enumerate}
\end{remark}

We can also obtain asymptotic normality in the ``borderline small'' case.
The definition of slowly varying is available as
Definition~\ref{definition:slowly-varying}.
\begin{theorem}
  \label{theorem:clt5.4}
  If \(2 \leq m \leq 26\) and the real toll sequence \((t_n)\)
  satisfies
  \begin{equation}
    \label{eq:6.1.5.8}
    t_n \sim \sqrt{n} L(n)
  \end{equation}
  with \(L\) slowly varying, then the mean \(\mu_n\) of the
  corresponding additive functional \( X_n \) on \(m\)-ary search
  trees with the random permutation model satisfies
  \begin{equation}
    \label{eq:61.5.9}
    \mu_n = \frac{K_1}{H_m-1}n - \frac{\rising{(3/2)}{m-1}}{m! -
      \rising{(3/2)}{m-1}} \sqrt{n}L(n) + o(\sqrt{n}L(n)),
  \end{equation}
  with \(K_1\) defined at~\eqref{eq:6K1}.
  \begin{enumerate}
  \item   If \( \sum^\infty
    \frac{L^2(k)}{k} < \infty\), then the variance \( \sigma_n^2 \)
    satisfies~\eqref{eq:6.1.5.3}--\eqref{eq:6.1.5.4} and we define
    \begin{equation*}
      s^2(n) := \sigma^2n.
    \end{equation*}
  \item If \( \sum^\infty \frac{L^2(k)}{k} = \infty \), then
    \begin{equation}
      \label{eq:6.1.5.10}
      \sigma_n^2 \sim s^2(n) := \sigma^2 n \sum_{k \leq n} \frac{L^2(k)}{k},
    \end{equation}
    where in this case we define
    \begin{equation}
      \label{eq:6.1.5.11}
      \sigma^2 :=
      \frac{(\rising{(3/2)}{m-1})^2[\tfrac{\pi}4(m-1)+1]-(m!)^2}{(H_m-1)[m!-
        \rising{(3/2)}{m-1}]^2}.
    \end{equation}
  \end{enumerate}
  Moreover, in either case
  \begin{equation*}
    \frac{X_n - \mu n}{s(n)} \stackrel{\mathcal{L}}{\rightarrow}
    \mathcal{N}(0,1)
  \end{equation*}
  and there is convergence of moments of every order.
\end{theorem}
\begin{proof}[Proof sketch]
  Given the similarity to the proof of Theorem~\ref{theorem:6.1.5.1},
  we will be brief here.  Again we use the method of moments together
  with asymptotic transfer results.

  Equation~\eqref{eq:61.5.9} simply repeats the transfer
  result~\eqref{eq:3.2.4.6}.  As before, \((\sigma_n^2)\) satisfies
  the recurrence~\eqref{eq:2.7} with \((b_n)\) replaced by \((r_n)\)
  of~\eqref{eq:6.1.5.16}--\eqref{eq:6.1.5.17}, where again \(\mu
  := K_1/(H_m-1)\) and \(r_j := 0\) for \(0 \leq j \leq m-2\).  Here
  the proofs diverge somewhat.  The analogue of
  Lemma~\ref{lemma:6.1.5.6} is Lemma~\ref{lemma:6.1.5.7} below.  Then
  the asymptotic variance assertions of the theorem follow immediately
  from Theorem~\ref{theorem:3.2.4.5}(c).

  If \(\sum^\infty \frac{L^2(k)}{k} < \infty\), then
  from~\eqref{eq:3.2.4.10} applied to \(L^2\) it follows
  that~\eqref{eq:6.1.5.12} satisfies~\eqref{eq:6.1.5.13} for \(k=1\).
  Then higher moments are treated exactly as in the proof of
  Theorem~\ref{theorem:6.1.5.1} to complete the proof.

  If \(\sum^\infty \frac{L^2(k)}{k} = \infty\), then one
  uses~\eqref{eq:3.2.4.10}, Theorem~\ref{theorem:3.2.4.5}(d), and
  induction to show that the moments~\eqref{eq:6.1.5.12} satisfy
  \begin{equation}
    \label{eq:6.1.5.22}
    \tilde\mu_n(2k) \sim \frac{(2k)!}{2^k k!}s^{2k}(n) \quad\text{ and
    }\quad \tilde\mu_n(2k-1) = o(s^{2k-1}(n)), \qquad k \geq 1
  \end{equation}
  and thereby complete the proof.  We omit the details.
\end{proof}

The following cousin to Lemma~\ref{lemma:6.1.5.6} was used in the
proof of Theorem~\ref{theorem:clt5.4} above.
\begin{lemma}
  \label{lemma:6.1.5.7}
  In the context of Theorem~\ref{theorem:clt5.4}, the sequence
  \((r_n)\) defined for \(n \geq m-1\) by
  \begin{equation*}
    r_n := \frac{1}{\binom{n}{m-1}} \sum_{\mathbf{j}} [ t_n +
    \tilde\mu_{j_1} + \cdots + \tilde\mu_{j_m} - \tilde\mu_n ]^2
  \end{equation*}
  satisfies
  \begin{equation*}
    r_n \sim (H_m-1) \sigma^2 n L^2(n).
  \end{equation*}
\end{lemma}
\begin{proof}
  By~\eqref{eq:61.5.9}, with \(\theta := \rising{(3/2)}{m-1}/(m! -
  \rising{(3/2)}{m-1})\), we have
  \begin{align*}
    r_n &\sim \frac{1}{\binom{n}{m-1}} \sum_{\mathbf{j}} [ n^{1/2}
    L(n) - \theta_{j_1}^{1/2} L(j_1) - \cdots - \theta j_m^{1/2}
    L(j_m) + \theta n^{1/2} L(n) ]^2\\
    &\sim n L^2(n) (m-1)! \int
    \left[
      (1+\theta) - \theta \sum_{i=1}^{m-1} x_i^{1/2} - \theta
      \left(
        1- \sum_{i=1}^{m-1} x_i
      \right)^{1/2}
    \right]^2 \, dx_1 \cdots dx_{m-1},
  \end{align*}
  where the integral, call it \(J\), is over \((x_1,\ldots,x_{m-1})
  \in [0,1]^{m-1}\) with the sum not exceeding unity.  To complete the
  proof we need only show \( J = (H_m-1)\sigma^2/(m-1)! \), with
  \(\sigma^2\) defined at~\eqref{eq:6.1.5.11}.

  Indeed, \(J\) equals
  \begin{align*}
    & \int
    \left[
      (1+\theta)^2 + \theta^2 \sum_{i=1}^{m-1} x_i + \theta^2
      \left(
        1 - \sum_{i=1}^{m-1} x_i
      \right)
      -2\theta(1+\theta) \sum_{i=1}^{m-1} x_i^{1/2}
      -
      2\theta(1+\theta)
      \left( 1 - \sum_{i=1}^{m-1} x_i
      \right)^{1/2}\right.\\
    &\qquad \left. {}+ \theta^2 \sum_{i,j:i \ne j} x_i^{1/2} x_j^{1/2} +
      2\theta^2
      \left(
        1 - \sum_{i=1}^{m-1} x_i
      \right)^{1/2} \sum_{i=1}^{m-1} x_i^{1/2}
    \right] dx_1 \cdots dx_{m-1}\\
    &\quad = [(1+\theta)^2 + \theta^2] \frac{1}{(m-1)!} - 2\theta(1+\theta)m
    \frac{\Gamma(\tfrac32)}{\Gamma(m+\tfrac12)} + \theta^2 m(m-1)
    \frac{[\Gamma(3/2)]^2}{\Gamma(m+1)}\\
    &\quad = [(1+\theta)^2 + \theta^2] \frac{1}{(m-1)!} - 2\theta(1+\theta)m
    \frac{1}{\rising{(3/2)}{m-1}} + \theta^2 \frac{\pi/4}{(m-2)!}.
  \end{align*}
  Plugging in the value of \(\theta\) and simplifying, we obtain \(
  J=(H_m-1)\sigma^2/(m-1)! \), as desired.
\end{proof}

\begin{remark}
  When $m=2$ the constant $\sigma^2$ in~\eqref{eq:6.1.5.11} equals
  $\tfrac92 \pi - 14$, and Theorem~\ref{theorem:clt5.4} reduces to
  ``Case~S2'' of Theorem~2 in~\cite{01808498} (for deterministic toll
  sequences):\ see especially their displays~(15) and~(17), with
  $\tau_2 = 1$.
\end{remark}

\subsection*{\texorpdfstring{Periodicity for $m \geq 27$}{Periodicity
    for m bigger than 27}}
\label{sec:periodicity-m-geq}

If \( t_n = o(\sqrt{n}) \) as in Theorem~\ref{theorem:6.1.5.1} but \(
m \geq 27 \), then the remainder term \( \tilde\mu_n := \mu_n -
\mu(n+1) \) for the mean---which by \eqref{eq:6.1.5.2} was \(
o(\sqrt{n}) \) when \( m \leq 26 \)---now satisfies, by the ETT
(Theorem~\ref{theorem:ETT}) and~\eqref{eq:3.2.4.3}
[compare~\eqref{eq:3.2.4.9}]
\begin{multline}
  \label{eq:6.1.x}
  \tilde\mu_n = c_2 \frac{n^{\lambda_2-1}}{\Gamma(\lambda_2)} + c_3
  \frac{n^{\lambda_3-1}}{\Gamma(\lambda_3)}
  + m! \sum_{j=2}^{m-1}
  \frac{1}{\psi'(\lambda_j)} [z^n]
  \left(
    (1-z)^{-\lambda_j} \int_{\zeta=0}^z \widehat{T}(\zeta)
    (1-\zeta)^{\lambda_j-1} \,d\zeta
  \right)\\
  + o(\sqrt{n}) + O(n^{\Re(\lambda_4)  -1}).
\end{multline}
\emph{Typically} this will lead to the negative result that
\((\tilde\mu_n)\) [and hence also \((r_n)\) and \((\sigma_n^2)\)]
suffers from periodicity and that there is no natural distributional
limit for normalized \(X_n\).  See, e.g., Corollary~1
of~\cite{MR1871558}.

But it is perhaps difficult
to establish a \emph{general} negative result, due to
cancellations.  For example, suppose $T(z)$ equals $(1 - z)^{-1}$, so that
$t_n \equiv 1$, as studied by Chern and Hwang~\cite{MR1871558}, except
perhaps that the initial
values $t_0, \ldots, t_{m - 2}$ are changed.  Then
\[
T^{(m - 1)}(z) \equiv (m - 1)! (1 - z)^{- m},
\]
whence
\begin{align*}
  A(z)
  &:= \sum_{n=0}^\infty \mu_n z^n = \sum_{j = 1}^{m - 1} c_j (1 - z)^{-
    \lambda_j} + (m -  1)! \sum_{j = 1}^{m - 1} \frac{(1 - z)^{-
      \lambda_j}}{\psi'(\lambda_j)} \int_{\zeta = 0}^z\!(1 -
  \zeta)^{\lambda_j -     2}\,d\zeta \\ 
  &= \sum_{j = 1}^{m - 1} c_j (1 - z)^{- \lambda_j} - \frac{1}{m - 1} (1 -
  z)^{-1} - (m - 1)! \sum_{j = 1}^{m - 1} \frac{(1 - z)^{- \lambda_j}}{(1
    - \lambda_j)\psi'(\lambda_j)}.
\end{align*}
Now it is possible to choose $t_0, \ldots, t_{m - 2}$ so that
\begin{equation}
  \label{ctrex}
  c_j = \frac{(m - 1)!}{(1 - \lambda_j) \psi'(\lambda_j)}, \qquad j = 1,
  \ldots, m - 1. 
\end{equation}
In that case $A(z) = - \frac{1}{m - 1} (1 - z)^{-c}$, whence $\tilde{\mu}_n \equiv
\mu_n \equiv - 1/ (m - 1)$  [and we see that the
chosen values of $t_0, \ldots, t_{m - 2}$ are all $ - 1 / (m - 1)$], and we get
linear variance and asymptotic normality, just as in
Theorem~\ref{theorem:6.1.5.1}, for \emph{every} $2 \leq m < \infty$.

One might object that the above example is artificial, in that the toll
sequence changes sign.  But the same calculation show that if the toll sequence
is chosen as above ($t_n \equiv 1$) but with inital values
\[
t_j := K (j + 1) - \frac{1}{m - 1}, \qquad 0 \leq j \leq m - 2,
\]
then still, for every $m \geq 2$, the sequence $(\tilde{\mu}_n)$ is
constant, the 
variance is linear, and we have asymptotic normality.  Further, $(t_n)$ is
nonnegative provided $K \geq 1 / (m - 1)$.  [We remark in passing that the
choice $K = 1 / (m - 1)$ leads to the degenerate case mentioned at the
beginning of the proof of Theorem~\ref{theorem:6.1.5.1}.]
The sequence $(t_n)$ is also nondecreasing (as in most real examples) provided
$K \leq m / (m - 1)^2$.

\section{Moderate and large toll functions}
\label{sec:moderate-large-toll}

In order to describe the limiting distribution of $X_n$ for moderate
and large tolls, we will introduce
a family of random variables $Y \equiv Y_\beta$ defined for
$\beta > 1/2$, $\beta \ne 1$.  Anticipating Lemma~\ref{lemma:6.1.1},
we need to consider the distributional equation
\begin{equation}
  \label{eq:6.2.1}
  Y \stackrel{\mathcal{L}}{=} \sum_{j=1}^{m} S_j^\beta
  Y_{j} + 1
\end{equation}
Here $(Y_j)_{j=1}^m$ are independent copies of $Y$ and $(S_1,\ldots,S_m)$
is uniformly distributed on the $(m-1)$-simplex,
independent of
$(Y_j)_{j=1}^m$.  Recall that the $(m-1)$-simplex is the set
\[
\{(s_1,\ldots,s_{m})\!: s_j \geq 0 \text{ for } 1 \leq j \leq m
\text{ and } \mathbf{s_+} = 1 \}, \]
where $\mathbf{s_+}$ denotes $\sum_{j=1}^{m} s_j$.

Let
$U_{(1)},\ldots,U_{(m-1)}$ be the order statistics of a sample of size $m-1$
from $\Unif(0,1)$ with joint density
\begin{equation*}
  f_{U_{(1)},\ldots,U_{(m-1)}}(x_1,\ldots,x_{m-1}) \equiv (m-1)!
  \, \indicator{0 < x_1 <  \cdots < x_{m-1} < 1}
\end{equation*}
with respect to Lebesgue measure on $\mathbb{R}^{m-1}$, 
where $\indicator{A}$ is the indicator of $A$.

By a change of variables,  we find that the
joint distribution of the spacings $S_1,\ldots,S_{m}$, defined, with
$U_{(0)} := 0$ and $U_{(m)} := 1,$  by
\begin{equation*}
  S_{i} := U_{(i)} - U_{(i-1)}, \; i=1,\ldots,m,
\end{equation*}
is uniform over the $(m-1)$-simplex:
\begin{equation*}
  f_{S_1,\ldots,S_{m-1}}(s_1,\ldots,s_{m-1}) \equiv (m-1)!
  \,\indicator{s_j > 0,
    j=1,\ldots,m-1; \; \mathbf{s_{+}} < 1}.
\end{equation*}
When $r_j > -1$ for $1 \leq j \leq m$, observe that
\begin{align}
  \EE\,\left[\prod_{j=1}^m S_j^{r_j}\right] & = (m-1)!
  \int_{\substack{s_1,\ldots,s_{m-1} > 0\\ \mathbf{s_{+}}<1}}
  s_1^{r_1} \cdots s_{m-1}^{r_{m-1}}(1 - \mathbf{s_{+}})^{r_m} \,
  ds_{m-1}\cdots ds_1\notag\\
  \label{eq:6.2.2}
  & =: (m-1)! B(r_1+1,\ldots,r_m + 1)\\
  & = (m-1)! \frac{\prod_{j=1}^{m}
    \Gamma(r_j + 1)}{\Gamma(r_1+\cdots+r_m + m)}. \notag 
\end{align}
\begin{lemma}
  \label{lemma:6.1.1}
  Fix $\beta > 1/2$ with $\beta \ne 1$.  Then there
  exists a unique distribution $\mathcal{L}(Y) \equiv
  \mathcal{L}(Y_\beta)$ with finite second
  moment satisfying the distributional identity~\eqref{eq:6.2.1}.
\end{lemma}
\begin{proof}
  We first observe the that mean of any such distribution is
  determined by~\eqref{eq:6.2.1}.  Indeed, by taking expectations
  in~\eqref{eq:6.2.1} and using~\eqref{eq:6.2.2}, we get
  \begin{equation*}
    \mu := \EE\,Y = \left( 1 - \frac{m!
        \Gamma(\beta+1)}{\Gamma(\beta+m)} \right)^{-1}
  \end{equation*}
  since $\beta \ne 1$.  Thus we can equivalently
  consider the distributional identity
  \begin{equation*}
    W \stackrel{\mathcal{L}}{=} \sum_{j=1}^m S_j^\beta W_j + H,
  \end{equation*}
  where
  \begin{equation*}
    H := 1 - \mu + \mu \sum_{j=1}^m S_j^\beta.
  \end{equation*} Here $W$ is restricted to have mean~0 and finite
  second moment, $(W_j)_{j=1}^m$ are independent copies of $W$, and
  $(S_1,\ldots,S_m)$ is uniformly distributed on the $(m-1)$-simplex,
  independent of $(W_j)_{j=1}^m$.

  We now employ a standard contraction-method
  argument~\cite{MR2003b:68217,MR2003c:68096}.  Let $d_2$
  denote the metric on $\mathcal{M}_2(0)$, the space of probability
  distributions with mean~0 and finite variance, defined by
  \begin{equation*}
    d_2( G_1, G_2) := \min \lVert X_2 - X_1 \rVert_2,
  \end{equation*}
  taking the minimum over all pairs of random variables $X_1$ and
  $X_2$ defined on a common probability space with $\mathcal{L}(X_1)
  = G_1$ and $\mathcal{L}(X_2) = G_2$.  Here $\lVert \cdot \rVert_2$
  denotes  $L_2$-norm.

  Let $T$ be the map
  \begin{equation*}
    T: \mathcal{M}_2(0) \to \mathcal{M}_2(0), \quad G \mapsto
    \mathcal{L} \left( \sum_{j=1}^{m} S_j^\beta X_j + H \right),
  \end{equation*}
  where $(X_j)_{j=1}^m$ are independent with $\mathcal{L}(X_j) = G$,
  $j=1,\ldots,m$, and $(S_1,\ldots,S_m)$ is uniformly distributed on
  the $(m-1)$-simplex, independent of $(X_j)_{j=1}^m$.  We show that
  $T$ is a contraction on
  $\mathcal{M}_2(0)$; more precisely, that there exists a $\rho < 1$
  such that
  \begin{equation*}
    d_2(T(\mathcal{L}(A)),T(\mathcal{L}(B))) \leq \rho
    d_2(\mathcal{L}(A),\mathcal{L}(B))
  \end{equation*}
  for all pairs $\mathcal{L}(A)$ and $\mathcal{L}(B)$ in
  $\mathcal{M}_2(0)$.  To bound
  $d_2(T(\mathcal{L}(A)),T(\mathcal{L}(B)))$, we couple $T(\mathcal{L}(A))$ and
  $T(\mathcal{L}(B))$ by taking $m$ independent copies  $(A_j, B_j)$
  of the $d_2$-optimally coupled $(A,B)$,  an independent
  $(S_1,\ldots,S_m)$, and defining
  \begin{equation*}
    A' := \sum_{j=1}^m S_j^\beta A_j + H \sim T(\mathcal{L}(A)),\quad
    B' := \sum_{j=1}^m S_j^\beta B_j + H \sim T(\mathcal{L}(B)).
  \end{equation*}
  Now, defining $\mathbf{S} := (S_1,\ldots,S_m)$ and  using the law of
  total variance, 
  \begin{align*}
    &d_2( \mathcal{L}(T(A)), \mathcal{L}(T(B)) )^2\\
    \quad &\leq \lVert B' - A'
    \rVert_2^2
    = \left \lVert \sum_{j=1}^m S_j^\beta (B_j - A_j)
    \right\rVert_2^2
    = \Var\, \left[ \sum_{j=1}^m S_j^\beta (B_j - A_j) \right]\\
    &= \EE\, \Var\, \left[ \left.\sum_{j=1}^m S_j^\beta (B_j - A_j) \right|
      \mathbf{S} \right]
    + \Var\,
    \EE\, \left[ \left. \sum_{j=1}^m S_j^\beta (B_j - A_j) \right|
      \mathbf{S} \right]\\
    &= \sum_{j=1}^m (\EE S_j^{2\beta}) \Var\, [ B_j - A_j ] =
    d_2(\mathcal{L}(A),\mathcal{L}(B))^2  \sum_{j=1}^m \EE S_j^{2\beta}
    = m! \frac{\Gamma(2\beta+1)}{\Gamma(2\beta+m)} d_2(A,B)^2.
  \end{align*}
  We need only verify that
  \begin{equation*}
    \rho^2 := m! \frac{\Gamma(2\beta+1)}{\Gamma(2\beta+m)} =
    \frac{m!}{(2\beta+m-1) \cdots (2\beta+1)} < 1,
  \end{equation*}
  which is true when $\beta > 1/2$.  The existence and uniqueness of
  $\mathcal{L}(Y)$ now follows from the Banach fixed point
  theorem~\cite[Theorem~2]{MR2003c:68096}.
\end{proof}

\subsection*{Moderate toll functions}
\label{sec:moder-toll-funct}
In the case of moderate toll functions, convergence in distribution
and convergence of all moments can be stated as
\begin{theorem}
  \label{thm:F2}
  If the real toll sequence~$(t_n)$ satisfies
  \[
  t_n \sim n^{\beta}, \qquad 1/2 < \beta < 1,
  \]
  and $\alpha < 1 + \beta$, then the mean of the corresponding linear
  functional $X_n$ on $m$-ary search trees with the random permutation
  model satisfies
  \begin{equation}
    \label{eq:F2.1}
    \mu_n = \mu n -
    \frac{\rising{(1+\beta)}{m-1}}{m!-\rising{(1+\beta)}{m-1}}
    n^{\beta} + o(n^{\beta}),  \qquad \mu := \frac{K_1}{H_m-1},
  \end{equation}
  with $K_1$ defined at~\eqref{eq:6K1}.  Moreover,
  \[
  \frac{X_n-\mu{}n}{n^{\beta}} \stackrel{\mathcal{L}}{\to}
  Y_{\beta},
  \]
  with convergence of all moments.
\end{theorem}
\begin{proof}
  We use the notation introduced in the proof of
  Theorem~\ref{theorem:6.1.5.1}.
  Equation~(\ref{eq:F2.1}) is simply a restatement of the asymptotic
  transfer result~\eqref{eq:3.2.4.6}.

  We show that the moments $\tilde\mu_n(k)$ satisfy
  \begin{equation}
    \label{eq:F2.5}
    \tilde\mu_n(k) = g_k n^{k\beta} + o(n^{k\beta}) \text{
      as } n   \to \infty.
  \end{equation}
  The claim holds for $k=1$ with
  \begin{equation}
    \label{eq:6.2.8}
    g_1 := -\frac{\rising{(1+\beta)}{m-1}}{m!-\rising{(1+\beta)}{m-1}} =
    \left( 1 - \frac{m!\Gamma(\beta+1)}{\Gamma(\beta+m)} \right)^{-1}.
  \end{equation}
  Using~\eqref{eq:6.1.5.16}
  and~\eqref{eq:6.1.5.15}, by induction we get, for $k \geq 2$,
  \begin{align*}
    r_n(k) &= o(n^{k\beta}) +  \sumstark
    \binom{k}{k_1,\ldots,k_m,k_{m+1}}
    (n^{\beta})^{k_{m+1}} \frac{1}{\binom{n}{m-1}}
    \sum_{\mathbf{j}} g_{k_1}j_1^{k_1\beta} \cdots
    g_{k_m}j_m^{k_m\beta} \\
    &= o(n^{k\beta})
    +  \sumstark
    \binom{k}{k_1,\ldots,k_m,k_{m+1}}
    n^{k\beta} g_{k_1} \cdots g_{k_m} \frac{1}{\binom{n}{m-1}}
    \sum_{\mathbf{j}} \left(\frac{j_1}{n}\right)^{k_1\beta}
    \cdots \left(\frac{j_m}{n}\right)^{k_m\beta}.
  \end{align*}
  But
  \[
  \frac{1}{\binom{n}{m-1}} \sum_{\mathbf{j}}
  \left(\frac{j_1}{n}\right)^{k_1\beta}  \cdots
  \left(\frac{j_m}{n}\right)^{k_m\beta} \to (m-1)!
  B(k_1\beta+1,\ldots k_m\beta+1)
  \]
  so that
  \[
  r_n(k) = o(n^{k\beta}) +  n^{k\beta} (m-1)! \sumstark
  \binom{k}{k_1,\ldots,k_m,k_{m+1}} g_{k_1}\cdots g_{k_m}
  B(k_1\beta+1,\ldots,k_m\beta+1).
  \]
  Using Theorem~\ref{theorem:att} (the ATT), with $v=k\beta > 1$, we get
  \begin{equation*}
    \tilde\mu_n(k) = o(n^{k\beta})
    +  n^{k\beta} \frac{(m-1)!}{1-\frac{m!
        \Gamma(k\beta+1)}{\Gamma(k\beta+m)}}
    \sumstark \binom{k}{k_1,\ldots,k_m,k_{m+1}} g_{k_1}
    \cdots g_{k_m} B(k_1\beta+1,\ldots,k_m\beta+1).
  \end{equation*}
  Thus, defining $g_k$ recursively as
  \begin{equation}
    \label{eq:6.2.9}
    g_k =
    \frac{(m-1)!}{1-\frac{m!\Gamma(k\beta+1)}{\Gamma(k\beta+m)}}
    \sumstark \binom{k}{k_1,\ldots,k_{m+1}} g_{k_1}\cdots
    g_{k_m} B(k_1\beta+1,\ldots,k_m\beta+1),
  \end{equation}
  we see that~\eqref{eq:F2.5} holds for all $k \geq 1$.

  By Lemma~\ref{lem:F7.1} (to follow) and the method of moments
  (cf., e.g.,~\cite[Sections~4.4 and~4.5]{MR49:11579}), the $g_k$'s are the
  moment of a uniquely determined distribution, say
  $\mathcal{L}(\widehat{Y})$, and
  \begin{equation*}
    \frac{X_n - \mu n}{n^\beta} \stackrel{\mathcal{L}}{\to} \widehat{Y}
  \end{equation*}
  with convergence of all moments.  It remains to show that $\widehat{Y}
  \stackrel{\mathcal{L}}{=} Y_\beta$.

  Define
  \begin{equation}
    \label{eq:6.2.7}
    \widetilde{Y} := \sum_{j=1}^m S_j^\beta \widehat{Y}_j + 1,
  \end{equation}
  where $(\widehat{Y}_j)_{j=1}^m$ are independent copies of $\widehat{Y}$ and
  $(S_1,\ldots,S_m)$ is uniformly distributed on the $(m-1)$-simplex,
  independent of $(\widehat{Y}_j)_{j=1}^m$.  We will show that $\widehat{Y}
  \stackrel{\mathcal{L}}{=} \widetilde{Y}$, and then, by~\eqref{eq:6.2.7},
  $\mathcal{L}(\widehat{Y})$ satisfies the distributional
  identity~\eqref{eq:6.2.1} and has finite second moment.  By
  Lemma~\ref{lemma:6.1.1}, $\widehat{Y} \stackrel{\mathcal{L}}{=}
  Y_\beta$, as desired.

  To show $\widehat{Y} \stackrel{\mathcal{L}}{=} \widetilde{Y}$, it
  suffices to show that
  $\widehat{Y}$ and $\widetilde{Y}$ have the same moments.  Letting
  $\sum_{\mathbf{k}}$ denote (as before) the sum over $(m+1)$-tuples
  $(k_1,\ldots,k_{m+1})$ of nonnegative integers summing to $k$, and
  using~\eqref{eq:6.2.7},
  \eqref{eq:6.2.2}, and~\eqref{eq:6.2.9},
  \begin{align*}
    \EE\,\widetilde{Y}^k &= \sum_{\mathbf{k}} \binom{k}{k_1,\ldots,k_{m+1}}
    \EE\,\left[\prod_{j=1}^m (S_j^{\beta} \widehat{Y}_j)^{k_j} \right] \\
    &= \sum_{\mathbf{k}} \binom{k}{k_1,\ldots,k_{m+1}} (m-1)! B(k_1\beta
    + 1, \ldots k_m\beta + 1) g_{k_1} \cdots g_{k_m}\\
    &= (m-1)! \sumstark
    \binom{k}{k_1,\ldots,k_{m+1}}  B(k_1\beta 
    + 1, \ldots k_m\beta + 1) g_{k_1} \cdots g_{k_m}\\
    & \quad + m!  B(k\beta+1,1,\ldots,1) g_k\\
    &= \left[ 1 - \frac{m! \Gamma(k\beta+1)}{\Gamma(k\beta+m)} \right]
    g_k + \frac{m! \Gamma(k\beta+1)}{\Gamma(k\beta+m)} g_k = g_k = \EE\,
    \widehat{Y}^k,
  \end{align*}
  where (as before) $\sum_{\mathbf{k}}^{\mathbf{*}}$ denotes the same
  sum as $\sum_{\mathbf{k}}$ with the additional restriction that $k_i <
  k$ for $i=1,\ldots,m$.
\end{proof}

\begin{lemma}
  \label{lem:F7.1}
  The moments $(g_k)$ uniquely determine the
  distribution $\mathcal{L}(Y)$.
\end{lemma}
\begin{proof}
  Define $\gamma_k := g_k/k!$.  It suffices to show that there
  exists an $M$ such that $\gamma_k \leq M^k$ for all $k$ sufficiently
  large.  We proceed by induction.  Indeed, by~\eqref{eq:6.2.9} we know
  \begin{align*}
    \gamma_k &= \frac{(m-1)!}{1 -
      \frac{m!\Gamma(k\beta+1)}{\Gamma(k\beta+m)}} \sumstark
    \frac{1}{k_{m+1}!} \left( \prod_{j=1}^m \gamma_{k_j} \right)
    B(k_1\beta+1, \ldots, k_m\beta+1)\\
    & \leq M^k \frac{(m-1)!}{1 -
      \frac{m!\Gamma(k\beta+1)}{\Gamma(k\beta+m)}} \sumstark
    \frac{M^{-k_{m+1}}}{k_{m+1}!} B(k_1\beta+1,\ldots,k_m\beta+1)
  \end{align*}
  by the induction hypothesis.  So it is certainly sufficient to show
  that
  \begin{align}
    \label{eq:6.2.3}
    &\sumstark
    \frac{M^{-k_{m+1}}}{k_{m+1}!} B(k_1\beta+1,\ldots,k_m\beta+1)
    \notag\\ 
    &= \sum_{k_{m+1}=0}^k \frac{M^{-k_{m+1}}}{k_{m+1}!} \Gamma(
    (k-k_{m+1})\beta + m)^{-1} \sum_{\substack{0 \leq k_1,\ldots,k_m <
        k\\ k_1 + \cdots + k_m = k - k_{m+1}}} \prod_{j=1}^m \Gamma(
    k_j\beta + 1) \notag\\
    &\to    0 \quad \text{ as $k \to \infty$}.
  \end{align}
  
  For this, fix a value of $k_{m+1} \in \{0,1,2,\ldots\}$, and
  consider the sum
  \begin{equation}
    \label{eq:6.2.5}
    \sum_{\substack{0 \leq
        k_1,\ldots,k_m < k\\ k_1+\cdots+k_m=k-k_{m+1}}} \prod_{j=1}^m
    \Gamma(k_j\beta + 1).
  \end{equation}
  By log-convexity of $\Gamma$~\cite[6.4.1]{MR29:4914},
  taking $I$ to be $(0,\infty)$ and $g$ to be $\ln \Gamma$ in
  Proposition~3.C.1 of~\cite{MR81b:00002}, the logarithm of the
  product in~\eqref{eq:6.2.5} is Schur-convex on $(0,\infty)^m$.  Thus applying
  Proposition~5.C.2 of~\cite{MR81b:00002} with $m=0$ there, the
  biggest terms in
  the sum correspond to $k_{j} = k-k_{m+1}$ for $j$ equal to some
  $j_0$ and $k_j=0$ otherwise; together, these $m$ terms contribute
  $m\Gamma((k-k_{m+1})\beta + 1)$ to the sum.  If $k > k_{m+1}$, there
  are other terms, the biggest of which corresponds to having one of
  the $k_j$'s be $k-k_{m+1}-1$, one be~1, and the rest be~0.  (This
  follows from Proposition~5.C.1 of~\cite{MR81b:00002} with $m=0$ and
  $M=k-k_{m+1}-1$ there.)  The  total number of terms in the
  sum~\eqref{eq:6.2.5} is at most   $\binom{m-1+(k-k_{m+1})}{m-1}$.  So
  the remaining  contribution  to~\eqref{eq:6.2.5} is at most
  \begin{equation*}
    {\binom{m-1+(k-k_{m+1})}{m-1}}
    \Gamma( (k-k_{m+1}-1)\beta + 1) \Gamma(\beta + 1).
  \end{equation*}
  We have found that the left side of~\eqref{eq:6.2.3} is bounded by
  \begin{equation*}
    \sum_{k_{m+1}=0}^{\infty} \frac{M^{-k_{m+1}}}{k_{m+1}!}
    \indicator{k_{m+1} \leq k} f(k-k_{m+1}),
  \end{equation*}
  where
  \begin{align}
    \label{eq:6.2.6}
    f(k) & \leq \frac{m}{(k\beta+1) \cdots (k\beta+m-1)} +
    \Gamma(\beta+1) \frac{\binom{m-1+k}{m-1}\Gamma((k-1)\beta +
      1)}{\Gamma(k\beta+m)} \\
    & \leq \frac{\binom{m-1+k}{m-1}\Gamma(k\beta+1)}{\Gamma(k\beta+m)}
    = \frac{1}{(m-1)!} \frac{(k+1) \cdots (k+(m-1))}{(k\beta+1) \cdots
      (k\beta + (m-1))},
  \end{align}
  which is a bounded function of $k$.  To apply the dominated
  convergence theorem, it suffices to show that the right side
  of~\eqref{eq:6.2.6} tends to 0 as $k \to \infty$, which follows from
  Stirling's approximation and the fact that $\beta > 0$.
\end{proof}

When $t_n$ satisfies the conditions in Theorem~\ref{thm:F2} but
$\alpha \geq 1+\beta$ then [compare~\eqref{eq:6.1.x}],
\begin{align*}
  \tilde\mu_n &= c_2\frac{n^{\lambda_2-1}}{\Gamma(\lambda_2)} + c_3
  \frac{n^{\lambda_3-1}}{\Gamma(\lambda_3)}
  + m!\sum_{j=1}^{m-1} \frac{1}{\psi'(\lambda_j)} [z^n]
  \left((1-z)^{-\lambda_j} \int_{\zeta=0}^z
    \widehat{T}(\zeta)(1-\zeta)^{\lambda_j-1}\,d\zeta\right) \\
  & \qquad + o(n^{\beta}) + O(n^{\Re(\lambda_4-1)}),
\end{align*}
and typically this leads to periodicity.

It is known that $\alpha < 3/2$ for $m \leq 26$.  Computations show
that $\alpha$ increases with $m$, at least for $m \leq
10000$~\cite{MR96j:68042}.  This suggests that for a fixed $\beta$,
the condition $\alpha < 1 + \beta$ (resp.\ $\alpha \geq 1 + \beta$) is
equivalent to $m \leq m_0$ (resp.\ $m > m_0$) for some $m_0 > 26$.

\subsection*{Large toll functions}
\label{sec:rpm-large-toll-functions}

If $t_n \sim n^\beta$, where $\beta > 1$, then we have convergence in
distribution for all values of $m$.  We state the result, omitting the
proof, as it is very similar to that of Theorem~\ref{thm:F2}.

\begin{theorem}
  \label{thm:F3}
  If the real toll sequence~$(t_n)$ satisfies
  \[
  t_n \sim n^\beta, \qquad \beta > 1,
  \]
  then
  \[
  \frac{X_n}{n^\beta} \stackrel{\mathcal{L}}{\to} Y_{\beta}
  \]
  with convergence of all moments, where $\mathcal{L}(Y_\beta)$ is the
  unique distribution satisfying~\eqref{eq:6.2.1}.
\end{theorem}

The borderline case $t_n \sim n$ [more specifically, $t_n \equiv n-(m-1)$]
corresponds to the well-studied 
total path length of a random $m$-ary search tree.  The corresponding
additive functional measures the 
number of key comparisons in $m$-ary
Quicksort.  As is well known~\cite{MR92f:68028}, the number
of key comparisons has mean $\Theta(n \ln n)$ and standard deviation
$\Theta(n)$.
See~\cite{hennequin91:_analy} for 
details in this case.


\chapter{Catalan model: the toll function
  \texorpdfstring{$n^\alpha$}{n to the alpha}}
\label{sec:toll-function-nalpha}

We will now use the extended singularity analysis toolkit from
Chapter~\ref{cha:sing-analys-hadam} to derive
limiting distributions for the additive functionals on binary search
trees under the uniform model for the class of toll functions
of the form \( n^\alpha \) where \( \alpha > 0 \).

\section{Preliminaries}
\label{sec:n_alpha_preliminaries}
In the sequel the notation [\(\cdots\)] is used both for Iverson's
convention~\cite[see~1.2.3(16)]{knuth97} and the coefficient of
certain terms in the succeeding expression.  The interpretation will
be clear from the context.  For example, \( [ \alpha > 0 ] \) has the
value 1 when \( \alpha > 0 \) and the value 0 otherwise.  In contrast, \( [z^n]
F(z) \) denotes the coefficient of \( z^n \) in the series expansion
of \( F(z) \).  For notational convenience we define $L(z) :=
\ln((1-z)^{-1})$.  Unless otherwise noted, henceforth all our singular
expansions are valid as \( z \to 1 \) in
the sector \( -\pi+\epsilon < \arg(1-z) < \pi-\epsilon \), for any
$\epsilon > 0$.

We will  make extensive use of the following consequences of
the singular expansion of the generalized polylogarithm.  Note that
neither this lemma nor the ones following make any claims about
uniformity in $\alpha$ or $r$.
\begin{lemma}
  \label{lem:litoomz}
  For real $\alpha < 1$ and $r$ a nonnegative integer,
  \[
  \Li_{\alpha,r}(z) = \sum_{k=0}^r \lambda_k^{(\alpha,r)}
  (1-z)^{\alpha-1} L^{r-k}(z) +
  O(|1-z|^{\alpha-\epsilon}) + (-1)^r \zeta^{(r)}(\alpha)[\alpha > 0],
  \]
  where $\lambda_k^{(\alpha,r)} \equiv \binom{r}{k}
  \Gamma^{(k)}(1-\alpha)$ and $\epsilon > 0$ is arbitrarily small.
\end{lemma}
\begin{proof}
  By Theorem~\ref{theorem:Li-sing},
  \begin{equation}\label{eq:37}
    \Li_{\alpha,0}(z) \sim \Gamma(1-\alpha) t^{\alpha-1} + \sum_{j
    \geq 0} \frac{(-1)^j}{j!} \zeta(\alpha-j) t^j, \qquad t = -\ln z
    = \sum_{l=1}^\infty \frac{(1-z)^l}{l},
  \end{equation}
  and for \( r \) a positive integer,
  \begin{equation*}
    \Li_{\alpha,r}(z) = (-1)^r \frac{\partial^r}{\partial\alpha^r}
    \Li_{\alpha,0}(z).
  \end{equation*}
  Moreover, the singular expansion for $\Li_{\alpha,r}$ is obtained by
  performing the indicated differentiation of~(\ref{eq:37}) term-by-term.
  To establish the claim we set \( f = \Gamma(1-\alpha) \) and \( g =
  t^{\alpha-1} \) in the general formula for the $r$th derivative of a
  product:
   \begin{equation*}
    (fg)^{(r)} = \sum_{k=0}^r \binom{r}{k} f^{(k)} g^{(r-k)}
   \end{equation*}
   to first obtain
   \begin{equation*}
     (-1)^r \frac{\partial^r}{\partial\alpha^r}
    [\Gamma(1-\alpha)t^{\alpha-1}] = (-1)^r \sum_{k=0}^r \binom{r}{k}
    (-1)^k\Gamma^{(k)}(1-\alpha) t^{\alpha-1} (\ln t)^{r-k}
   \end{equation*}
   The claim then follows easily.
\end{proof}
The following ``inverse'' of Lemma~\ref{lem:litoomz} is very useful
for computing with Hadamard products.
\begin{lemma}
  \label{lem:omztoli}
  For real $\alpha < 1$ and $r$ a nonnegative integer,
  \[
  (1-z)^{\alpha-1}L^r(z) = \sum_{k=0}^r
  \mu_k^{(\alpha,r)} \Li_{\alpha,r-k}(z) + O(|1-z|^{\alpha-\epsilon})
  + c_r(\alpha)[\alpha > 0],
  \]
  where $\mu_0^{(\alpha,r)}=1/\Gamma(1-\alpha)$, \( c_r(\alpha) \) is
  a constant, and
  $\epsilon > 0$ is  arbitrarily small.
\end{lemma}
\begin{proof}
  We use induction on~$r$. For $r=0$ we have
  \[
  \Li_{\alpha,0}(z) = \Gamma(1-\alpha)(1-z)^{\alpha-1} +
  O(|1-z|^{\alpha-\epsilon}) + \zeta(\alpha) [ \alpha > 0 ]
  \]
  and the claim is verified with
  \begin{equation*}
    \mu_0^{(\alpha,0)} = \frac{1}{\Gamma(1-\alpha)} \qquad  \text{ and
    } \qquad c_0(\alpha) = -\frac{\zeta(\alpha)}{\Gamma(1-\alpha)}.
  \end{equation*}
  Let $r \geq 1$.  Then using
  Lemma~\ref{lem:litoomz} and the induction hypothesis we get
  \begin{align*}
    \Li_{\alpha,r}(z) &=
    \Gamma(1-\alpha)(1-z)^{\alpha-1}
    L^r(z) \\
    & \quad + \sum_{k=1}^r \lambda_{k}^{(\alpha,r)} \left[
      \sum_{l=0}^{r-k}\mu_l^{(\alpha,r-k)} \Li_{\alpha,r-k-l}(z) +
    O(|1-z|^{\alpha-\epsilon}) +
      c_{r-k}(\alpha) [\alpha > 0 ] \right]\\
    & \quad + O(|1-z|^{\alpha-\epsilon})
    + (-1)^r\zeta^{(r)}(\alpha)[\alpha > 0]\\
    &=
    \Gamma(1-\alpha)(1-z)^{\alpha-1}L^r(z)
    + \sum_{k=1}^r \lambda_{k}^{(\alpha,r)} \sum_{s=0}^{r-k}
    \mu_{r-k-s}^{(\alpha,r-k)} \Li_{\alpha,s}(z) +
    O(|1-z|^{\alpha-\epsilon}) \\
    & \quad +
    \left(
      \sum_{k=1}^r \lambda_k^{(\alpha,r)} c_{r-k}(\alpha)
      + (-1)^r\zeta^{(r)}(\alpha)
    \right)[\alpha > 0]\\
    &=
    \Gamma(1-\alpha)(1-z)^{\alpha-1}
    L^r(z)
    + \sum_{s=0}^{r-1} \nu_s^{(\alpha,r)} \Li_{\alpha,s}(z) +
    O(|1-z|^{\alpha-\epsilon}) + \gamma_r(\alpha)[\alpha > 0],
  \end{align*}
  where, for $0 \leq s \leq r-1$,
  \[
  \nu_s^{(\alpha,r)} := \sum_{k=1}^{r-s} \lambda_k^{(\alpha,r)}
  \mu_{r-s-k}^{(\alpha,r-k)},
  \]
  and where
  \begin{equation*}
    \gamma_r(\alpha) := \sum_{k=1}^r
    \lambda_k^{(\alpha,r)}c_{r-k}(\alpha) + (-1)^r \zeta^{(r)}(\alpha).
  \end{equation*}
  Setting
  \begin{equation*}
    \mu_0^{(\alpha,r)} = \frac{1}{\Gamma(1-\alpha)},\qquad 
    \mu_k^{(\alpha,r)} =
    -\frac{\nu_{r-k}^{(\alpha,r)}}{\Gamma(1-\alpha)}, \quad 1 \leq k
    \leq r,
  \end{equation*}
  and
  \begin{equation*}
    c_r(\alpha) = - \frac{\gamma_r(\alpha)}{\Gamma(1-\alpha)},
  \end{equation*}
  the result follows.
\end{proof}
For the calculation of the mean, the following refinement of a special
case of Lemma~\ref{lem:litoomz} is required.  It is a simple
consequence of Theorem~\ref{theorem:Li-sing}.
\begin{lemma}
  \label{lem:li0omz}
  When \(\alpha < 0\),
  \[
  \Li_{\alpha,0}(z) = \Gamma(1-\alpha)(1-z)^{\alpha-1} -
  \Gamma(1-\alpha)\frac{1-\alpha}2 (1-z)^{\alpha} +
  O(|1-z|^{\alpha+1}) + \zeta(\alpha)[\alpha > -1] .
  \]
\end{lemma}

\section{Asymptotics of the mean}
\label{sec:cat_n_alpha_asymptotics-mean}

We start by reviewing some notation.
Let \(\beta_n\) denote the \(n\)th Catalan number
\[
\beta_n := \frac1{n+1}\binom{2n}{n},
\]
with generating function
\[
\CAT(z) := \sum_{n=0}^\infty \beta_n z^n = \frac1{2z}(1 - \sqrt{1-4z}).
\]
In our subsequent analysis we will make use of the identity
\begin{equation}
  \label{eq:7.50}
  z \CAT^2(z) = \CAT(z) - 1.
\end{equation}
The mean of the cost function \(a_n := \E{X_n} \) can be
obtained recursively by conditioning on the key stored at the root as
\begin{equation*}
  a_n = \sum_{j=1}^n \frac{\beta_{j-1}\beta_{n-j}}{\beta_n}(a_{j-1} +
  a_{n-j}) + t_n, \qquad n \geq 1.
\end{equation*}
This recurrence can be rewritten as
\begin{equation}
  \label{eq:7.35}
  (\beta_n a_n) = 2 \sum_{j=1}^n (\beta_{j-1} a_{j-1}) \beta_{n-j} +
  (\beta_n t_n), \qquad n \geq 1.
\end{equation}
Multiplying~(\ref{eq:7.35}) by $z^n/4^n$ and summing over $n \geq 1$ we get
\begin{equation}
  \label{eq:7.1}
A(z)\odot\CAT(z/4) = \frac{T(z)\odot\CAT(z/4)}{\sqrt{1-z}} +
\frac{a_0-t_0}{\sqrt{1-z}},
\end{equation}
where \(A(z)\) and \(T(z)\) are the ordinary generating functions of
\((a_n)\) and \((t_n)\) respectively.  In our analysis we will assume
$a_0=t_0$ so that the second term in~\eqref{eq:7.1} does not contribute.
(Our results can be easily adjusted when this is not the case.)  We
will use~(\ref{eq:7.1}) (with $a_0 = t_0$) as the starting point for our
analysis.

Since $t_n = n^\alpha$, by definition \(T=\Li_{-\alpha,0}\). Thus,
by Lemma~\ref{lem:li0omz},
\[
T(z) = \Gamma(1+\alpha)(1-z)^{-\alpha-1} -
\Gamma(1+\alpha)\frac{\alpha+1}2(1-z)^{-\alpha}  +
O(|1-z|^{-\alpha+1}) + \zeta(-\alpha)[\alpha < 1].
\]
We will now use~(\ref{eq:7.1}) to obtain the asymptotics of the mean.

First we treat the case
\(\alpha < 1/2\).
Since \(\CAT(z/4) = 2 + O(|1-z|^{1/2})\) as \(z \to 1\) in a suitable
indented crown, we have, by Theorem~\ref{thm:hadamard}(b),
\[
T(z)\odot\CAT(z/4)  = C_0 + O(|1-z|^{-\alpha+\tfrac12}),
\]
where
\begin{equation*}
  \label{eq:7.20}
  C_0 := T(z)\odot\CAT(z/4) \Bigr\rvert_{z=1} = \sum_{n=1}^\infty n^\alpha
  \frac{\beta_n}{4^n}.
\end{equation*}
We now already know the constant
term in the singular expansion of \(T(z)\odot\CAT(z/4)\) at \(z=1\)
and henceforth we need only compute lower-order terms.  The constant
\(\bar{c}\) is used in the sequel to denote an unspecified
(possibly~0) constant, possibly different at each appearance.

Let us write \(T(z) = L_1(z) + R_1(z)\), and \(\CAT(z/4) = L_2(z) +
R_2(z)\), where
\begin{align*}
  L_1(z) &:= \Gamma(1+\alpha)(1-z)^{-\alpha-1} -
  \Gamma(1+\alpha)\frac{\alpha+1}2(1-z)^{-\alpha} + \zeta(-\alpha),\\
  R_1(z) &:= T(z) - L_1(z) = O(|1-z|^{1-\alpha}),\\
  L_2(z) &:= 2(1 - (1-z)^{1/2}),\\
  R_2(z) &:= \CAT(z/4) - L_2(z) = O(|1-z|).
\end{align*}
We will analyze each of the four Hadamard products separately. First,
\begin{equation*}
  L_1(z)\odot{}L_2(z) = -
 2\Gamma(1+\alpha)\bigl[(1-z)^{-\alpha-1}\odot(1-z)^{1/2}\bigr]
  +  2\Gamma(1+\alpha)\frac{\alpha+1}2
 \bigl[(1-z)^{-\alpha}\odot(1-z)^{1/2}\bigr] + \bar{c}.
\end{equation*}
By Theorem~\ref{thm:hadamard},
\[
(1-z)^{-\alpha-1}\odot(1-z)^{1/2} = \bar{c} +
\frac{\Gamma(\alpha-\tfrac12)}{\Gamma(\alpha+1)\Gamma(-1/2)}
(1-z)^{-\alpha+\tfrac12} + O(|1-z|),
\]
and
\[
(1-z)^{-\alpha}\odot(1-z)^{1/2} = \bar{c} + O(|1-z|)
\]
by another application of 
Theorem~\ref{thm:hadamard}(b), this time with $k=1$. Hence
\[
L_1(z)\odot{}L_2(z) = \bigl[L_1(z)\odot{}L_2(z)\bigr]\Bigr\rvert_{z=1} +
\frac{\Gamma(\alpha-\tfrac12)}{\sqrt\pi}
(1-z)^{-\alpha+\tfrac12} + O(|1-z|).
\]
The other three Hadamard products are easily handled as
\begin{align*}
  L_1(z)\odot{}R_2(z) &= \bigl[L_1(z)\odot{}R_2(z)\bigr]\Bigr\rvert_{z=1} +
  O(|1-z|^{-\alpha+1}), \\
  L_2(z)\odot{}R_1(z) &= \bigl[L_2(z)\odot{}R_1(z)\bigr]\Bigr\rvert_{z=1} +
  O(|1-z|), \\
  R_1(z)\odot{}R_2(z) &= \bigl[R_1(z)\odot{}R_2(z)\bigr]\Bigr\rvert_{z=1} +
  O(|1-z|).
\end{align*}
Putting everything together, we get
\[ T(z)\odot\CAT(z/4) = C_0 +
\frac{\Gamma(\alpha-\tfrac12)}{\sqrt\pi}
(1-z)^{-\alpha+\tfrac12} + O(|1-z|^{-\alpha+1}).
\]
Using this in~\eqref{eq:7.1}, we get
\begin{equation}
  \label{eq:7.9}
  A(z)\odot\CAT(z/4) = C_0(1-z)^{-1/2} +
  \frac{\Gamma(\alpha-\tfrac12)}{\sqrt\pi} (1-z)^{-\alpha} +
  O(|1-z|^{-\alpha+\tfrac12}).
\end{equation}

To treat the case \(\alpha \geq 1/2\) we make use of the estimate
\begin{equation}
  \label{eq:7.21}
  (1-z)^{1/2} = \frac1{\Gamma(-1/2)}[ \Li_{3/2,0}(z) - \zeta(3/2) ]
  + O(|1-z|),
\end{equation}
a consequence of Theorem~\ref{theorem:Li-sing},
so that
\[
T(z) \odot (1-z)^{1/2} = \Li_{-\alpha,0}(z) \odot (1-z)^{1/2} =
\frac{1}{\Gamma(-1/2)} \Li_{\tfrac32-\alpha,0}(z) + R(z),
\]
where
\begin{equation}
  \label{eq:7.55}
  R(z) =
  \begin{cases}
    \bar{c} + O(|1-z|^{1-\alpha}) & 1/2 \leq \alpha < 1 \\
    O(|L(z)|)                     & \alpha=1            \\
    O(|1-z|^{1-\alpha})           & \alpha > 1.
  \end{cases}
\end{equation}
Hence
\[
T(z)\odot\CAT(z/4) = -\frac2{\Gamma(-1/2)}\Li_{\tfrac32-\alpha,0}(z) +
\widetilde{R}(z), 
\]
where $\widetilde{R}$, like $R$, satisfies~\eqref{eq:7.55} (with
a possibly different $\bar{c}$).  When \(\alpha=1/2\), this gives us
\[
T(z)\odot\CAT(z/4) = -\frac2{\Gamma(-1/2)}L(z) +
\bar{c} + O(|1-z|^{1/2}),
\]
so that
\begin{equation}
  \label{eq:7.7}
  A(z)\odot\CAT(z/4) =
  \frac1{\sqrt\pi}(1-z)^{-1/2}L(z) +
  \bar{c}(1-z)^{-1/2} + O(1).
\end{equation}
For \(\alpha > 1/2\) another singular expansion leads to the conclusion
that
\begin{equation}
  \label{eq:7.8}
  A(z)\odot\CAT(z/4) =
  \frac{\Gamma(\alpha-\tfrac12)}{\sqrt\pi}(1-z)^{-\alpha} + \widehat{R}(z),
\end{equation}
where
\begin{equation*}
  \widehat{R}(z) = 
  \begin{cases}
    O(|1-z|^{-\tfrac12})        & 1/2 < \alpha < 1 \\
    O(|1-z|^{-\tfrac12}|L(z)|)  & \alpha = 1       \\
    O(|1-z|^{-\alpha+\tfrac12}) & \alpha > 1.
  \end{cases}
\end{equation*}

We defer deriving the asymptotics of \(a_n\) until
Section~\ref{sec:asymptotics-moments}.

\section{Higher moments}
\label{sec:catalan_n_alpha_higher-moments}

We will analyze separately the cases \(0 < \alpha < 1/2\), \(\alpha=1/2\),
and \(\alpha > 1/2\).  The reason for this will become evident soon;
though the technique used to derive the asymptotics is induction in
each case, the induction hypothesis is different for each of these
cases.

\subsubsection{Small toll functions\texorpdfstring{ (${0 < \alpha < 1/2}$)}{}}
\label{sec:small-toll-functions}

We start by restricting ourselves to tolls of the form \(n^\alpha\)
where \(0 < \alpha < 1/2\).  In this case we observe that by
singularity analysis applied to~(\ref{eq:7.9}),
\begin{equation*}
  \frac{a_n \beta_n}{4^n} = \frac{C_0}{\sqrt\pi} n^{-1/2} + O(n^{-3/2}) +
  O(n^{\alpha-1}) = \frac{C_0}{\sqrt\pi} n^{-1/2} +
  O(n^{\alpha-1}),
\end{equation*}
so
\begin{equation*}
  a_n = n^{\tfrac32}[1 + O(n^{-1})][ C_0n^{-\tfrac12} +
  O(n^{\alpha-1})] = C_0 n + O(n^{\alpha+\tfrac12}) = (C_0 + o(1)) (n+1).
\end{equation*}
The lead order of
the mean $a_n = \E{X_n}$ does not depend on \(\alpha\).  We
perform an approximate centering to get to the dependence on \(\alpha\).


Define \(\widetilde{X}_n := X_n -
C_0(n+1)\), with \(X_0 := 0\); \(\tilde{\mu}_n(k) :=
\E{\widetilde{X}_n^k}\), with $\tilde{\mu}_n(0) = 1$ for all $n \geq
0$; and \( \hat{\mu}_n(k) :=
\beta_n\tilde{\mu}_n(k)/4^n \).  Let \(\widehat{M}_k(z)\) denote the ordinary
generating function of \(\hat{\mu}_n(k)\) in the argument $n$.

By an argument similar to the one that led to~(\ref{eq:7.35}), we get,
for \(k \geq 2\),
\[
\hat{\mu}_n(k) = \frac12 \sum_{j=1}^n \frac{\beta_{n-j}}{4^{n-j}}
\hat{\mu}_{j-1}(k) + \hat{r}_n(k), \qquad n \geq 1,
\]
where
\begin{align*}
\hat{r}_n(k) &:= \frac14\sum_{j=1}^n \sum_{\substack{k_1+k_2+k_3=k\\
    k_1,k_2 < k}} \binom{k}{k_1,k_2,k_3}
    \hat{\mu}_{j-1}(k_1)
    \hat{\mu}_{n-j}(k_2) t_n^{k_3}\\
    &= \frac14 \sum_{\substack{k_1+k_2+k_3=k\\
        k_1,k_2 < k}} \binom{k}{k_1,k_2,k_3} t_n^{k_3} \sum_{j=1}^n
    \hat{\mu}_{j-1}(k_1) \hat{\mu}_{n-j}(k_2),
\end{align*}
for \(n \geq 1\) and \(\hat{r}_0(k) := \hat{\mu}_0(k) = \tilde{\mu}_0(k) =
(-1)^kC_0^k\). Let \(\widehat{R}_k(z)\) denote the ordinary generating
function of \(\hat{r}_n(k)\) in the argument $n$.  Then,
mimicking~(\ref{eq:7.1}),
\begin{equation}
  \label{eq:7.4}
\widehat{M}_k(z) = \frac{\widehat{R}_k(z)}{\sqrt{1-z}}  
\end{equation}
with
\begin{equation}
  \label{eq:7.3}
\widehat{R}_k(z) = (-1)^kC_0^k + \sum_{\substack{k_1+k_2+k_3=k\\
    k_1,k_2<k}} \binom{k}{k_1,k_2,k_3} \bigl( T(z)^{\odot k_3}
\bigr) \odot
\bigl[\frac{z}4 \widehat{M}_{k_1}(z)\widehat{M}_{k_2}(z)\bigr],
\end{equation}
where for \(k\) a nonnegative integer
\[
T(z)^{\odot k} := \underbrace{T(z)\odot\cdots\odot T(z)}_{k}.
\]
Note that \(\widehat{M}_0(z) = \CAT(z/4)\).

We analyze separately the cases \(0 < \alpha \leq 1/4\) and \(1/4 <
\alpha < 1/2\).  The proof technique in either case is induction, though
the latter case is slightly simpler.

\begin{proposition}
  \label{thm:0alpha14}
  In the notation introduced above, when \(0 < \alpha \leq 1/4\), for
  $k \geq 1$,
  \[
  \widehat{M}_k(z) = C_k(1-z)^{-k(\alpha+\tfrac12) +\tfrac12} +
  O(|1-z|^{-k(\alpha+\tfrac{1}{2}) + \tfrac12 + (2\alpha-\epsilon)}),
  \]
  as $z \to 1$ in the sector \( -\pi+\epsilon' < \arg(1-z) <
  \pi-\epsilon' \), where
  \(\epsilon > 0\) and $\epsilon' > 0$ are arbitrarily small. The
  \(C_k\)'s here are defined by the  recurrence
  \begin{equation}
    \label{eq:7.10}
    C_k = \frac14 \sum_{j=1}^{k-1} \binom{k}{j} C_{j} C_{k-j} +
    kC_{k-1} \frac{ \Gamma(k\alpha+\tfrac{k}2-1)}{\Gamma((k-1)\alpha +
      \tfrac{k}2 -1)}, \quad k \geq 2; \qquad C_1 =
    \frac{\Gamma(\alpha-\tfrac12)}{\sqrt\pi}.
  \end{equation}
\end{proposition}
\begin{proof}
  For \(k=1\) the claim is true as shown in~\eqref{eq:7.9} with \(C_1\) as
  defined in~\eqref{eq:7.10}.  We will now analyze each term
  in~\eqref{eq:7.3} for \(k \geq 2\).  For notational convenience, define
  \(\alpha' := \alpha+\tfrac12\).  Also, observe that
  \[
  T(z)^{\odot{}k} = \Li_{-k\alpha,0}(z)
  = \Gamma(1+k\alpha)(1-z)^{-k\alpha-1} + O(|1-z|^{-k\alpha-\epsilon})
  \]
  by Lemma~\ref{lem:litoomz}.

  For this paragraph, consider the case that
  \(k_1\) and \(k_2\) are both nonzero.  It follows from the
  induction hypothesis that
  \begin{align*}
    \frac{z}4 \widehat{M}_{k_1}(z)\widehat{M}_{k_2}(z) &= \frac14(1-(1-z))
    \bigl[ C_{k_1}(1-z)^{-k_1\alpha'+\tfrac12} +
    O(|1-z|^{-k_1\alpha'+\tfrac12+(2\alpha-\epsilon)})
    \bigr]\\
    & \quad \times
    \bigl[ C_{k_2}(1-z)^{-k_2\alpha'+\tfrac12} +
    O(|1-z|^{-k_2\alpha'+\tfrac12+(2\alpha-\epsilon)})
    \bigr]\\
    &= \frac14C_{k_1}C_{k_2} (1-z)^{-(k_1+k_2)\alpha' + 1}
    + O(|1-z|^{-(k_1+k_2)\alpha'+1+(2\alpha-\epsilon)}).
  \end{align*}
  If \(k_3=0\) then the corresponding contribution to \(\widehat{R}_k(z)\) is
  \[
  \frac14\binom{k}{k_1}
  C_{k_1}C_{k_2} (1-z)^{-k\alpha' + 1}
    + O(|1-z|^{-k\alpha'+1+(2\alpha-\epsilon)}).
    \]
    If \(k_3 \ne 0\) we use Lemma~\ref{lem:omztoli} to express
    \begin{multline*}
       \frac{z}4 \widehat{M}_{k_1}(z)\widehat{M}_{k_2}(z) =
         \frac{C_{k_1}C_{k_2}}{4\Gamma((k_1+k_2)\alpha'-1)}
       \Li_{-(k_1+k_2)\alpha'+2,0}(z)\\
       +
       O(|1-z|^{-(k_1+k_2)\alpha'+1+(2\alpha-\epsilon)}) -
         \frac{C_{k_1}C_{k_2}}{4} 
       [(k_1+k_2)\alpha' <
         2]\frac{\zeta(-(k_1+k_2)\alpha'+2)}
         {\Gamma((k_1+k_2)\alpha'-1)}.
    \end{multline*}
  The corresponding contribution to \(\widehat{R}_{k}(z)\) is then
  \(\binom{k}{k_1,k_2,k_3}\) times:
  \begin{equation*}
  \frac{C_{k_1}C_{k_2}}{4\Gamma((k_1+k_2)\alpha'-1)}
  \Li_{-k\alpha'+\tfrac{k_3}2+2,0}(z)
  + \Li_{-k_3\alpha,0}(z)
  \odot O(|1-z|^{-(k_1+k_2)\alpha'+1+(2\alpha-\epsilon)}).
  \end{equation*}
  Now \(k_3 \leq k-2\) so \(-k\alpha'+\tfrac{k_3}2 + 2 < 1\).
  Hence the contribution when \(k_3 \ne 0\) is
  \[
  O(|1-z|^{-k\alpha'+\tfrac{k_3}2+1}) =
  O(|1-z|^{-k\alpha'+\tfrac32}) =
  O(|1-z|^{-k\alpha'+1+(2\alpha-\epsilon)}).
  \]

  Next we consider the case when \(k_1\) is nonzero but \(k_2=0\).  In
  this case using the induction hypothesis we see that
  \begin{align*}
    \frac{z}4 \widehat{M}_{k_1}(z)\widehat{M}_{k_2}(z) &=
    \frac{z}4 \CAT(z/4)\widehat{M}_{k_1}(z)\\
    &= \frac{1 - (1-z)^{1/2}}{2} \bigl[
    C_{k_1}(1-z)^{-k_1\alpha' + \tfrac12}\bigr]
    +
    O(|1-z|^{-k_1\alpha' + \tfrac12 +
      (2\alpha-\epsilon)})\\
    &= \frac{C_{k_1}}2 (1-z)^{-k_1\alpha'+\tfrac12} +
    O(|1-z|^{-k_1\alpha' + \tfrac12 + (2\alpha-\epsilon)}).
  \end{align*}
  Applying Lemma~\ref{lem:omztoli} to the last expression we get
  \begin{multline*}
  \frac{z}4 \widehat{M}_{k_1}(z)\widehat{M}_{k_2}(z) =
  \frac{C_{k_1}}{2\Gamma(k_1\alpha'-\tfrac12)}
  \Li_{-k_1\alpha'+\tfrac32,0}(z) \\
  + O(|1-z|^{-k_1\alpha' + \tfrac12 + (2\alpha-\epsilon)})
  - \frac{C_{k_1}}{2}[k_1\alpha'-\tfrac12 < 1]
  \frac{\zeta(-k_1\alpha'+\tfrac32)} 
    {\Gamma(k_1\alpha'-\tfrac12)}.
   \end{multline*}
   The contribution to \(\widehat{R}_{k}(z)\) is hence \(\binom{k}{k_1}\)
   times:
   \[
   \frac{C_{k_1}}{2\Gamma(k_1\alpha'-\tfrac12)}
   \Li_{-k\alpha'+\tfrac{k_3}2+\tfrac32,0}(z) +
   \Li_{-k_3\alpha,0}(z) \odot O(|1-z|^{-k_1\alpha' +
     \tfrac12 + (2\alpha-\epsilon)}).
   \]
   Using the fact that \(\alpha > 0\) and \(k_3 \leq k-1\), we conclude
   that \(-k\alpha'+\tfrac{k_3}2+\tfrac32 < 1\) so that, by
   Lemma~\ref{lem:litoomz} and Theorem~\ref{thm:hadamard}(a), the
   contribution is
   \[
   O(|1-z|^{-k\alpha'+\tfrac{k_3}2+\tfrac12}) =
   O(|1-z|^{-k\alpha'+\tfrac32})
   \]
   where the displayed equality holds unless \(k_3=1\).  When
   \(k_3=1\) we get a corresponding contribution to $\widehat{R}_k(z)$
   of \(\binom{k}{k-1}\) times:
   \begin{equation*}
   \frac{C_{k-1}\Gamma(k\alpha'-1)}
   {2\Gamma((k-1)\alpha'-\tfrac12)}
   (1-z)^{-k\alpha'+1} +
     O(|1-z|^{-k\alpha'+1+(2\alpha-\epsilon)}),
   \end{equation*}
   since for \(k \geq 2\) we have \(k\alpha' > 1 +
   (2\alpha-\epsilon)\).  The introduction of \(\epsilon\) handles the
   case when \(k\alpha' = 1 + 2\alpha\), which would have otherwise,
   according to Theorem~\ref{thm:hadamard}(c)
   introduced a logarithmic remainder.  In either case the
   remainder is $ O(|1-z|^{-k\alpha'+1+(2\alpha-\epsilon)})$.
   The case when \(k_2\) is nonzero but \(k_1=0\) is handled similarly by
   exchanging the roles of \(k_1\) and \(k_2\).
   
   The final contribution comes from the single term where both
   \(k_1\) and \(k_2\) are zero.  In this case the contribution to
   \(\widehat{R}_k(z)\) is, recalling~\eqref{eq:7.50},
   \begin{equation}
     \label{eq:7.5}
   \Li_{-k\alpha,0}(z) \odot [\frac{z}4 \CAT^2(z/4)] =
   \Li_{-k\alpha,0}(z) \odot 
   (\CAT(z/4)-1)= \Li_{-k\alpha,0}(z) \odot \CAT(z/4).
   \end{equation}
   Now, using Theorem~\ref{theorem:Li-sing},
   \begin{equation*}
     \CAT(z/4)              = 2-2(1-z)^{1/2} + O(|1-z|)
      = 2  +   2\frac{\zeta(3/2)}{\Gamma(-1/2)}
     - \frac{2}{\Gamma(-1/2)}\Li_{3/2,0}(z)  +
     O(|1-z|),
   \end{equation*}
   so that~\eqref{eq:7.5} is
   \[
   -\frac2{\Gamma(-1/2)} \Li_{\tfrac32-k\alpha,0}(z) +
   O(|1-z|^{1-k\alpha}) +
   \begin{cases}
     0                     & 1-k\alpha < 0,              \\
     O(|1-z|^{-\epsilon}) & 1-k\alpha=0,                \\
     O(1)                  & 1-k\alpha > 0.
   \end{cases}
   \]
   When \(\tfrac32-k\alpha < 1\) this is
   \(O(|1-z|^{-k\alpha+\tfrac12})\); when \(\tfrac32-k\alpha \geq 1\), it
   is \(O(1)\).  In either case we get a contribution
   which is \(O(|1-z|^{-k\alpha'+1+(2\alpha-\epsilon)})\).

   Hence
   \begin{align*}
     \widehat{R}_k(z) &= O(1) + \Biggl[ \sum_{\substack{k_1+k_2=k\\k_1,k_2<k}}
     \binom{k}{k_1} \frac{C_{k_1}C_{k_2}}4
     + 2k
     \frac{C_{k-1}}2
     \frac{\Gamma(k\alpha+\tfrac{k}2-1)}{\Gamma((k-1)\alpha+\tfrac{k}2-1)}
     \Biggr] (1-z)^{-k\alpha'+1}\\
     & \qquad\qquad + O(|1-z|^{-k\alpha'+1+(2\alpha-\epsilon)})\\
     &= C_k(1-z)^{-k\alpha'+1} + O(|1-z|^{-k\alpha'
       + 1 + (2\alpha-\epsilon)}),
   \end{align*}
   with the \(C_k\)'s defined by the recurrence~\eqref{eq:7.10}.  Now
   using~\eqref{eq:7.4}, the claim follows.
 \end{proof}

 Next, we handle the case when \(1/4 < \alpha < 1/2\).
\begin{proposition}
  \label{thm:14alpha12}
  In the notation introduced above Proposition~\ref{thm:0alpha14}, when
  \(1/4 < \alpha < 1/2\),
  \[
  \widehat{M}_k(z) = C_k(1-z)^{-k(\alpha+\tfrac12) +\tfrac12} +
  O(|1-z|^{-k(\alpha+\tfrac{1}{2}) + 1}),
  \]
  for \(k \geq 1\), where the \(C_k\)'s are defined by  the
  recurrence~\eqref{eq:7.10}.
\end{proposition}
\begin{proof}
  For \(k=1\) the claim is true as shown in
  Section~\ref{sec:cat_n_alpha_asymptotics-mean} with \(C_1\) defined
  at~\eqref{eq:7.10}.  We will now analyze each term in~\eqref{eq:7.3} for
  \(k \geq 2\) in a manner similar to that in the proof of
  Proposition~\ref{thm:0alpha14}.
  
  First we consider the case that \(k_1\) and \(k_2\) are both nonzero. It
  follows from the induction hypothesis that
  \begin{equation*}
    \frac{z}4 \widehat{M}_{k_1}(z)\widehat{M}_{k_2}(z) =
    \frac{C_{k_1}C_{k_2}}4 (1-z)^{-(k_1+k_2)\alpha' + 1} 
    + O(|1-z|^{-(k_1+k_2)\alpha'+\tfrac32}).
  \end{equation*}
  If \(k_3=0\) then the corresponding contribution to \(\widehat{R}_k(z)\) is
  \[
  \frac14 \binom{k}{k_1} C_{k_1}C_{k_2} (1-z)^{-k\alpha' + 1} +
  O(|1-z|^{-k\alpha'+\tfrac32}).
    \]
    If \(k_3 \ne 0\) we use Lemma~\ref{lem:omztoli} to express
    \begin{multline*}
       \frac{z}4 \widehat{M}_{k_1}(z)\widehat{M}_{k_2}(z) =
         \frac{C_{k_1}C_{k_2}}{4\Gamma((k_1+k_2)\alpha'-1)}
       \Li_{-(k_1+k_2)\alpha'+2,0}(z)\\
       +
       O(|1-z|^{-(k_1+k_2)\alpha'+\tfrac32}) - \frac{C_{k_1}C_{k_2}}{4}
       [(k_1+k_2)\alpha' < 2]
         \frac{\zeta(-(k_1+k_2)\alpha'+2)}
         {\Gamma((k_1+k_2)\alpha'-1)}.
       \end{multline*}
  The contribution to \(\widehat{R}_{k}(z)\) is then
  \(\binom{k}{k_1,k_2,k_3}\) times:
  \begin{equation*}
  \frac{C_{k_1}C_{k_2}}{4\Gamma((k_1+k_2)\alpha'-1)}
  \Li_{-k\alpha'+\tfrac{k_3}2+2,0}(z)
  + \Li_{-k_3\alpha,0}(z)
  \odot O(|1-z|^{-(k_1+k_2)\alpha'+\tfrac32}).
  \end{equation*}
  Now \(k_3 \leq k-2\) so \(-k\alpha'+\tfrac{k_3}2 + 2 < 1\).
  Hence the contribution when \(k_3 \ne 0\) is
  \[
  O(|1-z|^{-k\alpha'+\tfrac{k_3}2+1}) =
  O(|1-z|^{-k\alpha'+\tfrac32})
  \]
  by another application of Theorem~\ref{thm:hadamard}(a).

  Next we consider the case when \(k_1\) is nonzero but \(k_2=0\).  In
  this case using the induction hypothesis we see that
  \begin{equation*}
    \frac{z}4 \widehat{M}_{k_1}(z)\widehat{M}_{k_2}(z) = \frac{z}4
    \CAT(z/4)\widehat{M}_{k_1}(z)
    = \frac{C_{k_1}}2 (1-z)^{-k_1\alpha'+\tfrac12} +
    O(|1-z|^{-k_1\alpha' + 1}).
  \end{equation*}
  Applying Lemma~\ref{lem:omztoli} to the last expression we get
  \begin{multline*}
  \frac{z}4 \widehat{M}_{k_1}(z)\widehat{M}_{k_2}(z) =
  \frac{C_{k_1}}{2\Gamma(k_1\alpha'-\tfrac12)}
  \Li_{-k_1\alpha'+\tfrac32,0}(z) \\
  + O(|1-z|^{-k_1\alpha' + 1})
  - \frac{C_{k_1}}{2} [k_1\alpha' - \tfrac12 < 1]
    \frac{\zeta(-k_1\alpha'+\tfrac32)}
    {\Gamma(k_1\alpha'-\tfrac12)}
   \end{multline*}
   The contribution to \(\widehat{R}_{k}(z)\) is hence
   \(\binom{k}{k_1}\) times:
   \[
   \frac{C_{k_1}}{2\Gamma(k_1\alpha'-\tfrac12)}
   \Li_{-k\alpha'+\tfrac{k_3}2+\tfrac32,0}(z) +
   \Li_{-k_3\alpha,0}(z) \odot O(|1-z|^{-k_1\alpha' +
     1}).
   \]
   Using the fact that \(\alpha > 1/4\) and \(k_3 \leq k-1\), we conclude
   that \(-k\alpha'+\tfrac{k_3}2+\tfrac32 < \tfrac12\) so that the
   contribution is
   \[
   O(|1-z|^{-k\alpha'+\tfrac{k_3}2+\tfrac12}) =
   O(|1-z|^{-k\alpha'+\tfrac32})
   \]
   by Theorem~\ref{thm:hadamard}(a),
   unless \(k_3=1\).  When \(k_3=1\) we get a contribution of
   \begin{equation*}
     \binom{k}{k-1}
   \frac{C_{k-1}\Gamma(k\alpha'-1)}
   {2\Gamma((k-1)\alpha'-\tfrac12)}
   (1-z)^{-k\alpha'+1} +
     O(|1-z|^{-k\alpha'+\tfrac32}),
   \end{equation*}
   since \(k \geq 2\) and \(\alpha > 1/4\) imply \(k\alpha' > 3/2\).
   The case when \(k_2\) is nonzero but \(k_1=0\) is handled similarly.
   
   The final contribution comes from the single term where both \(k_1\)
   and \(k_2\) are zero.  In this case the contribution to
   \(\widehat{R}_k(z)\) is
   \begin{equation*}
     \label{eq:6}
   \Li_{-k\alpha,0}(z) \odot [\frac{z}4 \CAT^2(z/4)] =
   \Li_{-k\alpha,0}(z) \odot 
   (\CAT(z/4)-1)= \Li_{-k\alpha,0}(z) \odot \CAT(z/4)
   \end{equation*}
   which is easily seen to be \(O(|1-z|^{-k\alpha'+\tfrac32})\).

   Hence, as in the proof of Proposition~\ref{thm:0alpha14}, we see that
   \[
   \widehat{R}_k(z) = C_k(1-z)^{-k(\alpha+\tfrac{1}2)+1} +
       O(|1-z|^{-k(\alpha+\tfrac{1}2) 
       + \tfrac32})
   \]
   with the \(C_k\)'s defined by the recurrence~\eqref{eq:7.10}, and the
   claim follows using~\eqref{eq:7.4}.
 \end{proof}

\subsubsection{Large toll functions\texorpdfstring{ ($ \alpha \geq
    1/2$)}{}}
\label{sec:large-toll-functions}
When \(\alpha \geq 1/2\) there is no need to apply the centering
techinques.  Define \({\mu}_n(k) := \E{{X}_n^k}\) and
\( \bar{\mu}_n(k) := \beta_n{\mu}_n(k)/4^n \).  Let
\(\overline{M}_k(z)\) denote the 
ordinary generating function of \(\bar{\mu}_n(k)\) in $n$.  Observe that
\(\overline{M}_0(z) = \CAT(z/4)\).  As earlier, conditioning on the
key stored at the root, we get, for \(k \geq 2\),
\[
\bar{\mu}_n(k) = \frac12 \sum_{j=1}^n \frac{\beta_{n-j}}{4^{n-j}}
  \bar{\mu}_{j-1}(k) + \bar{r}_n(k),
\qquad n \geq 1,
\]
where
\begin{equation*}
\bar{r}_n(k) := 
     \frac14 \sum_{\substack{k_1+k_2+k_3=k\\
        k_1,k_2 < k}} \binom{k}{k_1,k_2,k_3} t_n^{k_3} \sum_{j=1}^n
    \bar{\mu}_{j-1}(k_1) \bar{\mu}_{n-j}(k_2),
\end{equation*}
for \(n \geq 1\) and \(\bar{r}_0(k) := \bar{\mu}_0(k) = \mu_0(k) = 0\).
Let \(\overline{R}_k(z)\) denote the ordinary generating
function of \(\bar{r}_n(k)\) in $n$.  Then
\begin{equation*}
\overline{M}_k(z) = \frac{\overline{R}_k(z)}{\sqrt{1-z}}  
\end{equation*}
and
\begin{equation}
\label{eq:7.12}
\overline{R}_k(z) = \sum_{\substack{k_1+k_2+k_3=k\\
    k_1,k_2<k}} \binom{k}{k_1,k_2,k_3} \bigl( T(z)^{\odot k_3}
\bigr) \odot
\bigl[ \frac{z}4 \overline{M}_{k_1}(z)\overline{M}_{k_2}(z)\bigr].
\end{equation}

We can now state the result about the asymptotics of the generating
function~\(\overline{M}_k\) when \(\alpha > 1/2\).  The case  \(\alpha=1/2\)
will be handled subsequently, in Proposition~\ref{thm:alpha=12}.
\begin{proposition}
  \label{thm:12alpha}
  Let $\epsilon > 0$ be arbitrary, and define
  \begin{equation}
    \label{eq:7.56}
    c :=
    \begin{cases}
      \alpha - \frac12      & \frac12 < \alpha < 1 \\
      \frac12 - \epsilon & \alpha = 1       \\
      \frac12            & \alpha > 1.
    \end{cases}
  \end{equation}
  Then the generating function $\overline{M}_k(z)$
  of~\(\bar{\mu}_n(k)\)  satisfies
  \[
  \overline{M}_k(z) = C_k(1-z)^{-k(\alpha+\tfrac12)+\tfrac12} +
  O(|1-z|^{-k(\alpha+\tfrac12)+ \tfrac12 + c})
  \]
  for \(k \geq 1\), where the \(C_k\)'s are defined by the
  recurrence~\eqref{eq:7.10}.
\end{proposition}
\begin{proof}
  The proof is very similar to those of Propositions~\ref{thm:0alpha14}
  and~\ref{thm:14alpha12}.  We present a sketch.  The reader is
  invited to compare the cases enumerated below to those in the
  earlier proofs.
  
  When \(k=1\) the claim is true by~\eqref{eq:7.8}.  We analyze the
  various terms in~\eqref{eq:7.12} for \(k \geq 2\), employing the
  notational convenience \(\alpha' := \alpha+\tfrac12\). 

  When both \(k_1\) and \(k_2\) are nonzero then the contribution to
  \(\overline{R}_k(z)\) is
  \[
  \frac14\binom{k}{k_1} C_{k_1}C_{k_2} (1-z)^{-k\alpha'+1} +
  O(|1-z|^{-k\alpha'+c+1})
  \]
  when \(k_3=0\) and is \(O(|1-z|^{-k\alpha'+c+1})\) otherwise.

  When \(k_1\) is nonzero and \(k_2=0\) the contribution to
  \(\overline{R}_k(z)\)
  is
  \[
  k \frac{C_{k-1} \Gamma(k\alpha'-1)}{2\Gamma((k-1)\alpha'-\tfrac12)}
  (1-z)^{-k\alpha'+1} + O(|1-z|^{-k\alpha'+c+1})
  \]
  when \(k_3=1\) and \(O(|1-z|^{-k\alpha'+c+1})\) otherwise.  The
  case when \(k_2\) is nonzero and \(k_1=0\) is identical.

  The final contribution comes from the single term when both \(k_1\)
  and \(k_2\) are zero.  In this case we get a contribution of
  \(O(|1-z|^{-k\alpha+\tfrac12})\) which is
  \(O(|1-z|^{-k\alpha'+c+1})\).  Adding all these contributions
  yields the desired result.
\end{proof}

The result when \(\alpha=1/2\) is as follows.  Recall that $L(z) :=
\ln((1-z)^{-1})$.
\begin{proposition}
  \label{thm:alpha=12}
  Let $\alpha=1/2$.
  In the notation of Proposition~\ref{thm:12alpha},
  \[
  \overline{M}_k(z) = (1-z)^{-k+\tfrac12} \sum_{l=0}^k C_{k,l}
  L^{k-l}(z) + O(|1-z|^{-k+1-\epsilon})
  \]
  for \(k \geq 1\) and any \(\epsilon > 0\), where the $C_{k,l}$'s are
  constants.  The constant multiplying the
  lead-order term is given by
  \begin{equation}
    \label{eq:7.11}
    C_{k,0} = \frac{(2k-2)!}{2^{2k-2}(k-1)!\pi^{k/2}}.
  \end{equation}
\end{proposition}
\begin{proof}
  For \(k=1\) the claim is true by~\eqref{eq:7.7}.  For \(k \geq 2\) we
  estimate each term in~\eqref{eq:7.12}.  

  For this paragraph, suppose that \(k_1\) and \(k_2\) are both nonzero.
  Then using the induction hypothesis we get
  \begin{multline}
    \label{eq:7.13}
    \frac{z}4 \overline{M}_{k_1}(z)\overline{M}_{k_2}(z) = \frac14
    (1-z)^{-(k_1+k_2)+1}
    \sum_{l=0}^{k_1} \sum_{m=0}^{k_2} C_{k_1,l}C_{k_2,m}
    L^{(k_1+k_2)-(l+m)}(z)
     \\ + O(|1-z|^{-(k_1+k_2)+\tfrac32-2\epsilon}).
  \end{multline}
  If \(k_3=0\) the corresponding contribution to \(\overline{R}_k(z)\) is
  \[
  \frac14\binom{k}{k_1} (1-z)^{-k+1} \sum_{l=0}^k A_{k,l}
  L^{k-l}(z) + O(|1-z|^{-k+\tfrac32-2\epsilon}),
  \]
  where \(A_{k,0} = C_{k_1,0} C_{k-k_1,0}\) and \((A_{k,l})_{l=1}^k\)
  can be determined easily from~\eqref{eq:7.13}.  If \(k_3 \ne 0\) then we
  use Lemma~\ref{lem:omztoli} once again to get
  \[
  \frac{z}4 \overline{M}_{k_1}(z)\overline{M}_{k_2}(z) = \sum_{l=0}^{k_1+k_2}
  \widetilde{A}_{k_1+k_2,l} \Li_{-(k_1+k_2)+2,k_1+k_2-l}(z) +
  O(|1-z|^{-(k_1+k_2)+\tfrac32-2\epsilon}),
  \]
  where \((\widetilde{A}_{k_1+k_2,l})_{l=0}^{k_1+k_2}\) can be determined
  using~\eqref{eq:7.13} and Lemma~\ref{lem:omztoli}.  The corresponding
  contribution to \(\overline{R}_k(z)\) is then \(\binom{k}{k_1,k_2,k_3}\)
  times:
  \[
  \sum_{l=0}^{k_1+k_2} \widetilde{A}_{k_1+k_2,l}
  \Li_{-k+\tfrac{k_3}2+2,k_1+k_2-l}(z) + \Li_{-\tfrac{k_3}2,0}(z)
  \odot O(|1-z|^{-(k_1+k_2)+\tfrac32-2\epsilon}).
  \]
  Invoking Lemma~\ref{lem:litoomz} we see that the contribution when
  \(k_3 \ne 0\) is \(O(|1-z|^{-k+\tfrac{k_3}2+1-\epsilon})\), which is
  \(O(|1-z|^{-k+\frac32-\epsilon})\).

  Next suppose that \(k_1\) is nonzero but \(k_2=0\).  In this case using
  the induction hypothesis and we see that
  \[
  \frac{z}4 \overline{M}_{k_1}(z)\overline{M}_{k_2}(z) = \frac12
  (1-z)^{-k_1+\tfrac12} 
  \sum_{l=0}^{k_1} C_{k_1,l} L^{k_1-l}(z) +
  O(|1-z|^{-k_1+1-\epsilon}).
  \]
  By Lemma~\ref{lem:omztoli} this is
  \[
  \sum_{l=0}^{k_1} B_{k_1,l} \Li_{-k_1+\tfrac32,k_1-l}(z) +
  O(|1-z|^{-k_1+1-\epsilon}),
  \]
  where \((B_{k_1,l})_{l=0}^{k_1}\) are constants.
  The corresponding contribution to $\overline{R}_k(z)$ is then
  \(\binom{k}{k_1}\) times:
  \[
  \sum_{l=0}^{k_1} B_{k_1,l}
  \Li_{-k+\tfrac{k_3}2+\tfrac32,k_1-l}(z) + \Li_{-\tfrac{k_3}2,0}(z)
  \odot O(|1-z|^{-k_1+1-\epsilon}).
  \]
  When \(k_3 \geq 2\) the contribution is
  \(O(|1-z|^{-k+\tfrac{k_3}2+\tfrac12-\epsilon})\), which is
  \(O(|1-z|^{-k+\tfrac32-\epsilon})\).  When \(k_3=1\) the contribution
  to \(\overline{R}_k(z)\) is
  \[
  (1-z)^{-k+1}\sum_{l=0}^{k-1}
  \widetilde{B}_{k-1,l}L^{k-1-l}(z) + O(|1-z|^{-k+\tfrac32-\epsilon}),
  \]
  for some suitably defined constants $\widetilde{B}_{k,l}$.
  Note that there is no lead-order contribution from this term.  The
  terms where \(k_1=0\) and \(k_2\) is nonzero are similarly handled.

  The final contribution from the single term with \(k_1=k_2=0\) is
  \(O(|1-z|^{-\tfrac{k}{2}+\tfrac12})\), which is
  \(O(|1-z|^{-k+\tfrac32})\) for \(k \geq 2\).

  Summing all the contributions we get
  \begin{equation*}
    \overline{R}_{k}(z) = (1-z)^{-k+1} \sum_{l=0}^k C_{k,l}
    L^{k-l}(z) + O(|1-z|^{-k+\frac32-\epsilon}),
  \end{equation*}
  where, for \(k \geq 2\),
  \begin{align*}
    C_{k,0} &= \frac14 \sum_{k_1=1}^{k-1}\binom{k}{k_1}
    C_{k_1,0}C_{k-k_1,0}, \qquad C_{1,0} = \frac1{\sqrt\pi},\\ 
    C_{k,l} &= \frac14 \sum_{k_1=1}^{k-1}\binom{k}{k_1} A_{k,l} +
    2\widetilde{B}_{k-1,l}, \qquad l \geq 1.
  \end{align*}
  The recurrence for \(C_{k,0}\) can be easily solved to yield
  \begin{equation*}
    \label{eq:7.14}
    C_{k,0} = 2\cdot k!\:[z^k]
    \left[1-\left(1-\frac{z}{\sqrt{\pi}} \right)^{1/2} \right] =
    \frac{(2k-2)!}{2^{2k-2}(k-1)!\pi^{k/2}}.
  \end{equation*}
  The result follows immediately.
\end{proof}

\section{Asymptotics of moments}
\label{sec:asymptotics-moments}
For \(0 < \alpha< 1/2 \), we have seen in Propositions~\ref{thm:0alpha14}
and~\ref{thm:14alpha12} that the generating function $\widehat{M}_k(z)$ of
\( \hat{\mu}_n(k) = \beta_n\tilde{\mu}_n(k)/4^n \) satisfies
\[
\widehat{M}_k(z) = C_k(1-z)^{-k(\alpha+\tfrac{1}2)+\tfrac12} +
O(|1-z|^{-k(\alpha+\tfrac{1}{2})+\tfrac12 + c}),
\]
where \(c := \min\{2\alpha-\epsilon,1/2\}\).
By Theorem~\ref{theorem:transfer},
\[
\frac{\beta_n\tilde{\mu}_n(k)}{4^n} =
C_k
\frac{n^{k(\alpha+\tfrac12)-\tfrac32}}{\Gamma(k(\alpha+\tfrac12)-\tfrac12)}
  + O(n^{k(\alpha+\tfrac{1}{2})-\tfrac32-c}).
\]
Recall that
\[
\beta_n = 
\frac{4^n}{\sqrt{\pi}n^{3/2}} \left(1 + O(\tfrac1n)\right),
\]
so that
\begin{equation}
  \label{eq:7.15}
  \tilde{\mu}_n(k)
  = \frac{C_k\sqrt{\pi}}{\Gamma(k(\alpha+\tfrac{1}2)-\tfrac12)}
  n^{k(\alpha+\tfrac{1}2)} +
  O(n^{k(\alpha+\tfrac{1}2) - c}).
\end{equation}

For \(\alpha > 1/2\) a similar analysis using
Proposition~\ref{thm:12alpha} yields
\begin{equation}
  \label{eq:7.16}
  \mu_n(k) = \frac{C_k\sqrt{\pi}}{\Gamma(k(\alpha+\tfrac{1}2)-\tfrac12)}
  n^{k(\alpha+\tfrac{1}2)} +
  O(n^{k(\alpha+\tfrac{1}2) - c}),
\end{equation}
with now $c$ as defined at~\eqref{eq:7.56}.  Finally, when \(\alpha=1/2\) the
asymptotics of the moments are given by
\begin{equation}
  \label{eq:7.17}
  \mu_n(k) =
  \left( \frac{1}{\sqrt\pi} \right)^k(n\ln{n})^k +
  O(n^k(\ln{n})^{k-1}).
\end{equation}

\section{The limiting distributions}
\label{sec:limit-distr}
In Section~\ref{sec:alpha-ne-12} we will use our moment
estimates~\eqref{eq:7.15} and \eqref{eq:7.16} with the method of moments
to derive limiting distributions for our additive functionals.  The
case $\alpha=1/2$ requires a somewhat delicate analysis, which we will
present separately in Section~\ref{sec:alpha=12}.

\subsection{\texorpdfstring{$\alpha \ne 1/2$}{alpha not 1/2}}
\label{sec:alpha-ne-12}

We first handle the case \(0 < \alpha < 1/2\). (We assume this
restriction until just before
Proposition~\ref{thm:limit_law_alpha_ne_12}.) We have
\begin{equation}
  \label{eq:7.39}
  \tilde\mu_n(1) = \E\widetilde{X}_n = \E[X_n - C_0(n+1)]
  = \frac{C_1\sqrt\pi}{\Gamma(\alpha)} n^{\alpha+\tfrac12} +
  O(n^{\alpha+\tfrac12-c})
\end{equation}
with $c := \min\{2\alpha-\epsilon,1/2\}$ and
\begin{equation*}
  \tilde\mu_n(2) = \E\widetilde{X}_n^2 =
  \frac{C_2\sqrt\pi}{\Gamma(2\alpha+\tfrac12)} n^{2\alpha+1} +
  O(n^{2\alpha+1-c}).
\end{equation*}
So
\begin{equation}\label{eq:7.40}
  \Var{X_n} = \Var{\widetilde{X}_n} = \tilde\mu_n(2) -
  [\tilde\mu_n(1)]^2 = \sigma^2 n^{2\alpha+1} + O(n^{2\alpha+1-c}),
\end{equation}
where
\begin{equation}\label{eq:7.46}
  \sigma^2 := \frac{C_2\sqrt\pi}{\Gamma(2\alpha+\tfrac12)} -
  \frac{C_1^2\pi}{\Gamma^2(\alpha)}.
\end{equation}
We also have, for \(k \geq 1\),
\begin{equation}\label{eq:7.41}
  \E \left[ \frac{\widetilde{X}_n}{n^{\alpha+\tfrac12}} \right]^k
  = \frac{\tilde{\mu}_n(k)}{n^{k(\alpha+\tfrac12)}}
  = \frac{C_k\sqrt\pi}{\Gamma(k(\alpha+\tfrac12)-\tfrac12)} + O(n^{-c}).
\end{equation}

The following lemma provides a sufficient bound on the moments
facilitating the use of the method of moments.
\begin{lemma}
  \label{lem:Ckbound}
  Define \(\alpha' := \alpha+\tfrac12\).  There exists a constant \(A
  < \infty \)
  depending only on \(\alpha\) such that
  \begin{equation*}
    \left|\frac{C_k}{k!}\right| \leq A^k k^{\alpha'k}
  \end{equation*}
  for all \(k \geq 1\).
\end{lemma}
\begin{proof}
  The proof is by induction.  For the basis case $k=1$, we simply choose
  \(A\) large enough 
  so that the claim is true.  For the induction step, assume $k \geq
  2$.  Define \(s_k := {C_k}/{k!}\).
  Dividing~\eqref{eq:7.10} by \(k!\) we get
  \begin{equation}
    \label{eq:7.18}
    s_k = \frac14 \sum_{j=1}^{k-1} s_j s_{k-j} + s_{k-1}
    \frac{\Gamma(k\alpha'-1)}{\Gamma(k\alpha'-1-\alpha)}.
  \end{equation}
  Using Stirling's approximation we can find a constant~\(\gamma < \infty\)
  depending only on \(\alpha\) so that, for all $k \geq 2$,
  \begin{equation*}
    \left|\frac{\Gamma(k\alpha'-1)}{\Gamma(k\alpha'-1-\alpha)}\right|
    \leq \gamma k^\alpha.
  \end{equation*}
  Using this in~\eqref{eq:7.18} we get
  \begin{equation*}
    |s_k| \leq \frac14 \sum_{j=1}^{k-1} |s_j| |s_{k-j}| + \gamma |s_{k-1}|
     k^\alpha
   \end{equation*}
   By the induction hypothesis
   \begin{equation*}
     |s_k| \leq \frac{A^k}{4} \sum_{j=1}^{k-1} [j^j (k-j)^{k-j}
      ]^{\alpha'} + \gamma A^{k-1} (k-1)^{\alpha'(k-1)} k^\alpha.
    \end{equation*}
    Now \(x^x (k-x)^{k-x}\) decreases as \(x\) increases for $0 < x
    \leq k/2$ so that we can bound the sum above by the sum of two
    times the $j=1$ term and $k$ times the $j=2$ term.  This and
    trivial bounds yield
    \begin{align*}
      |s_k| &\leq \frac{A^k}2 \left[ (k-1)^{k-1} + \frac{k}2 \times
      4(k-2)^{k-2} \right]^{\alpha'} + \gamma A^{k-1} k^{\alpha'k-\tfrac12}\\
      & \leq \frac{A^k}{2} [ k^{k-1} + 2k^{k-1} ]^{\alpha'} + A^k
      \frac\gamma{A} k^{\alpha'k-\tfrac12}\\
      & = \frac{A^k}2 (3k^{k-1})^{\alpha'} + A^k
      \frac\gamma{A} k^{\alpha'k-\tfrac12} \\
      &= \left[ \frac{3^{\alpha'}}2 k^{-\alpha'} + \frac{\gamma}A k^{-\tfrac12}
      \right] A^k k^{\alpha'k} \leq A^k k^{\alpha'k},
    \end{align*}
    the last inequality justified for all $k \geq 2$ by choosing \(A\)
    large enough.  This proves the claim.
\end{proof}
It follows from Lemma~\ref{lem:Ckbound} and Stirling's approximation that
\begin{equation}
  \label{eq:7.19}
  \left|\frac{C_k\sqrt\pi}{k! \Gamma(k(\alpha+\tfrac12)-\tfrac12)}\right|
  \leq B^k 
\end{equation}
for large enough \(B\) depending only on $\alpha$.  Using standard
arguments~\cite[Theorem 30.1]{MR95k:60001} it follows that \(X_n\)
suitably normalized has a limiting distribution that is characterized
by its moments.  Before we state the result, we observe that the
argument presented above can be adapted with minor modifications to
treat the case \(\alpha > 1/2\), with \(\widetilde{X}_n\) replaced by
\(X_n\).  We can now state a result for $\alpha \ne 1/2$.  We will use
the notation $\stackrel{\mathcal{L}}{\to}$ to denote convergence in
law (or distribution).

\begin{proposition}\label{thm:limit_law_alpha_ne_12}
  Let $X_n$ denote the additive functional on Catalan trees
  induced by the toll sequence $(n^\alpha)_{n \geq 0}$.  Define the
  random variable $Y_n$ as follows:
  \begin{equation*}
    Y_n := 
    \begin{cases}
      \displaystyle\frac{X_n-C_0(n+1)}{n^{\alpha+\tfrac12}} & 0 <
      \alpha < 1/2,\\
      \displaystyle\frac{X_n}{n^{\alpha+\tfrac12}} & \alpha > 1/2,
    \end{cases}
  \end{equation*}
  where
  \begin{equation*}
    C_0 := \sum_{n=0}^\infty n^\alpha \frac{\beta_n}{4^n}, \qquad
    \beta_n = \frac{1}{n+1} \binom{2n}{n}.
  \end{equation*}
  Then
  \begin{equation*}
    Y_n \stackrel{\mathcal{L}}{\to} Y;
  \end{equation*}
  here $Y$ is a random variable with the unique distribution whose
  moments are
  \begin{equation}\label{eq:7.49}
    \E Y^k = \frac{C_k
      \sqrt\pi}{\Gamma(k(\alpha+\tfrac12)-\tfrac12)}, 
  \end{equation}
  where the $C_k$'s satisfy the recurrence
  \begin{equation*}
    C_k = \frac{1}{4} \sum_{j=1}^{k-1} \binom{k}{j} C_j C_{k-j} + k
    \frac{\Gamma(k\alpha + \tfrac{k}2 -1)}{\Gamma((k-1)\alpha +
    \tfrac{k}2 - 1)} C_{k-1}, k \geq 2; \quad C_1 =
    \frac{\Gamma(\alpha-\tfrac12)}{\sqrt\pi}.
  \end{equation*}
\end{proposition}
The case $\alpha=1/2$ is handled in Section~\ref{sec:alpha=12},
leading to Proposition~\ref{prop_alpha_12}, and a unified result for
all cases is stated as Theorem~\ref{thm:limit-dist}.

\begin{remark}\label{remark_y_alpha_properties}
We now consider some properties of the limiting random variable $Y
\equiv Y(\alpha)$ defined by its moments at~\eqref{eq:7.49} for $\alpha
\ne 1/2$.
\begin{enumerate}[(a)]
\item 
  When \( \alpha=1 \) setting $\Omega_k := C_k/2$ we see immediately
  that 
  \begin{equation*}
    \E Y^k =  \frac{-\Gamma(-1/2)}{\Gamma((3k-1)/2)}\Omega_k,
  \end{equation*}
  where
  \begin{equation*}
    2\Omega_k = \sum_{j=1}^{k-1} \binom{k}{j} \Omega_j \Omega_{k-j} + k(3k - 4)
    \Omega_{k-1}, \qquad \Omega_1 = \frac12.
  \end{equation*}
  Thus \(Y\) has the ubiquitous Airy distribution and we have
  recovered the limiting distribution of path length in Catalan
  trees~\cite{MR92m:60057,MR92k:60164,MR93e:60175}.  The Airy
  distribution arises in many contexts including parking allocations,
  hashing tables, trees, discrete random walks, mergesorting,
  etc.---see, e.g., the introduction of~\cite{MR2002j:68115}
  which contains numerous references to the Airy distribution.
\item When $\alpha = 2$, setting $\eta := Y/\sqrt2$ and $a_{0,l} :=
  2^{2l-1} C_l$, we see that
  \begin{equation*}
    \E{\eta^l} = \frac{\sqrt\pi}{2^{(5l-2)/2} \Gamma((5l-1)/2)} a_{0,l},
  \end{equation*}
  where
  \begin{equation*}
    a_{0,l} = \frac{1}{2} \sum_{j=1}^{l-1} \binom{l}{j} a_{0,j}
    a_{0,l-j} + l(5l-4)(5l-6), \qquad a_{0,1} = 1.
  \end{equation*}
  We have thus recovered the recurrence for the moments of the
  distribution $\mathcal{L}({\eta})$, which arises in the study of the
  Wiener index  of Catalan trees~\cite[proof of Theorem~3.3 in
  Section~5]{janson:_wiener}.
\item Consider the variance $\sigma^2$ defined at~\eqref{eq:7.46}.
  \begin{enumerate}[(i)]
  \item   Figure~\ref{fig:variance}, plotted using
    \texttt{Mathematica}, suggests that $\sigma^2$ is positive for
    all $\alpha > 0$.  We will prove this fact in
    Theorem~\ref{thm:limit-dist}.
    \begin{figure}[htbp]
      \centering
      \includegraphics{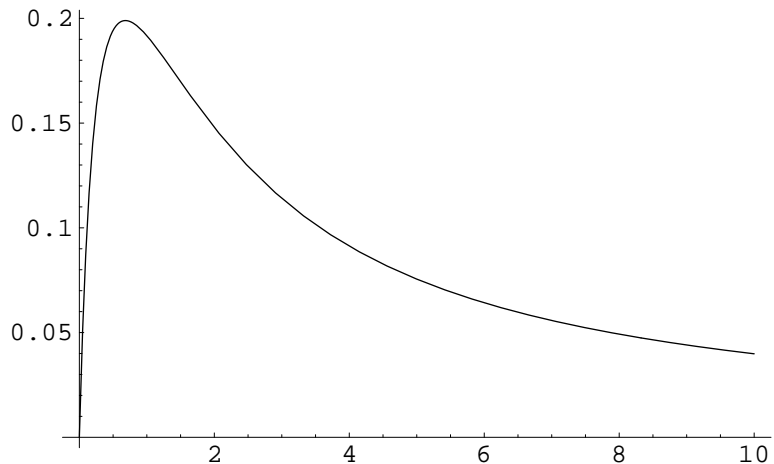}
      \caption{$\sigma^2$ of~\eqref{eq:7.46} as a function of $\alpha$.}
      \label{fig:variance}
    \end{figure}
    There is also numerical evidence that $\sigma^2$ is unimodal with
    $\max_{\alpha} \sigma^2(\alpha) \doteq 0.198946$ achieved at $\alpha \doteq
    0.682607$. (Here $\doteq$ denotes approximate equality.)
  \item As $\alpha \to \infty$, using Stirling's approximation one can
    show that $\sigma^2 \sim (\sqrt{2}-1)\alpha^{-1}$.
  \item As $\alpha \downarrow 0$, using a Laurent series expansion of
    $\Gamma(\alpha)$ we see that $\sigma^2 \sim 4(1-\ln2) \alpha$.
  \item Though the random variable $Y(\alpha)$ has been defined only
    for $\alpha \ne 1/2$, the variance $\sigma^2$ has a limit at $\alpha=1/2$:
    \begin{equation}\label{eq:7.51}
      \lim_{\alpha \to 1/2} \sigma^2(\alpha) = \frac{8 \ln 2}{\pi} -
      \frac{\pi}2.
    \end{equation}
  \end{enumerate}
\item Figure~\ref{fig:mc3} shows the third central moment
  $\E[Y - \E{Y}]^3$ as a function of $\alpha$. The plot
  suggests that the third central moment is positive for each $\alpha
  > 0$, which would also establish that $Y(\alpha)$ is not normal for any
  $\alpha > 0$.  However we do not know a proof of this positive
  skewness. [Of course, the law of
  $Y(\alpha)$ is not normal for any $\alpha > 1/2$, since its support
  is a subset of $[0,\infty)$.]
  \begin{figure}[htbp]
    \centering
    \includegraphics{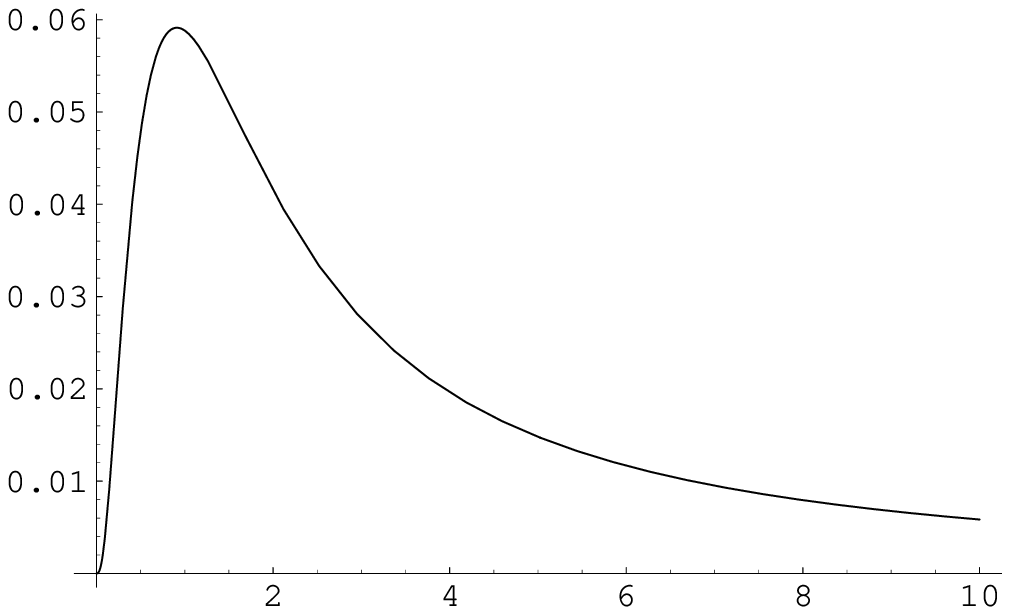}
    \caption{$\E[Y - \E{Y}]^3$ of
    Proposition~\ref{thm:limit_law_alpha_ne_12} as a function of $\alpha$.}
    \label{fig:mc3}
  \end{figure}\label{item:7.1}
\item When $\alpha=0$, the additive functional with toll sequence
  $(n^\alpha=1)_{n \geq 1}$ is $n$ for all trees with $n$ nodes.
  However, if one considers the random variable
  $\alpha^{-1/2}Y(\alpha)$ as $\alpha \downarrow 0$,
  using~\eqref{eq:7.49} and induction one can show that
  $\alpha^{-1/2}Y(\alpha)$ converges in distribution to the normal
  distribution with mean~0 and variance $4(1-\ln2)$.
\item Finally, if one considers the random variable
  $\alpha^{1/2}Y(\alpha)$ as $\alpha \to \infty$, again
  using~\eqref{eq:7.49} and induction we find that
  $\alpha^{1/2}Y(\alpha)$ converges in distribution to the unique
  distribution with $k$th moment $\sqrt{k!}$.
\end{enumerate}
\end{remark}

\begin{figure}[htbp]
  \centering
  \includegraphics[width=2.7in,keepaspectratio]{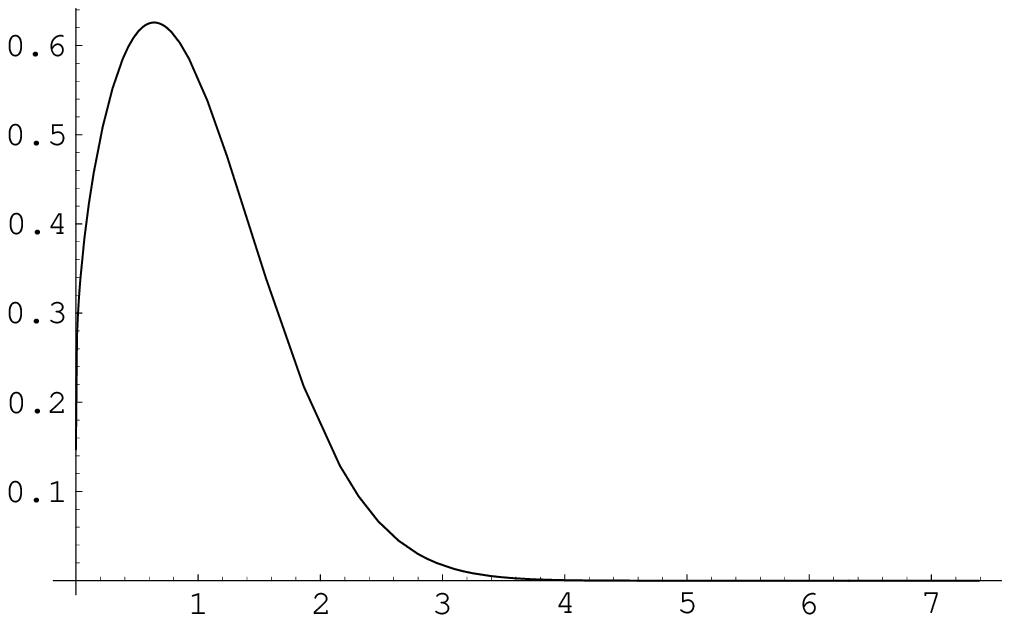}
  \includegraphics[width=2.7in,keepaspectratio]{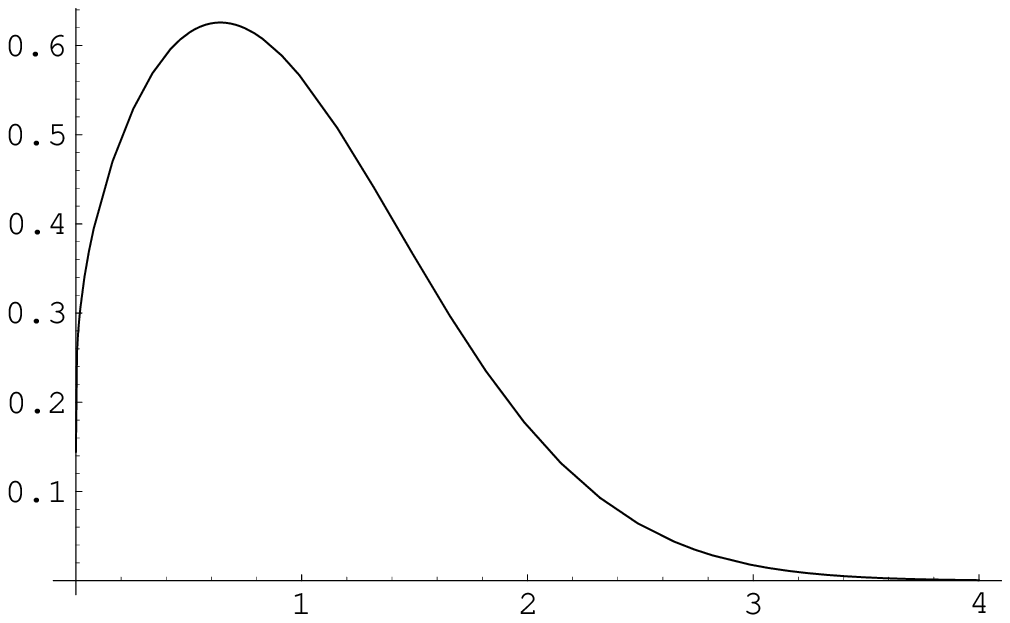}
  \caption{$f_Y$ of Remark~\ref{remark:sqrkfact}.}
  \label{fig:sqrtfactdensity}
\end{figure}

\begin{remark}
  \label{remark:sqrkfact}
  Let $Y$ be the unique distribution whose $k$th moment is
  $\sqrt{k!}$ for $k=1,2,\ldots$.  Taking $Y^*$ to be an independent
  copy of $Y$ and defining
  $X := Y Y^*$, we see immediately that $X$ is Exponential with unit
  mean.  It follows by taking logarithms that the distribution of
  $\ln Y$ is a convolution square root of the distribution of $\ln
  X$.  In particular, the characteristic function $\phi$ of $\ln Y$ has
  square equal to $\Gamma(1 + it)$ at $t \in (-\infty, \infty)$.  By
  exponential decay of $\Gamma(1 + it)$ as $t \to \pm\infty$ and
  standard theory (see, e.g.,~\cite[Chapter~XV]{MR42:5292}), $\ln Y$
  has an infinitely smooth  density on   $(-\infty, \infty)$,
  and the density and each of its derivatives are bounded.

  So $Y$ has an infinitely smooth density on $(0, \infty)$.  By change
  of variables, the density $f_Y$ of $Y$ satisfies
  \begin{equation*}
    f_Y(y) = \frac{f_{\ln Y}(\ln y)}{y}.
  \end{equation*}
  Clearly $f_Y(y)$ is bounded for $y$
  \emph{not} near~0.  (We shall drop further  consideration of
  derivatives.) To determine the behavior near~0, we need to
  know the behavior of $f_{\ln Y}(\ln y)/y$ as
  $y \to 0$.   Using the Fourier inversion formula, we may
  equivalently study
  \begin{equation*}
    e^x f_{\ln{Y}}(-x) = \frac{1}{2\pi} \int_{-\infty}^\infty
    e^{(1+it)x} \phi(t)\, dt,
  \end{equation*}
  as $x \to \infty$.  By an application of the method of steepest
  descents [(7.2.11) in~\cite{MR89d:41049}, with $g_0=1$, $\beta =
  1/2$, $w$ the identity map, $z_0=0$, and $\alpha=0$], we get
  \begin{equation*}
    f_Y(y) \sim \frac{1}{\sqrt{\pi\ln{(1/y)} }} \quad \text{as $y
      \downarrow 0$}.
  \end{equation*}
  Hence $f_Y$ is bounded everywhere.

  Using the Cauchy integral formula and simple estimates, it is easy
  to show that
  \begin{equation*}
    f_Y(y) = o( e^{-My} ) \quad \text{as $y \to \infty$}
  \end{equation*}
  for any $M < \infty$.
  Computations using the \textsc{WKB} method~\cite{MR30:3694} suggest
  \begin{equation}
    \label{eq:7.57}
    f_Y(y) \sim (2/\pi)^{1/4} y^{1/2} \exp(-y^2/2) \quad \text{as $y
      \to \infty$},
  \end{equation}
  in agreement with numerical calculations using
  \texttt{Mathematica}. [In fact, the right-side of~\eqref{eq:7.57}
  appears to be a highly accurate approximation to $f_Y(y)$ for all $y
  \geq 1$.]  Figure~\ref{fig:sqrtfactdensity} depicts the salient
  features of $f_Y$.  In particular, note the steep descent of
  $f_Y(y)$ to 0 as $y \downarrow 0$ and the quasi-Gaussian tail.
\end{remark}

\subsection{\texorpdfstring{$\alpha=1/2$}{alpha = 1/2}}
\label{sec:alpha=12}
For \(\alpha=1/2\), from~(\ref{eq:7.17}) we see immediately that
\[
\E\left[ \frac{X_n}{n\ln{n}} \right]^k = \left(
  \frac1{\sqrt\pi} \right)^k + O\left( \frac{1}{\ln n} \right).
\]
Thus the random variable $X_n/(n\ln n)$ converges in distribution to
the degenerate random variable $1/\sqrt\pi$.  To get a nondegenerate
distribution, we carry out an analysis similar to the one that led
to~(\ref{eq:7.7}), getting more precise asymptotics for the mean of
$X_n$.  The refinement of~(\ref{eq:7.7}) that we need is the following,
whose proof we omit:
\begin{equation*}
  A(z) \odot \CAT(z/4) = \frac{1}{\sqrt\pi} (1-z)^{-1/2} L(z) +
  D_0(1-z)^{-1/2} + O(|1-z|^{\tfrac12-\epsilon}),
\end{equation*}
where $D_0$ is a constant whose value we will not need in our
subsequent analysis.  By Theorem~\ref{theorem:transfer}, this leads to
\begin{equation}\label{eq:7.43}
  \E{X_n} = \frac{1}{\sqrt\pi} n \ln n + D_1 n + O(n^\epsilon),
\end{equation}
where $D_1$ is another constant whose value we will not need in our
subsequent analysis.  Now analyzing the random variable $X_n -
\pi^{-1/2} n \ln n$ in a manner similar to that of
Section~\ref{sec:small-toll-functions} we obtain
\begin{equation}\label{eq:7.44}
  \Var{[ X_n - \pi^{-1/2} n \ln n ]} =
  \left(
    \frac{8}{\pi} \ln 2 - \frac{\pi}{2}
  \right)
  n^2 + O(n^{\tfrac32 + \epsilon}).
\end{equation}
Using~(\ref{eq:7.43}) and~(\ref{eq:7.44}) we conclude that
\begin{equation*}
  \E
  \left[
    \frac{X_n - \pi^{-1/2} n \ln n - D_1 n}{n}
  \right]
  = o(1) \end{equation*}
and
\begin{equation}\label{eq:7.45}
  \Var{
  \left[
    \frac{X_n - \pi^{-1/2} n \ln n - D_1 n}{n}
  \right]}
  = \frac{8}{\pi} \ln 2 - \frac{\pi}{2} = \lim_{\alpha \to 1/2}
  \sigma^2(\alpha),
\end{equation}
where $\sigma^2 \equiv \sigma^2(\alpha)$ is defined at~(\ref{eq:7.46})
for $\alpha \ne 1/2$. [Recall~\eqref{eq:7.51} of
Remark~\ref{remark_y_alpha_properties}.]

It is possible to carry out a program similar to that of
Section~\ref{sec:catalan_n_alpha_higher-moments} to derive asymptotics
of higher order
moments using singularity analysis.  However we choose to sidestep
this arduous, albeit mechanical, computation.  Instead we will derive
the asymptotics of higher moments using a somewhat more direct
approach akin to the one employed in~\cite{MR97f:68021}.  The approach
involves approximation of sums by Riemann integrals.  To that end,
define
\begin{equation}\label{eq:7.47}
\widetilde{X}_n := X_n - \pi^{-1/2}(n+1) \ln (n+1) -
D_1(n+1),  \quad\text{ and }\quad  \hat\mu_n(k) := \frac{\beta_n}{4^{n+1}}
\E\widetilde{X}_n^k.
\end{equation}
Note that $\widetilde{X}_0 = -D_1$, $\hat\mu_n(0) = \beta_n/4^{n+1}$,
and $\hat\mu_0(k) = (-D_1)^k/4$.  Then, in a now familiar manner, for $n
\geq 1$ we find
\begin{equation*}
  \hat\mu_n(k) = 2 \sum_{j=1}^n \frac{\beta_{j-1}}{4^j}
  \hat\mu_{n-j}(k) + \hat{r}_n(k),
\end{equation*}
where now we define
\begin{multline*}
  \hat{r}_n(k) := \sum_{\substack{k_1+k_2+k_3=k\\k_1,k_2 < k}}
  \binom{k}{k_1,k_2,k_3} \sum_{j=1}^n \hat\mu_{j-1}(k_1)
  \hat{\mu}_{n-j}(k_2)\\
  \times \left[
    \frac{1}{\sqrt\pi} ( j\ln j + (n+1-j) \ln (n+1-j) - (n+1) \ln
  (n+1) + \sqrt\pi n^{1/2})
  \right]^{k_3}
\end{multline*}
Passing to generating functions and then back to sequences one gets,
for $n \geq 0$,
\begin{equation*}
  \hat{\mu}_n(k) = \sum_{j=0}^n (j+1) \frac{\beta_j}{4^j}
  \hat{r}_{n-j}(k). 
\end{equation*}
Using induction on $k$, we can approximate $\hat{r}_n(k)$ and
$\hat{\mu}_n(k)$ above by integrals and obtain the following result.
We omit the proof, leaving it as an exercise for the ambitious reader.
\begin{proposition}\label{prop_alpha_12}
  Let $X_n$ be the additive functional induced by the toll sequence
  $(n^{1/2})_{n \geq 1}$ on Catalan trees.  Define $\widetilde{X}_n$
  as in~\textrm{(\ref{eq:7.47})}.  Then
  \begin{equation*}
    \E
    [
      {\widetilde{X}_n}/{n}
    ]^k
    = m_k + o(1) \text{ as $n \to \infty$},
  \end{equation*}
  where $m_0=1$, $m_1=0$, and, for $k \geq 2$,
  \begin{equation}\label{eq:7.48}
    m_k = \frac{1}{4\sqrt\pi} \frac{\Gamma(k-1)}{\Gamma(k-\tfrac12)}
    \left[
      \sum_{\substack{k_1+k_2+k_3=k\\k_1,k_2 < k}}
    \binom{k}{k_1,k_2,k_3} m_{k_1} m_{k_2}
    \left(
      \frac{1}{\sqrt\pi}
    \right)^{k_3}
    J_{k_1,k_2,k_3}
    + 4 \sqrt\pi k m_{k-1}
    \right],
  \end{equation}
  where
  \begin{equation*}
    J_{k_1,k_2,k_3} := \int_{x=0}^1 x^{k_1-\tfrac32}
    (1-x)^{k_2-\tfrac32} [ x \ln x + (1-x) \ln (1-x) ]^{k_3} \, dx.
  \end{equation*}
  Furthermore $\widetilde{X}_n/(n+1) \stackrel{\mathcal{L}}{\to} Y$, where
  $Y$ is a random variable with the unique distribution whose
  moments are  $\E{Y^k} = m_k$, $k \geq 0$.
\end{proposition}

\subsection{A unified result}
\label{sec:unified-result}
The approach outlined in the preceding section can also be used for
the case $\alpha \ne 1/2$.  For completeness, we state the result for
that case here (without proof). 
\begin{proposition}\label{prop:alpha_ne_12_riemann}
  Let $X_n$ be the additive functional induced by the toll sequence
  $(n^{\alpha})_{n \geq 1}$ on Catalan trees.  Let $\alpha' := \alpha
  + \tfrac12$.  Define $\widetilde{X}_n$ as
  \begin{equation}
    \label{eq:7.52}
    \widetilde{X}_n :=
    \begin{cases}
      X_n - C_0(n+1) -
      \displaystyle\frac{ \Gamma(\alpha-\tfrac12
      )}{\Gamma(\alpha)}(n+1)^{\alpha'}
      & 0 < \alpha < 
      1/2, \\
      X_n - \displaystyle\frac{ \Gamma(\alpha-\tfrac12
      )}{\Gamma(\alpha)}(n+1)^{\alpha'}            & \alpha > 1/2,
    \end{cases}
  \end{equation}
  where
  \begin{equation*}
    C_0 := \sum_{n=1}^\infty n^\alpha \frac{\beta_n}{4^n}.
  \end{equation*}
  Then, for $k=0,1,2,\ldots$,
  \begin{equation*}
    \E \left[ {\widetilde{X}_n}/{n^{\alpha'}} \right]^k =
    m_k + o(1) \text{ as $n \to \infty$},
  \end{equation*}
  where $m_0=1$, $m_1=0$, and, for $k \geq 2$,
  \begin{multline}\label{eq:7.42}
    m_k = \frac{1}{4\sqrt\pi}
    \frac{\Gamma(k\alpha'-1)}{\Gamma(k\alpha'-\tfrac12)}\\
    \times
    \left[
      \sum_{\substack{k_1+k_2+k_3=k\\k_1,k_2 < k}}
      \binom{k}{k_1,k_2,k_3} m_{k_1} m_{k_2}
      \left( \frac{\Gamma(\alpha-\tfrac12)}{\Gamma(\alpha)}\right)^{k_3}
    J_{k_1,k_2,k_3}
    + 4 \sqrt\pi k m_{k-1}
    \right],
  \end{multline}
  with
  \begin{equation*}
    J_{k_1,k_2,k_3} := \int_{0}^1 x^{k_1\alpha'-\tfrac32}
    (1-x)^{k_2\alpha'-\tfrac32} [ x^{\alpha'} + (1-x)^{\alpha'} - 1
    ]^{k_3} \, dx.
  \end{equation*}
  Furthermore, $\widetilde{X}_n \stackrel{\mathcal{L}}{\to} Y$, where
  $Y$ is a random variable with the unique distribution whose moments
  are $\E{Y^k} = m_k$.
\end{proposition}

[The reader may wonder as to why we have chosen to state
Proposition~\ref{prop:alpha_ne_12_riemann} using several instances of
$n+1$, rather than $n$, in~\eqref{eq:7.52}.  The reason is that use of
$n+1$ is somewhat more natural in the calculations that establish the
proposition.]

In light of Propositions~\ref{thm:limit_law_alpha_ne_12},
\ref{prop_alpha_12}, and~\ref{prop:alpha_ne_12_riemann}, there are a variety of
ways to state a unified result.  We state one such version here.
\begin{theorem}
  \label{thm:limit-dist}
  Let \(X_n\) denote the additive functional induced by the toll sequence
  \((n^\alpha)_{n \geq 1}\) on Catalan trees. Then
  \begin{equation*}
    \frac{X_n - \E{X_n}}{\sqrt{\Var{X_n}}}
    \stackrel{\mathcal{L}}{\to} W,
  \end{equation*}
  where the distribution of~$W$ is described as follows:
  \bigskip\par\noindent
  (a) For $\alpha \ne 1/2$,
    \[
    W = \frac{1}{\sigma}
    \left(
      Y - \frac{C_1\sqrt\pi}{\Gamma(\alpha)}
    \right), \quad \text{ with } \quad  \sigma^2 :=
    \frac{C_2\sqrt\pi}{\Gamma(2\alpha+\tfrac12)} -
    \frac{C_1^2\pi}{\Gamma^2(\alpha)} > 0,
    \]
    where $Y$ is a random variable with the unique distribution whose
    moments are
    \begin{equation*}
      \E{Y^k} =
      \frac{C_k\sqrt\pi}{\Gamma(k(\alpha+\tfrac12)-\tfrac12)},
    \end{equation*}
    and the \(C_k\)'s satisfy the recurrence~\textrm{(\ref{eq:7.10})}.
    \bigskip\par\noindent
    (b) For $\alpha=1/2$,
    \begin{equation*}
      W = \frac{Y}{\sigma}, \quad \text{ with } \quad \sigma^2 :=
      \frac{8}{\pi} \ln 2 - \frac{\pi}{2},
    \end{equation*}
    where $Y$ is a random variable with the unique distribution whose
    moments $m_k = \E{Y^k}$ are given by~\textrm{(\ref{eq:7.48})}.
\end{theorem}
\begin{proof}
  Define
  \begin{equation*}
    W_n := \frac{X_n - \E{X_n}}{\sqrt{\Var{X_n}}}
  \end{equation*}
  \bigskip\par\noindent
  (a) Consider first the case $\alpha < 1/2$ and let $\alpha' := \alpha +
    \tfrac12$.  By~\eqref{eq:7.39},
    \begin{equation}\label{eq:7.53}
      \E{X_n} = C_0(n+1) + \frac{C_1\sqrt\pi}{\Gamma(\alpha)}
      n^{\alpha'} + o(n^{\alpha'}).
    \end{equation}
    Since $\widetilde{X}_n$ defined at~\eqref{eq:7.52} and $X_n$ differ
    by a deterministic amount, $\Var{X_n}=\Var{\widetilde{X}_n}$. Now by
    Proposition~\ref{prop:alpha_ne_12_riemann},
    \begin{equation}\label{eq:7.54}
      \Var{\widetilde{X}_n} = \E{\widetilde{X}_n^2} -
      (\E{\widetilde{X}_n})^2 = (m_2 + o(1))n^{2\alpha'} - (m_1^2 +
      o(1))n^{2\alpha'} = (m_2 + o(1))n^{2\alpha'}.
    \end{equation}
    So $\sigma^2$ equals $m_2$ 
    defined at~(\ref{eq:7.42}), namely,
    \begin{equation*}
      \frac{1}{4\sqrt\pi}\frac{\Gamma(2\alpha'-1)}{\Gamma(2\alpha'-\tfrac12)}
      \left(\frac{\Gamma(\alpha-\tfrac12)}{\Gamma(\alpha)}\right)^2 J_{0,0,2}.
    \end{equation*}
    Thus to show $\sigma^2 > 0$ it is enough to show that $J_{0,0,2} > 0$.  But
    \begin{equation*}
      J_{0,0,2} = \int_{x}^1 x^{-3/2} (1-x)^{-3/2} [ x^{\alpha'} +
      (1-x)^{\alpha'} - 1]^2 \,dx,
    \end{equation*}
    which is clearly positive.  Using~\eqref{eq:7.53} and~\eqref{eq:7.54},
    \begin{equation*}
      W_n = \frac{X_n - C_0(n+1) - \frac{C_1\sqrt\pi}{\Gamma(\alpha)}
      n^{\alpha'} + o(n^{\alpha'})}{ (1 + o(1))\sigma
      n^{\alpha'}},
    \end{equation*}
    so, by Proposition~\ref{thm:limit_law_alpha_ne_12} and Slutsky's
    theorem~\cite[Theorem~25.4]{MR95k:60001}, the claim
    follows.
    
    The case $\alpha > 1/2$ follows similarly.
    \bigskip\par\noindent
  (b) When $\alpha=1/2$, 
  \begin{equation*}
    \E{X_n} = \frac{1}{\sqrt\pi} n\ln{n} + D_1 n + o(n)
  \end{equation*}
  by~\eqref{eq:7.43} and
  \begin{equation*}
    \Var{X_n} =
    \left(
      \frac{8}{\pi}\ln2 - \frac{\pi}{2} + o(1)
    \right)n^2
  \end{equation*}
  by~\eqref{eq:7.45}. The claim then follows easily from
    Proposition~\ref{prop_alpha_12} and Slutsky's
    theorem.
\end{proof}


\chapter{Catalan model: the shape functional}
\label{sec:catalan_shape-functional}

We now turn our attention to the shape functional for Catalan trees.
The shape functional is the cost induced by the toll function
$t_n \equiv \log{n}$, $n \geq 1$.  For background and results on the shape
functional, we refer the reader to~\cite{MR97f:68021} and~\cite{MR99j:05171}.

In the sequel we will improve on the mean and variance estimates
obtained in~\cite{MR97f:68021} and derive a central limit theorem for
the shape functional for Catalan trees.  The technique employed is
singularity analysis followed by the method of moments.

\section{Mean}
\label{sec:shape_function_mean}

We use the notation and techniques of
Section~\ref{sec:cat_n_alpha_asymptotics-mean} again.  Observe that
now $T(z) = \Li_{0,1}(z)$ and~\eqref{eq:7.21} gives
\begin{equation*}
  \CAT(z/4) = 2 - \frac2{\Gamma(-1/2)}[ \Li_{3/2,0}(z) -
  \zeta(3/2)] + 2\left(1-\frac{\zeta(1/2)}{\Gamma(-1/2)} \right)(1-z) +
  O(|1-z|^{3/2}).
\end{equation*}
So
\begin{equation*}
  T(z)\odot\CAT(z/4) = -\frac2{\Gamma(-1/2)} \Li_{3/2,1}(z)
  + \bar{c} + \bar{\bar{c}}(1-z) +  O(|1-z|^{\tfrac32-\epsilon}),
\end{equation*}
where \( \bar{c} \) and \( \bar{\bar{c}} \) denote unspecified (possibly 0)
constants.  The constant term in the singular expansion of
\( T(z)\odot\CAT(z/4) \) is already known to be, using
Theorem~\ref{thm:hadamard},
\begin{equation*}
  C_0 = T(z)\odot\CAT(z/4) \Bigr\rvert_{z=1} = \sum_{n=1}^\infty
  (\log{n}) \frac{\beta_n}{4^n}.
\end{equation*}
Now using Theorem~\ref{theorem:Li-sing} to obtain the singular expansion
of \( \Li_{3/2,1}(z) \), we get
\begin{equation*}
  T(z) \odot \CAT(z/4) =  C_0 - 2(1-z)^{1/2}L(z) -
  2(2(1-\log(2))-\gamma)(1-z)^{1/2} + O(|1-z|),
\end{equation*}
so that
\begin{equation}
  \label{eq:8.22}
  A(z)\odot\CAT(z/4) = C_0(1-z)^{-1/2} - 2L(z) -
  2(2(1-\log2)-\gamma) + O(|1-z|^{1/2}).
\end{equation}
Using Theorem~\ref{theorem:transfer} and the asymptotics of the Catalan numbers
we get that the mean \( a_n \) of the shape functional  is given by
\begin{equation}
  \label{eq:8.24}
  a_n = C_0(n+1) - 2\sqrt\pi n^{1/2} + O(1),
\end{equation}
which agrees with the estimate in Theorem~3.1 of~\cite{MR97f:68021}
and improves the remainder estimate.

\section{Second moment and variance}
\label{sec:shape_functional_variance}

We now derive the asymptotics of the approximately centered second
moment and the variance of the shape functional.  These estimates will
serve as the basis for the induction to follow.  We will use the
notation of Section~\ref{sec:small-toll-functions}, centering the cost
function as before by $C_0(n+1)$. 

It is clear from~\eqref{eq:8.22} that
\begin{equation}\label{eq:8.38}
  \widehat{M}_1(z) = -2L(z) - 2(2(1-\log2)-\gamma) +
  O(|1-z|^{1/2}),
\end{equation}
and~\eqref{eq:7.3} with \( k=2 \) gives us, recalling~\eqref{eq:7.50},
\begin{equation}
  \label{eq:8.34}
  \widehat{R}_2(z) = C_0^2 + \CAT(z/4)\odot\Li_{0,2}(z) +  4\Li_{0,1}(z)
  \odot [\frac{z}4 \CAT(z/4)
  \widehat{M}_1(z)] + \frac{z}2\widehat{M}_1^2(z).
\end{equation}
We analyze each of the terms in this sum.  For the last term, observe
that $z/2 \to 1/2$ as $z \to 1$, so that
\begin{equation*}
  \frac{z}2 \widehat{M}_1^2(z) = 2L^2(z) + 4(2(1-\log2)-\gamma)L(z) +
  2(2(1-\log2)-\gamma)^2 + O(|1-z|^{\tfrac12-\epsilon}),
\end{equation*}
the \( \epsilon \) introduced to avoid logarithmic remainders.  The first
term is easily seen to be
\begin{equation*}
 \CAT(z/4) \odot \Li_{0,2}(z) = K + O(|1-z|^{\tfrac12-\epsilon})
\end{equation*}
by Theorem~\ref{thm:hadamard}, where
\begin{equation*}
  K := \sum_{n=1}^\infty (\log{n})^2 \frac{\beta_n}{4^n}.
\end{equation*}
For the middle term, first observe that
\begin{equation*}
  \frac{z}4 \CAT(z/4)\widehat{M}_1(z) = -L(z) - (2(1-\log2)-\gamma) +
  (1-z)^{1/2}L(z) + O(|1-z|^{1/2})
\end{equation*}
and that \( L(z) = \Li_{1,0}(z) \). Thus the third term on the right
in~\eqref{eq:8.34} is 4 times:
\begin{equation*}
  -\Li_{1,1}(z) + \bar{c} + O(|1-z|^{\tfrac12-2\epsilon}) = -\frac12
   L^2(z) + \gamma L(z) + \bar{c} + O(|1-z|^{\tfrac12-\epsilon}).
\end{equation*}
[The singular expansion for \( \Li_{1,1}(z) \)
was obtained using~\eqref{eq:4.1.posintalpha}.  We state it here for
the reader's convenience:  
\begin{equation*}
  \Li_{1,1}(z) = \frac12 L^2(z) - \gamma L(z) + \bar{c} + O(|1-z|),
\end{equation*}
where $\bar{c}$ is again an unspecified constant.]  Hence
\begin{equation*}
  \widehat{R}_2(z) = 8(1-\log2)L(z) + \bar{c} +
  O(|1-z|^{\tfrac12-\epsilon}), 
\end{equation*}
which leads to
\begin{equation}
  \label{eq:8.23}
  \widehat{M}_2(z) = 8(1-\log2)(1-z)^{-1/2}L(z) +
  \bar{c}(1-z)^{-1/2} + O(|1-z|^{-\epsilon}).
\end{equation}
We draw the attention of the reader to the cancellation
of the ostensible lead-order term \( L^2(z) \).  This kind of
cancellation will appear again in the next section when we deal with
higher moments.

Now using singularity analysis and estimates for the Catalan numbers
we get
\begin{equation}
  \label{eq:8.25}
  \tilde\mu_n(2) = 8(1-\log2)n\log{n} + \bar{c}n + O(n^{\tfrac12+\epsilon}).
\end{equation}
Using~\eqref{eq:8.24},
\begin{equation*}
  \Var{X_n} = \tilde\mu_n(2) - \tilde\mu_n(1)^2 =
  8(1-\log2)n\log{n} + \bar{c}n + O(n^{\tfrac12+\epsilon}),
\end{equation*}
which agrees with Theorem~3.1 of~\cite{MR97f:68021} (after a
correction pointed out in~\cite{MR99j:05171}) and improves the
remainder estimate.  In our subsequent analysis we will not need to
evaluate the unspecified constant $\bar{c}$.

\section{Higher moments}
\label{sec:shape_function_higher-moments}

We now turn our attention to deriving the asymptotics of higher
moments of the shape functional.  The main result is as follows.
\begin{proposition}
  \label{thm:shape_moments}
  Define \(\widetilde{X}_n := X_n - C_0(n+1)\), with \(X_0 := 0\);
  \(\tilde{\mu}_n(k) := \E{\widetilde{X}_n^k} \), with
  $\tilde{\mu}_n(0) = 1$ for
  all $n \geq 0$; and \(\hat{\mu}_n(k)
  := \beta_n\tilde{\mu}_n(k)/4^n \).  Let \(\widehat{M}_k(z)\) denote
  the ordinary 
  generating function of \(\hat{\mu}_n(k)\) in the argument $n$.  For
  \( k \geq 2 \), \( \widehat{M}_k(z) \) satisfies
  \begin{equation*}
    \widehat{M}_k(z) = (1-z)^{-\tfrac{k-1}2}
    \sum_{j=0}^{\lfloor {k}/2 \rfloor} C_{k,j}
    L^{\lfloor {k}/2 \rfloor-j}(z) +
    O(|1-z|^{-\tfrac{k}2+1-\epsilon}), 
  \end{equation*}
  with
  \begin{equation*}
    C_{2l,0} = \frac14 \sum_{j=1}^{l-1} \binom{2l}{2j}
    C_{2j,0}C_{2l-2j,0}, \qquad C_{2,0} = 8(1-\log2).
  \end{equation*}
\end{proposition}
\begin{proof}
  The proof is by induction.  For \( k=2 \) the claim is true
  by~\eqref{eq:8.23}.  We note that the claim is \emph{not} true for
  \( k=1 \). Instead, recalling~(\ref{eq:8.38}),
  \begin{equation}
    \label{eq:8.28}
    \widehat{M}_1(z) = -2L(z) -
    2(2(1-\log2)-\gamma) + O(|1-z|^{1/2}).
  \end{equation}
  For the induction step, let \( k \geq 3 \).  We will first get the
  asymptotics of \( \widehat{R}_k(z) \) defined at~\eqref{eq:7.3} with
  $T(z) = \Li_{0,1}(z)$.  In order to do that we will obtain the
  asymptotics of each term in the defining sum.  We remind the reader
  that we are only interested in the form of the asymptotic expansion
  of \( \widehat{R}_k(z) \) and the coefficient of the lead order term
  when \( k \) is even.  This allows us to ``define away'' all other
  constants, their determination delayed to the time when the need
  arises.

  For this paragraph suppose that \( k_1 \geq2 \) and \( k_2 \geq 2 \).  Then
  by the induction hypothesis
  \begin{multline}
    \label{eq:8.26}
    \frac{z}4 \widehat{M}_{k_1}(z)\widehat{M}_{k_2}(z) = \frac14
    (1-z)^{-\tfrac{k_1+k_2}2+1}
    \sum_{l=0}^{\lfloor {k_1}/2 \rfloor +
    \lfloor {k_2}/2 \rfloor} A_{k_1,k_2,l}
    L^{\lfloor {k_1}/2 \rfloor +  \lfloor {k_2}/2
      \rfloor - l}(z)\\
    {}+ O(|1-z|^{-\tfrac{k_1+k_2}2+\tfrac32-\epsilon}),
  \end{multline}
  where \( A_{k_1,k_2,0}= C_{k_1,0}C_{k_2,0} \).  (a) If \( k_3=0 \)
  then \( k_1+k_2=k \) and the corresponding contribution to \(
  \widehat{R}_k(z) \) is given by
  \begin{equation}
    \label{eq:8.27}
    \frac14 \binom{k}{k_1} (1-z)^{-\tfrac{k}2+1}
    \sum_{l=0}^{\lfloor {k_1}/2 \rfloor +
      \lfloor ({k-k_1})/2 \rfloor} A_{k_1,k-k_1,l}
    L^{\lfloor {k_1}/2 \rfloor +  \lfloor ({k-k_1})/2 \rfloor -
      l}(z) + O(|1-z|^{-\tfrac{k}2+\tfrac32-\epsilon}).
  \end{equation}
  Observe that if \( k \) is even and \( k_1 \) is odd the highest power of
  \( L(z) \) in~\eqref{eq:8.27} is \( \lfloor {k}/2 \rfloor-1 \).  In all
  other cases the the highest power of \( L(z) \) in~\eqref{eq:8.27} is
  \( \lfloor {k}/2 \rfloor \).  (b) If \( k_3 \ne 0 \) then we  use
  Lemma~\ref{lem:omztoli} to express~\eqref{eq:8.26} as a linear
  combination of
  \begin{equation*}
    \left\{
      \Li_{-\tfrac{k_1+k_2}2+2,l}(z) \right\}_{l=0}^{\lfloor {k_1}/2 \rfloor
      + \lfloor {k_2}/2 \rfloor}
  \end{equation*}
  with a remainder that is
  $O(|1-z|^{-\tfrac{k_1+k_2}2+\tfrac32-\epsilon})$.  When we take the 
  Hadamard product of such a term with \( \Li_{0,k_3}(z) \) we will get a
  linear combination of
    \begin{equation*}
    \left\{
      \Li_{-\tfrac{k_1+k_2}2+2,l+k_3}(z) \right\}_{l=0}^{\lfloor
      {k_1}/2 \rfloor + \lfloor {k_2}/2 \rfloor}
  \end{equation*}
  and a smaller remainder.  Such terms are all
  \( O(|1-z|^{-\tfrac{k_1+k_2}2+1-\epsilon}) \), which is
  \( O(|1-z|^{-\tfrac{k}2+\tfrac32-\epsilon}) \).

  Next, consider the case when \( k_1=1 \) and \( k_2 \geq 2 \).  Using the
  induction hypothesis and~\eqref{eq:8.28} we get
  \begin{equation}
    \label{eq:8.29}
    \frac{z}4 \widehat{M}_{k_1}(z) \widehat{M}_{k_2}(z) = -\frac12
    (1-z)^{-\tfrac{k_2-1}2} \sum_{j=0}^{\lfloor {k_2}/2 \rfloor+1}
    B_{k_2,j} L^{\left\lfloor \tfrac{k_2}2 \right\rfloor + 1 - j}(z) +
    O(|1-z|^{-\tfrac{k_2}2+1-2\epsilon}),
  \end{equation}
  with \( B_{k_2,0} = C_{k_2,0} \).  (a) If \( k_3=0 \) then \(
  k_2=k-1 \) and the 
  corresponding contribution to \( \widehat{R}_k(z) \) is given by
  \begin{equation}
    \label{eq:8.30}
    -\frac{k}{2} (1-z)^{-\tfrac{k}2+1}
     \sum_{j=0}^{\lfloor ({k-1})/{2} \rfloor+1} B_{k-1,j}
     L^{\left\lfloor \tfrac{k-1}2 \right\rfloor+1-j}(z) +
     O(|1-z|^{-\tfrac{k}2+\tfrac32-2\epsilon}).
   \end{equation}
   (b) If \( k_3 \ne 0 \) then Lemma~\ref{lem:omztoli} can be used once again
   to express~\eqref{eq:8.29} in terms of generalized polylogarithms,
   whence an argument similar to that at the end of the preceding paragraph
   yields that the contributions from such terms is
   \( O(|1-z|^{-\tfrac{k_2-1}2-\epsilon}) \), which is
   \( O(|1-z|^{-\tfrac{k}2+\tfrac32-\epsilon}) \).   The case when \(k_1
   \geq 2\)  and  \(k_2 =1\) is handled symmetrically.

   When \( k_1=k_2=1 \) then \(
   (z/4)\widehat{M}_{k_1}(z)\widehat{M}_{k_2}(z) \) is 
   \( O(|1-z|^{-\epsilon}) \) and when one takes the Hadamard product of
   this term with \( \Li_{0,k_3}(z) \) the contribution will be
   \( O(|1-z|^{-2\epsilon}) \).

   Now consider the case when \( k_1=0 \) and \( k_2 \geq 2 \).  Now
   \( \widehat{M}_0(z) = \CAT(z/4) \), so that
   \begin{equation}
     \label{eq:8.31}
     \frac{z}4 \widehat{M}_{k_1}(z)\widehat{M}_{k_2}(z) = \frac12
     (1-z)^{-\tfrac{k_2-1}2} \sum_{j=0}^{\lfloor {k_2}/2 \rfloor}
     C_{k_2,j} L^{\lfloor {k_2}/2 \rfloor-j}(z) +
     O(|1-z|^{-\tfrac{k_2}2+1-\epsilon}).
   \end{equation}
   By Lemma~\ref{lem:omztoli} this can be expressed as a linear
   combination of
   \begin{equation*}
     \left\{ \Li_{-\tfrac{k_2-1}2+1,j}(z)
     \right\}_{j=0}^{\lfloor {k_2}/2 \rfloor}
   \end{equation*}
   with a \( O(|1-z|^{-\tfrac{k_2}2+1-\epsilon}) \) remainder.  When we
   take the Hadamard product of such a term with \( \Li_{0,k_3}(z) \) we
   will get a linear combination, call it $S(z)$, of
   \begin{equation*}
     \left\{ \Li_{-\tfrac{k_2-1}2+1,j+k_3}(z)
     \right\}_{j=0}^{\lfloor {k_2}/2 \rfloor}
   \end{equation*}
   with a remainder of \( O(|1-z|^{-\tfrac{k_2}2 + 1 - 2\epsilon})
   \), which is  \( O(|1-z|^{-\tfrac{k}2 + \tfrac32 - 2\epsilon}) \)
   unless \( k_2=k-1 \).  When 
   \( k_2 =  k-1 \), by Lemma~\ref{lem:omztoli} the constant multiplying
   the lead-order term 
   \( \Li_{-\tfrac{k}{2}+2,\lfloor\tfrac{k-1}{2}\rfloor}(z) \) in $S(z)$ is
   $\frac{C_{k-1,0}}2
   \mu_0^{(-\tfrac{k}{2}+2,\lfloor\tfrac{k-1}{2}\rfloor)}$.  When we
   take the Hadamard product of this term with \( \Li_{0,k_3}(z) \) we get
   a lead-order term of
   \begin{equation*}
   \frac{C_{k-1,0}}2
   \mu_0^{(-\tfrac{k}{2}+2,\lfloor\tfrac{k-1}{2}\rfloor)}
   \Li_{-\tfrac{k}{2}+2,\lfloor\tfrac{k-1}{2}\rfloor+1}(z).
   \end{equation*}
   Now we use Lemma~\ref{lem:litoomz} and the observation that
   \( \lambda_0^{(\alpha,r)}\mu_0^{(\alpha,s)}=1 \) to conclude that the
   contribution to $\widehat{R}_k(z)$ from the term with \( k_1=0 \)
   and \( k_2=k-1 \) is 
   \begin{equation}
     \label{eq:8.32}
     \frac{k}2 (1-z)^{-\tfrac{k}2+1}
     \sum_{j=0}^{\lfloor\frac{k-1}{2}\rfloor+1} D_{k,j}
     L^{\lfloor\frac{k-1}{2}\rfloor+1-j}(z)
     + O(|1-z|^{-\tfrac{k}2+\tfrac{3}2-\epsilon}),
   \end{equation}
   with \( D_{k,0}=C_{k-1,0} \).  Notice that the lead order from this
   contribution is precisely that from~\eqref{eq:8.30} but with opposite
   sign; thus the two contributions cancel each other to lead order.
   The case $k_2 = 0$ and $k_1 \geq 2$ is handled symmetrically.
   
   The last two cases are \( k_1=0 \), \( k_2=1 \) (or vice-versa)
   and \( k_1=k_2=0 \).  The contribution from these cases can be
   easily seen to be \( O(|1-z|^{-\tfrac{k}2+\tfrac32-2\epsilon}) \).
   
   We can now deduce the asymptotic behavior of \( \widehat{R}_k(z)
   \).  The three contributions are~\eqref{eq:8.27}, \eqref{eq:8.30},
   and~\eqref{eq:8.32}, with only~\eqref{eq:8.27} contributing in net
   a term of the form \( (1-z)^{-\tfrac{k}{2}+1} L^{\lfloor {k}/{2}
     \rfloor}(z) \) when \(k\) is even.  The coefficient of this term
   when \( k \) is even is given by
   \begin{equation*}
     \frac14 \sum_{\substack{0 < k_1 < k\\k_1\text{ even}}}
     \binom{k}{k_1} C_{k_1,0}C_{k_2,0}.
   \end{equation*}
   Finally we can sum up the rest of the contribution, define
   \( C_{k,j} \) appropriately and use~\eqref{eq:7.4} to claim the result.
\end{proof}

\section{A central limit theorem}
\label{sec:shape_functional_centr-limit-theor}

Proposition~\ref{thm:shape_moments} and singularity analysis
(Theorem~\ref{theorem:transfer}) allows us to
get the asymptotics of the moments of the ``approximately centered'' shape
functional.  Using arguments identical to those in
Section~\ref{sec:asymptotics-moments} it is clear that for \( k \geq 2 \)
\begin{equation*}
  \tilde{\mu}_n(k) = \frac{C_{k,0}\sqrt\pi}{\Gamma(\tfrac{k-1}2)}
  n^{{k}/{2}} [\log{n}]^{\lfloor{k}/{2}\rfloor} +
  O(n^{{k}/{2}} [\log{n}]^{\lfloor{k}/{2}\rfloor-1}).
\end{equation*}
This and the asymptotics of the mean derived in
Section~\ref{sec:shape_function_mean} give us, for \( k \geq 1 \),
\begin{equation*}
  E\left[ \frac{\tilde{X}_n}{\sqrt{n\log{n}}} \right]^{2k} \to
  \frac{C_{2k,0}\sqrt\pi}{\Gamma(k-\tfrac12)}, \qquad E\left[
  \frac{\tilde{X}_n}{\sqrt{n\log{n}}} \right]^{2k-1}  = o(1)
\end{equation*}
as $n \to \infty$.
The recurrence for \( C_{2k,0} \) can be solved easily to yield, for $k
\geq 1$,
\begin{equation*}
  C_{2k,0} = \frac{(2k)!(2k-2)!}{2^k2^{2k-2}k!(k-1)!} \sigma^{2k},
\end{equation*}
where \( \sigma^2 := 8(1-\log2) \). Then using the identity
\begin{equation*}
  \frac{\Gamma(k-\tfrac12)}{\sqrt\pi} = \rising{\left(\frac12\right)}{k-1} =
  \left[2^{2k-2} \frac{(k-1)!}{(2k-2)!}\right]^{-1}
\end{equation*}
we get
\begin{equation*}
  \frac{C_{2k,0}\sqrt\pi}{\Gamma(k-\tfrac12)} = \frac{(2k)!}{2^k
  k!}\sigma^{2k}.
\end{equation*}
It is clear now that both the ``approximately centered'' and the
normalized shape functional are asymptotically normal.
\begin{theorem}
  \label{thm:shape_clt}
  Let \( X_n \) denote the shape functional, induced by the toll
  sequence $(\log{n})_{n \geq 1}$, for Catalan trees.  Then
  \begin{equation*}
    \frac{X_n-C_0(n+1)}{\sqrt{n\log{n}}} \stackrel{\mathcal{L}}{\to}
    \mathcal{N}(0,\sigma^2) \quad\text{ and }\quad
    \frac{X_n - \E{X_n} }{\sqrt{ \Var{X_n} }}
    \stackrel{\mathcal{L}}{\to}  \mathcal{N}(0,1),
  \end{equation*}
  where
  \begin{equation*}
  C_0 := \sum_{n=1}^\infty (\log{n}) \frac{\beta_n}{4^n}, \qquad
  \beta_n = \frac1{n+1}\binom{2n}{n},
  \end{equation*}
  and \( \sigma^2 := 8(1-\log2) \).
\end{theorem}

\section{Sufficient conditions for asymptotic normality}
\label{sec:suff-cond-asympt}

In this speculative final section we briefly examine the behavior of a
general additive functional $X_n$ induced by a given ``small'' toll
sequence $(t_n)$.  We have seen evidence
[Remark~\ref{remark_y_alpha_properties}(\ref{item:7.1})] that if
$(t_n)$ is the ``large'' toll sequence $n^{\alpha}$ for any fixed $\alpha
>  0$, then the limiting behavior is non-normal.  When $t_n = \log n$ (or
$t_n = n^\alpha$ and $\alpha \downarrow 0$), the (limiting) random
variable is normal.  Where is the interface between normal and non-normal
asymptotics?  We have carried out arguments similar to those leading to
Propositions~\ref{prop_alpha_12} and~\ref{prop:alpha_ne_12_riemann} (see
also~\cite{MR97f:68021}) that suggest a sufficient condition for asymptotic
normality, but our ``proof'' is somewhat heuristic, and further technical
conditions on $(t_n)$ may be required.  Nevertheless, to inspire further
work, we present our preliminary indications.

We assume that $t_n \equiv t(n)$, where $t(\cdot)$ is a function of a
nonnegative real argument.  Suppose that $x^{-3/2} t(x)$ is (ultimately)
nonincreasing and that $x t'(x)$ is slowly varying at infinity.  Then
\begin{equation*}
   \E{X_n} = C_0 (n+1) - (1 + o(1)) 2 \sqrt{\pi} n^{3/2} t'(n),
\end{equation*}
where
\begin{equation*}
   C_0 = \sum_{n=1}^\infty t_n \frac{\beta_n}{4^n}.
\end{equation*}
Furthermore,
\begin{equation*}
   \Var{X_n} \sim 8 (1 - \log 2) [n t'(n)]^2 n\log n,
\end{equation*}
and
\begin{equation*}
   \frac{X_n - C_0 (n + 1)}{n t'(n) \sqrt{n \log n}}
   \stackrel{\mathcal{L}}{\to} \mathcal{N}(0,\sigma^2), \text{ where }
   \sigma^2 = 8(1 - \log 2).
\end{equation*}
This asymptotic normality can also be stated in the form
\begin{equation*}
   \frac{X_n - \E{X_n} }{\sqrt{\Var{X_n}}}
   \stackrel{\mathcal{L}}{\to} \mathcal{N}(0,1).
\end{equation*}



\part{Appendices and Bibliography}

\appendix
\counterwithin{equation}{chapter}
\counterwithin{figure}{chapter}
\counterwithin{table}{chapter}
\counterwithin{theorem}{chapter}

\chapter{Solution of an Euler differential equation}
\label{appendix:diffeq}

We now provide the proof of Theorem~\ref{thm:diffeq-soln}, which states
the general solution of the differential equation~\eqref{eq:2.8} with
initial conditions \( a_j = b_j \), \( 0 \leq j \leq m-2 \).  This
linear differential equation can be written in the form
\begin{equation}
  \label{eq:A.1}
  \mathbf{L}g = h
\end{equation}
where the operator \( \mathbf{L} \) is defined as
\begin{equation}
  \label{eq:A.1.5}
  (\mathbf{L}g)(z) := g^{(m-1)}(z) - m! (1-z)^{-(m-1)} g(z).
\end{equation}
We seek the solution $g=A$ corresponding to input
$h = B^{(m-1)}$.

Equations of the form~\eqref{eq:A.1}--\eqref{eq:A.1.5} are members of a
class known as \emph{Euler differential equations}. In this appendix
we discuss a general method for solving Euler equations, restricting
attention, for the sake of definiteness and practicality,
to~\eqref{eq:A.1}--\eqref{eq:A.1.5}.  We assume that the reader is
familiar with the theory of linear differential equations with
constant coefficients (see,
e.g.,~\cite{boyce86:_elemen_differ_equat}).

\section*{The homogeneous solution}
\label{sec:homogeneous-solution}

The technique for solving \( \mathbf{L}g = 0 \) is quite easily
summarized: make the change of variable \( z = 1 - e^{-x} \), that is,
\( x = \ln{((1-z)^{-1}}) \).  For notational convenience we will
abbreviate $\ln{((1-z)^{-1}})$ as $L(z)$.  Then consider the function
\( \tilde{g} \) defined by
\begin{equation}
  \label{eq:A.2}
  \tilde{g}(x) := g(1-e^{-x}), \quad \text{ i.e., } \quad g(z) =
  \tilde{g}(L(z)).
\end{equation}
\begin{lemma}
  \label{lemma:A.1}
  The derivatives of \( g \) are related to those of \( \tilde{g} \)
  by
  \begin{equation}
    \label{eq:A.3}
    {g}^{(k)}(z) = (1-z)^{-k} \sum_{j=0}^k \stirlingone{k}{j}
    \tilde{g}^{(j)}(L(z)),
  \end{equation}
  where \( \stirlingone{k}{j} \) denotes a signless Stirling number
  of the first kind.
\end{lemma}
\begin{proof}
  Since \( \stirlingone{0}{0}=1 \), the basis of the induction
  is~\eqref{eq:A.2}.  For the induction step we differentiate the
  induction hypothesis~\eqref{eq:A.3} to get
  \begin{align*}
    g^{(k+1)}(z) &= (1-z)^{-(k+1)} \left[ k \sum_{j=0}^k
      \stirlingone{k}{j} \tilde{g}^{(j)}(L(z)) +
      \sum_{j=0}^k \stirlingone{k}{j}
      \tilde{g}^{(j+1)}(L(z)) \right]\\
    &= (1-z)^{-(k+1)} \Biggl\{ k \stirlingone{k}{0}
      \tilde{g}^{(0)}(L(z)) \\
    & \quad + \sum_{j=1}^k \left(
        k\stirlingone{k}{j} + \stirlingone{k}{j-1} \right)
      \tilde{g}^{(j)}(L(z)) + \stirlingone{k}{k}
      \tilde{g}^{(k+1)}(L(z)) \Biggr\}.
  \end{align*}
  We now use standard identities for the Stirling numbers~\citep[\S
  1.2.6]{knuth97}:
  \begin{equation*}
    k \stirlingone{k}{0} = 0 = \stirlingone{k+1}{0},\; k \geq 0,
    \qquad
    \stirlingone{k}{k} = 1 = \stirlingone{k+1}{k+1},\; k \geq 0,
  \end{equation*}
  and
  \begin{equation*}
  k \stirlingone{k}{j} +
    \stirlingone{k}{j-1} = \stirlingone{k+1}{j},\; 1 \leq j \leq k.
  \end{equation*}
  Thus
  \begin{equation*}
    g^{(k+1)}(z) = (1-z)^{-(k+1)} \sum_{j=0}^{k+1}
    \stirlingone{k+1}{j} \tilde{g}^{(j)}(L(z)),
  \end{equation*}
  completing the induction.
\end{proof}
The left-hand side of~\eqref{eq:A.1} can hence be expressed as
\begin{equation*}
  (\mathbf{L}g)(z) = (1-z)^{-(m-1)} \left\{ \sum_{j=0}^{m-1}
  \stirlingone{m-1}{j} \tilde{g}^{(j)}(L(z)) - m!
  \tilde{g}(L(z)) \right\}
\end{equation*}
so that solving \( \mathbf{L}g = 0 \) is equivalent to solving \(
\mathbf{\widetilde{L}}\tilde{g} = 0 \), where
\begin{equation}
  \label{eq:A.4}
  \mathbf{\widetilde{L}}\tilde{g}(x) := \sum_{j=0}^{m-1} \stirlingone{m-1}{j}
  \tilde{g}^{(j)}(x) - m! \tilde{g}(x).
\end{equation}
But this is a linear differential equation with constant coefficients.
Its \emph{indicial polynomial}, or \emph{characteristic polynomial},
is
\begin{equation}
  \label{eq:A.5}
  \psi_m(\lambda) \equiv \psi(\lambda) := \sum_{j=0}^{m-1}
  \stirlingone{m-1}{j} \lambda^j - m! = \rising{\lambda}{m-1} - m!,
\end{equation}
the last equality following from~\cite[1.2.6-(44)]{knuth97}.  This
polynomial is discussed in some detail
in~\cite[Chapter~3]{MR93f:68045} and~\cite{MR96j:68042}; some useful
identities involving it are stated in
Appendix~\ref{appendix:indicial}.  From the discussion
in~\cite{MR93f:68045} it follows that there are \( m - 1\) distinct
(and nonzero)
roots of \( \psi \), call them \( \lambda_1,\ldots, \lambda_{m-1} \)
arranged in nonincreasing order of real parts.  Thus the functions \(
\exp{( \lambda_j x )} \) are \( m -1 \) linearly independent solutions
of~\eqref{eq:A.4} and hence the functions \( (1-z)^{-\lambda_j} \)
form a basis of solutions to \( \mathbf{L}g =0 \).

\section*{A particular solution}
\label{sec:particular-solution}

We shall find it convenient to continue working in the \( x \)-domain
of~\eqref{eq:A.4}.  There the differential equation equivalent to \(
\mathbf{L}g = h \) is \( \mathbf{\widetilde{L}}\tilde{g} = \tilde{h} \) where
\begin{equation}
  \label{eq:A.5.1}
  \tilde{h}(x) := e^{-(m-1)x} h(1-e^{-x}), \quad \text{ i.e., } \quad
  h(z) = (1-z)^{-(m-1)} \tilde{h}(L(z)).
\end{equation}

From the elementary theory of linear differential equations (see,
e.g.,~\cite[\S 5.5]{boyce86:_elemen_differ_equat}), a particular
solution \( \tilde{g}_p(x) \) to \( \mathbf{\widetilde{L}}\tilde{g} =
\tilde{h} \) in an interval near $x=0$ is given by
\begin{equation*}
  \tilde{g}_p(x) = \sum_{j=1}^{m-1} y_j(x) \int_{t=0}^x
  \frac{\tilde{h}(t) \widetilde{W}_j(t)}{\widetilde{W}(t)} \,dt,
\end{equation*}
where \( y_1, \ldots, y_{m-1} \) form any basis of solutions to \(
\mathbf{\widetilde{L}}\tilde{g} = 0 \), \( \widetilde{W}(x) \) is the Wronskian
\begin{equation*}
  \det\left( \left( y_j^{(i-1)}(x) \right)_{i,j=1,\ldots,m-1} \right),
\end{equation*}
and \( \widetilde{W}_j(x) \) is the determinant obtained from \(
\widetilde{W}(x) \) by replacing the \(j\)th column by the column \(
\mathbf{e}_{m-1} := (0, 0, \dots, 0, 1)^T \).  By Cramer's rule,
the column vector \( (\widetilde{W}_j(x)/\widetilde{W}(x) )_{j=1,\ldots,m-1}
\) is the unique solution \( \mathbf{a} \) to \( \mathcal{W}(x)
\mathbf{a} = \mathbf{e}_{m-1} \), i.e.,
\begin{equation}
  \label{eq:A.5.2}
  \mathbf{a} = \mathcal{W}(x)^{-1} \mathbf{e}_{m-1}
\end{equation}
where \( \mathcal{W}(x) := \left( y_j^{(i-1)}(x) \right)_{i,j=1,\ldots,m-1} \).
The general theory also shows that the initial conditions for \(
\tilde{g}_p \) all vanish, i.e.,
\begin{equation}
  \label{eq:A.5.3}
  \tilde{g}_p^{(k)}(0) = 0, \quad 0 \leq k \leq m-2.
\end{equation}

Here \( y_j(x) = \exp{(\lambda_j x)}\), so that \( y_j^{(i-1)}(x) =
\lambda_j^{i-1} \exp{(\lambda_j x)} \) for \( i,j=1,\ldots,m-1 \).
It follows that
\begin{equation*}
  \mathcal{W}(x)^{-1} = ( v^{ij} e^{-\lambda_i x} )_{i,j=1,\ldots.m-1},
\end{equation*}
where
\begin{equation}
  \label{eq:A.5.5}
  V := (v_{ij}) := ( \lambda_j^{i-1} )_{i,j=1,\ldots.m-1}
\end{equation}
is the Vandermonde matrix and \( V^{-1} =: (v^{ij})_{i,j=1,\ldots,m-1}
\) is its inverse.  The inverse is well known.  By Exercise~1.2.3-40
in~\cite{knuth97} (but note slight differences in definitions),
\begin{align}
  v^{ij} &= [ \lambda^{j-1} ]
  \frac{(\lambda_1-\lambda)\ldots(\lambda_{m-1} - \lambda)}{\lambda_i
    - \lambda} 
  \Biggl/
  \prod_{1 \leq k \leq m-1:\,k \ne i} (\lambda_k - \lambda_i)\notag\\
  &= [ \lambda^{j-1} ] \frac{\psi(\lambda)}{\lambda-\lambda_i}
  \frac1{\psi'(\lambda_i)} \qquad \text{[see~\eqref{eq:B.1}
  and~\eqref{eq:B.3}]}\notag\\
  &= - \frac1{\lambda_i \psi'(\lambda_i)} [ \lambda^{j-1} ] \left(
    \psi(\lambda) \left(1 - \frac{\lambda}{\lambda_i} \right)^{-1}
  \right)\notag\\
  &= - \frac1{\lambda_i \psi'(\lambda_i)} [ \lambda^{j-1} ]
  \left(-m! + \sum_{k=1}^m \stirlingone{m-1}{k} \lambda^k \right)
  \left(\sum_{l=0}^{\infty} \lambda_i^{-l} \lambda^l \right)\notag\\
  &= - \frac1{\lambda_i \psi'(\lambda_i)} \left( -m!\lambda_i^{-(j-1)} +
   \sum_{k=1}^{j-1} \stirlingone{m-1}{k} \lambda_i^{-(j-1-k)}
 \right)\notag\\
   \label{eq:A.5.6}
  &= - \frac1{\lambda_i^{j} \psi'(\lambda_i)} \left( \sum_{k=1}^{j-1}
  \stirlingone{m-1}{k} \lambda_i^k - m! \right) = \frac{1}{\lambda_i^j
  \psi'(\lambda_i)} \sum_{k=j}^{m-1} \stirlingone{m-1}{k} \lambda_i^k.
\end{align}
In particular, \( \mathbf{a} = (a_i)_{i=1,\ldots,m-1} \)
of~\eqref{eq:A.5.2} satisfies
\begin{equation*}
  a_i = v^{i,m-1} e^{-\lambda_i x} = \frac1{\lambda_i^{m-1}
  \psi'(\lambda_i)} \lambda_i^{m-1} e^{-\lambda_i x} =
  \frac1{\psi'(\lambda_i)} e^{-\lambda_i x }.
\end{equation*}

We have thus established
\begin{lemma}
  \label{lemma:A.5.1} The particular solution to \(
  \mathbf{\widetilde{L}}\tilde{g} = \tilde{h} \) with vanishing
  initial conditions (through order \( m-2 \)) is
  \begin{equation*}
    \tilde{g}_p(x) = \sum_{j=1}^{m-1} \frac{e^{\lambda_j
    x}}{\psi'(\lambda_j)} \int_{t=0}^x \tilde{h}(t) e^{-\lambda_j t}\,dt.
  \end{equation*}
\end{lemma}
If desired, this can be transformed back to the \( z \)-domain
via~\eqref{eq:A.2}, \eqref{eq:A.3}, and~\eqref{eq:A.5.1} and a change
of variables in the integral:
\begin{corollary}
  \label{corollary:A.5.2}
  The particular solution to \( \mathbf{L}g=h \) with vanishing initial
  conditions (through order \( m-2 \)) is
  \begin{equation*}
    g_p(z) = \sum_{j=1}^{m-1}
    \frac{(1-z)^{-\lambda_j}}{\psi'(\lambda_j)} \int_{0}^z
    h(\zeta) (1-\zeta)^{\lambda_j+m-2} \,d\zeta.
  \end{equation*}
\end{corollary}

\section*{The initial conditions}
\label{sec:initial-conditions}

Having computed a basis of solutions to the homogeneous equation
and a particular solution to the inhomogeneous equation, so far we have
established Theorem~\ref{thm:diffeq-soln} modulo determination of
the coefficients \( c_1,\ldots,c_{m-1} \) at~\eqref{eq:2.9}.

The fact that the initial conditions for \( \tilde{g}_p \) and for \(
g_p \) vanish make it simple to solve \( \mathbf{L}g = h \) or \(
\mathbf{\widetilde{L}}\tilde{g} = \tilde{h} \) for specified initial
conditions:\ One need only match up the initial conditions of the
homogeneous solutions.

In the $x$-domain the general complementary solution
\begin{equation*}
  \tilde{g}_{\text{c}}(x) = \sum_{j=1}^{m-1} A_j e^{\lambda_j x}
\end{equation*}
has initial conditions
\begin{equation}
\label{eq:A.5.8}
\tilde{g}_{\text{c}}^{(k-1)}(0) = \sum_{j=1}^{m-1} A_j \lambda_j^{k-1} =
\sum_{j=1}^{m-1} v_{kj} A_j, \quad k=1,\ldots,m-1,
\end{equation}
where \( V \) is the Vandermonde matrix at~\eqref{eq:A.5.5}.  In
matrix notation, letting \[ \mathbf{\tilde{g}_0} := (
\tilde{g}_{\text{c}}^{(0)}(0), \ldots, \tilde{g}_{\text{c}}^{(m-2)}(0)
)^T \] and \( \mathbf{A} = ( A_1, \ldots, A_{m-1} )^T \), we can
write~\eqref{eq:A.5.8} as \( \mathbf{\tilde{g}_0} = V \mathbf{A}.  \)
Thus to obtain a specified \( \mathbf{\tilde{g}_0} \) we must choose
\( \mathbf{A} = V^{-1} \mathbf{\tilde{g}_0} \).
Recalling~\eqref{eq:A.5.6},
\begin{equation}
  \label{eq:A.5.100}
  A_j = \frac1{\psi'(\lambda_j)} \sum_{l=1}^{m-1} \left(
  \sum_{i=l}^{m-1} \stirlingone{m-1}{i} \lambda_j^{i-l} \right)
  \tilde{g}_c^{(l-1)}(0), \quad j=1,\ldots,m-1.
\end{equation}
The general complementary solution to $\mathbf{L}g = h$ is then given by
\begin{equation}
  \label{eq:A.12.25}
  g_{\text{c}}(z) = \sum_{j=1}^{m-1} A_j (1-z)^{-\lambda_j}.
\end{equation}

Finally we turn our attention back to \( g \).  According to
Lemma~\ref{lemma:A.1}, \( \mathbf{g_0} = S \mathbf{\tilde{g}_0} \), where
\[ \mathbf{{g}_0} := ( {g}_{\text{c}}^{(0)}(0), \ldots,
{g}_{\text{c}}^{(m-2)}(0) )^T \] are 
the initial conditions for \( g \) and \( S \) is the matrix
\begin{equation*}
  S := (s_{kj}) = \left( \stirlingone{k-1}{j-1} \right)_{k,j=1,\ldots,m-1}
\end{equation*}
of (signless) Stirling numbers  of the first kind.  As is well
known~\cite[\S1.2.6]{knuth97},
\begin{equation}
  \label{eq:A.12.5}
  S^{-1} := (s^{kj}) =  \left( (-1)^{j+k} \stirlingtwo{k-1}{j-1}
  \right)_{k,j=1,\ldots,m-1},
\end{equation}
where \( \stirlingtwo{\cdot}{\cdot} \) denotes (signless) Stirling
numbers of the second kind.  Thus in order for the general
complementary solution~\eqref{eq:A.12.25}
to have a specified vector \( \mathbf{g_0} \) of initial conditions, we
must take \( \mathbf{A} = V^{-1} \mathbf{\tilde{g}_0} = V^{-1} S^{-1}
\mathbf{g_0} \).  We can now establish~\eqref{eq:2.9}.  Observe that
if \( g_{\text{c}}(z) \) is specified as a power series
\[
g_{\text{c}}(z) = \sum_{k=0}^\infty g_k z^k
\]
then \( g_{\text{c}}^{(k)}(0) = k! g_k
\) for \( k = 0, \ldots, m-2 \).

\begin{proposition}
  \label{lemma:A.2.1}
  The constants \( A_j \) at~\eqref{eq:A.5.100} are given by
  \begin{equation*}
    A_j = \frac{m!}{\psi'(\lambda_j)} \sum_{k=0}^{m-2} 
    \frac{g_c^{(k)}(0)}{\rising{\lambda_j}{k+1}}, \quad j=1,\ldots,m-1.
  \end{equation*}
\end{proposition}
\begin{proof}
  Observe that, by~\eqref{eq:A.5.6} and~\eqref{eq:A.12.5},
  \begin{align*}
    A_j &= \sum_{k=1}^{m-1} \sum_{l=1}^{m-1} v^{jl}
    s^{lk} g_{\text{c}}^{(k-1)}(0) = \sum_{k=1}^{m-1} \sum_{l=1}^{m-1}
    \frac1{\psi'(\lambda_j)} \sum_{i=l}^{m-1}
    \stirlingone{m-1}{i} \lambda_j^{i-l} (-1)^{l+k}
    \stirlingtwo{l-1}{k-1} g_{\text{c}}^{(k-1)}(0)\\
    &= \frac{1}{\psi'(\lambda_j)} \sum_{l=1}^{m-1} \sum_{i=l}^{m-1}
    \stirlingone{m-1}{i}
    \lambda_j^{i-l}  \sum_{k=0}^{m-2} (-1)^{l+k+1}
    \stirlingtwo{l-1}{k} g_{\text{c}}^{(k)}(0).
  \end{align*}
  Hence it is enough to establish the curious identity
  \begin{equation}
    \label{eq:A.2.2}
    \sum_{l=1}^{m-1} \sum_{i=l}^{m-1} \stirlingone{m-1}{i}
    \lambda_j^{i-l} (-1)^{l+k+1} \stirlingtwo{l-1}{k} =
    \frac{m!}{\rising{\lambda_j}{k+1}}
  \end{equation}
  for $k=0,\ldots,m-2$.
  
  For this, denote the left-hand side of~\eqref{eq:A.2.2} by \( f_k \)
  (fixing \( m \) and \( j \)).  First, since \(\stirlingtwo{n}{0} \equiv
  \delta_{n,0} \), we have
  \begin{equation}
    \label{eq:A.2.3}
    f_0 = \sum_{i=1}^{m-1} \stirlingone{m-1}{i} \lambda_j^{i-1} =
    \frac{\rising{\lambda_j}{m-1}}{\lambda_j} =
    \frac{m!}{\lambda_j}
  \end{equation}
  as claimed.  For \( k \geq 1 \), using the recurrence
  relation~\citep[\S1.2.6]{knuth97}
  \begin{equation*}
    \stirlingtwo{n}{k} = k \stirlingtwo{n-1}{k} +
    \stirlingtwo{n-1}{k-1}, \quad k \geq 1
  \end{equation*}
  for Stirling numbers of the second kind and the
  connection~\cite[1.2.6-(47)]{knuth97}  between
  the two kinds of Stirling numbers
  we see that
  \begin{equation}
    \label{eq:A.2.4}
    f_k = -\frac{k}{\lambda_j} f_k + \frac1{\lambda_j} f_{k-1} +
    \delta_{m-1,k}, \quad
    k \geq 1.
  \end{equation}
  With initial condition~\eqref{eq:A.2.3} the
  recurrence~\eqref{eq:A.2.4} can be solved to give \( f_k =
  m!/\rising{\lambda_j}{k+1} \) for \( k = 0,\ldots,m-2 \), as desired.
\end{proof}


\chapter[Identities involving the indicial polynomial]{Identities
  involving the indicial polynomial}
\label{appendix:indicial}

The indicial polynomial
\begin{equation}
  \label{eq:B.0}
  \psi(\lambda) \equiv \psi_m(\lambda) := \rising{\lambda}{m-1} - m!
\end{equation}
plays an important role in the analysis of \( m \)-ary search trees
under the random permutation model.  We will enumerate a few useful
identities involving the polynomial in this appendix.

It is well known~\cite[Chapter~3]{MR93f:68045} that \( \psi_m \) has
\(m-1\) distinct roots \( 2 = \lambda_1, \lambda_2, \ldots,
\lambda_{m-1} \) listed in nonincreasing order of real part.
As in~\cite[Chapter~3]{MR93f:68045} we introduce
\begin{equation*}
  \alpha := \max_{2 \leq j \leq m-1} \Re{(\lambda_j)};
\end{equation*}
that is, \(\alpha\) is the real part of the root having the second
largest real part among all the roots of the indicial
polynomial~\eqref{eq:B.0}.  We list some important properties of the
roots of~\eqref{eq:B.0} stated in~\cite[\S3.3]{MR93f:68045}:
\begin{enumerate}
\item The number \(-m\) is a root if and only if \(m\) is odd.  
  All other roots of \(\psi(\lambda)\) are simple, non-real roots.
\item No two roots have the same real part unless they are
  mutually conjugate.  (This follows from the strict increasingness of
  $|\rising{(s+it)}{m-1}|$ in $|t|$.)
\end{enumerate}

The indicial polynomial can be expressed as
\begin{equation}
  \label{eq:B.1}
  \psi(\lambda) = \prod_{j=1}^{m-1} (\lambda-\lambda_j),
\end{equation}
so that
\begin{equation}
  \label{eq:B.2}
  \psi'(\lambda) = \sum_{j=1}^{m-1} \prod_{1 \leq k \leq m-1:\,k \ne
  j} (\lambda-\lambda_k).
\end{equation}
If \( \lambda = \lambda_j \) for some \( j=1,\ldots,m-1 \),
then from~\eqref{eq:B.2} we get
\begin{equation}
  \label{eq:B.3}
  \psi'(\lambda_j) = \prod_{1 \leq k \leq m-1:\,k \ne j}
  (\lambda_j-\lambda_k)
\end{equation}
since all other terms in the sum for $\psi'(\lambda_j)$ contain a
factor that is~0.
\begin{identity}
  \label{identity:B.1}
  When \( \lambda \notin \{ \lambda_1,\ldots,\lambda_{m-1} \} \),
  \begin{equation*}
    \sum_{j=1}^{m-1} \frac1{(\lambda-\lambda_j)\psi'(\lambda_j)} =
    \frac{1}{\psi(\lambda)}
  \end{equation*}
\end{identity}
\begin{proof}
  Using~\eqref{eq:B.1} and partial fraction expansion we get, for \(
  \lambda \notin \{ \lambda_1,\ldots,\lambda_{m-1} \} \),
  \begin{equation*}
    \frac{1}{\psi(\lambda)} = \frac{1}{\prod_{j=1}^m
      (\lambda-\lambda_j)} = \sum_{j=1}^{m-1} \frac{A_j}{\lambda-\lambda_j}
  \end{equation*}
  where, using~\eqref{eq:B.3},
  \begin{equation*}
    A_j = \frac{1}{\prod_{1 \leq k \leq m-1:\,k \ne j}
    (\lambda_j-\lambda_k)} = \frac{1}{\psi'(\lambda_j)},
  \end{equation*}
  as required.
\end{proof}
For \( r \) and \( n \) positive integers, let \( H_n^{(r)} \) denote
the $r$th-order harmonic number
\begin{equation*}
  \label{eq:B.4}
  H_n^{(r)} := \sum_{j=1}^n \frac{1}{j^r}.
\end{equation*}
When \( r=1 \) we will use \( H_n := H_n^{(1)} \) for the usual
(1st-order) harmonic number.
\begin{identity}
  \label{identity:B.2}
  \begin{equation*}
    \sum_{j=1}^{m-1}\frac{1}{(\lambda_j+r)^2\psi'(\lambda_j)} = 
    \begin{cases}
      \displaystyle
       \frac{(m-1)!H_{m-1}}{[(m-1)!]^2(m-1)^2} =
      \frac{H_{m-1}}{(m-1)!(m-1)^2} & r = -1\\\\
      \displaystyle\frac{(-1)^{r}r!(m-2-r)!}{[\rising{(-r)}{m-1}-m!]^2}
      & 0 \leq r \leq m-2\\\\
      \displaystyle \frac{(-1)^{m}\falling{r}{m-1}(H_r -
      H_{r-(m-1)})}{[(-1)^{m-1}\falling{r}{m-1} - m!]^2} & r \geq m-1.
    \end{cases}
  \end{equation*}
\end{identity}
\begin{proof}
  Differentiate Identity~\ref{identity:B.1} with respect to \( \lambda
  \) to get
  \begin{equation}
    \label{eq:B.30}
    \sum_{j=1}^{m-1} \frac1{(\lambda-\lambda_j)^2 \psi'(\lambda_j)}
     = \frac{\psi'(\lambda)}{\psi^2(\lambda)},
    \qquad \lambda \notin \{\lambda_1,\ldots,\lambda_{m-1}\}.
  \end{equation}
  Differentiating the defining equation~\eqref{eq:B.0} we get
  \begin{equation}
    \label{eq:B.20}
    \psi'(\lambda) = \sum_{j=0}^{m-2} \prod_{0 \leq k \leq m-2:\,k \ne
    j} (\lambda + k ) = \rising{\lambda}{m-1} \sum_{j=0}^{m-2}
    \frac1{\lambda+j},
  \end{equation}
  the last equality being valid only for \(
  \lambda \notin \{0,-1,\ldots,-(m-2)\} \).
  Hence for \( r \geq -1 \),
  \begin{equation*}
    \sum_{j=1}^{m-1} \frac1{(\lambda_j+r)^2 \psi'(\lambda_j)}
     = \frac{\psi'(-r)}{\psi^2(-r)} =
     \frac{\psi'(-r)}{(\rising{(-r)}{m-1}-m!)^2}.
  \end{equation*}
  Now using~\eqref{eq:B.20} we see that
  \begin{equation*}
    \psi'(1) = \rising{1}{m-1} \sum_{j=0}^{m-2} \frac{1}{1+j} =
    (m-1)!H_{m-1},
  \end{equation*}
  \begin{equation*}
    \psi'(-r) = \prod_{0 \leq j\leq m-2:\,j \ne r} (j-r) = (-1)^r
    r! (m-2-r)!, \qquad 0 \leq r \leq m-2,
  \end{equation*}
  \begin{equation*}
    \psi'(-r) = \rising{(-r)}{m-1} \sum_{j=0}^{m-2}
    \frac1{j-r} = (-1)^m \falling{r}{m-1} (H_r - H_{r-(m-1)}), \qquad
    r \geq m-1,
  \end{equation*}
  and the claim follows.
\end{proof}
\begin{identity}
  \label{identity:B.2.2}
  For \( 0 \leq k \leq m-3 \),
  \begin{equation*}
    \sum_{j=1}^{m-1} \frac{\lambda_j^k}{\psi'(\lambda_j)} = 0.
  \end{equation*}
\end{identity}
\begin{proof}
  We take asymptotic expansions of both sides of
  Identity~\ref{identity:B.1} as \( \lambda \to \infty \) in powers of
  \( 1/\lambda \) and equate coefficients.  The left-hand side
  yields
  \begin{equation*}
    \frac1\lambda \sum_{j=1}^{m-1} \frac1{\left( 1 -
    \frac{\lambda_j}{\lambda}\right) \psi'(\lambda_j)} = \frac1\lambda
    \sum_{j=1}^{m-1} \frac1{\psi'(\lambda_j)} \sum_{k=0}^\infty
    \frac{\lambda_j^k}{\lambda^k} = \sum_{k=0}^\infty
    \frac1{\lambda^{k+1}} \sum_{j=1}^{m-1}
    \frac{\lambda_j^k}{\psi'(\lambda_j)}.
  \end{equation*}
  On the other hand, the right-hand side equals
  \begin{equation*}
    \frac{1}{\rising{\lambda}{m-1}-m!} = \left(\prod_{k=0}^{m-2}(\lambda+k) -
      m!\right)^{-1} = \frac1{\lambda^{m-1}} \left( \prod_{k=0}^{m-2} \left(1 +
        \frac{k}{\lambda} \right) - \frac{m!}{\lambda^{m-1}} \right)^{-1}.
  \end{equation*}
  Hence the coefficient of \( 1/\lambda^{k+1} \) for \( k = 0,\ldots,m-3
  \) in the left-hand side must vanish, proving the claim.
\end{proof}
\begin{identity}
  \label{identity:B.5.22}
  \begin{equation*}
    \sum_{j=2}^{m-1} \frac1{(\lambda_j-2)\psi'(\lambda_j)} =
    \frac1{2(m!)} \left[ 1 - \frac{H_m^{(2)}-1}{(H_m-1)^2} \right].    
  \end{equation*}
\end{identity}
\begin{proof}
  When \( \lambda_k \) is a root of \( \psi \) we have
  from~\eqref{eq:B.20}
  \begin{equation*}
    \psi'(\lambda_k) = m! \sum_{j=0}^{m-2} \frac1{\lambda_k+j}.
  \end{equation*}
  Also, differentiating~\eqref{eq:B.20} once more gives us
  \begin{equation}
    \label{eq:B.7}
    \psi''(\lambda) = \rising{\lambda}{m-1} \left[ \left(
    \sum_{j=0}^{m-2} \frac{1}{\lambda+j} \right)^2 - \sum_{j=0}^{m-2}
    \frac1{(\lambda+j)^2} \right],
  \end{equation}
  so that
  \begin{equation}
    \label{eq:psp2}
    \psi'(2) = m!( H_m-1 ) \qquad\text{ and }\qquad \psi''(2) = m![
    (H_m-1)^2 -  (H_m^{(2)}-1) ].
  \end{equation}
  Thus, as \( \lambda \to 2 \),
  \begin{equation*}
    \psi(\lambda) = (\lambda-2) m!( H_m-1 ) + \frac12 m![
    (H_m-1)^2 -  (H_m^{(2)}-1) ] (\lambda-2)^2 + O(|\lambda-2|^3).
  \end{equation*}
  By Identity~\ref{identity:B.1}, for $\lambda \notin
  \{\lambda_1,\ldots,\lambda_{m-1}\}$ we have
  \begin{align*}
    &\sum_{j=2}^{m-1} \frac{1}{(\lambda - \lambda_j) \psi'(\lambda_j)} \\
    &=\frac1{\psi(\lambda)} - \frac1{(\lambda-2)m!(H_m-1)}\\
    &= \frac1{(\lambda-2)m!(H_m-1)} \left\{ \left[ 1 + \frac12 \left[
          (H_m-1) - \frac{H_m^{(2)}-1}{H_m-1} \right](\lambda-2) +
        O(|\lambda-2|^2) \right]^{-1} - 1 \right\}\\
    &= \frac1{(\lambda-2)m!(H_m-1)} \left \{ -\frac12 \left[ (H_m-1) -
    \frac{H_m^{(2)}-1}{H_m-1} \right](\lambda-2) + O(|\lambda-2|^2)
    \right\}\\
    &= -\frac1{2(m!)}\left[ 1 - \frac{H_m^{(2)}-1}{(H_m-1)^2} \right] +
    O(|\lambda-2|);
  \end{align*}
  in particular, letting $\lambda \to 2$ yields
  \begin{equation*}
    \sum_{j=2}^{m-1} \frac1{(\lambda_j-2)\psi'(\lambda_j)} =
    \frac1{2(m!)}\left[ 1 - \frac{H_m^{(2)}-1}{(H_m-1)^2} \right],
  \end{equation*}
  as claimed.
\end{proof}

For the next two identities we assume that \( m \) is odd.  Note that
when \( m \) is odd, \(-m\) is a root of \( \psi_m \)---in fact, \(
\lambda_{m-1} = -m \), as noted at the beginning of this appendix.
\begin{identity}
  \label{identity:B.3}
  When \( m \) is odd,
  \begin{equation*}
    \sum_{j=1}^{m-2} \frac1{(\lambda_j+m)\psi'(\lambda_j)} =
    \frac1{2(m!)} \left[ 1 - \frac{H_m^{(2)}-1}{(H_m-1)^2}\right].
  \end{equation*}
\end{identity}
\begin{proof}
  By Identity~\ref{identity:B.1} the left-hand side above equals
  \begin{align*}
    & -\lim_{\lambda \to -m} \left[ \frac1{\psi(\lambda)} -
      \frac1{(\lambda+m)\psi'(-m)} \right]\\
    &=  -\lim_{\lambda \to -m} \left[ \frac1{(\lambda + m)\psi'(-m) +
        \tfrac12(\lambda+m)^2\psi''(-m) + O(|\lambda+m|^3)} -
      \frac1{(\lambda+m)\psi'(-m)} \right]\\
    &= \lim_{\lambda \to -m} \frac1{(\lambda+m)\psi'(-m)} \left[
      \frac12(\lambda+m) \frac{\psi''(-m)}{\psi'(-m)} +
      O(|\lambda+m|^2) \right]
    = \frac12 \frac{\psi''(-m)}{[\psi'(-m)]^2}.
  \end{align*}
  Now by~\eqref{eq:B.20},
  \begin{equation}
    \label{eq:psipmm}
   \psi'(-m) = -m!(H_m-1),
  \end{equation}
   and by~\eqref{eq:B.7}
  \begin{equation}
    \label{eq:psippmm}
    \psi''(-m) = m! \left[ \left( \sum_{j=0}^{m-2} \frac1{j-m}
    \right)^2 - \sum_{j=0}^{m-2} \frac1{(j-m)^2} \right] = m![
    (H_m-1)^2 - (H_m^{(2)}-1)],
  \end{equation}
  and the result follows.
\end{proof}
\begin{identity}
  \label{identity:B.4}
  When \(m\) is odd
  \begin{equation*}
  \sum_{j=1}^{m-2} \frac1{(\lambda_j+m)^2\psi'(\lambda_j)} =
  \frac1{m!}\left[ \frac1{12} (H_m-1) - \frac13 \frac{ H_m^{(3)} -
  1}{(H_m-1)^2} + \frac14 \frac{(H_m^{(2)}-1)^2}{(H_m-1)^3} \right].
\end{equation*}
\end{identity}
\begin{proof}
  By~\eqref{eq:B.30} the left-hand-side above equals
  \begin{align*}
    & \lim_{\lambda \to -m}
    \left[\frac{\psi'(\lambda)}{\psi^2(\lambda)} -
      \frac{1}{(\lambda+m)^2\psi'(-m)}\right]\\
    &= \lim_{\lambda\to-m}
    \Biggl[ \frac{\psi'(-m) + \psi''(-m)(\lambda+m) +
      \tfrac{1}{2}\psi'''(-m)(\lambda+m)^2 +
      O(|\lambda+m|^3)}{[\psi'(-m)(\lambda+m) +
      \tfrac{1}{2}\psi''(-m)(\lambda+m)^2 +
      \tfrac{1}{6}\psi'''(-m)(\lambda+m)^3 + 
      O(|\lambda+m)^4|]^2}\\
    & \qquad - \frac{1}{(\lambda+m)^2\psi'(-m)}\Biggr]\\
    &= \lim_{\lambda\to-m}\biggl[\frac{1}{(\lambda+m)^2\psi'(-m)}\\
      & \qquad \times
      \left(\frac{1 + \tfrac{\psi''(-m)}{\psi'(-m)}(\lambda+m) +
          \tfrac{1}{2}\frac{\psi'''(-m)}{\psi'(-m)}(\lambda+m)^2 +
          O(|\lambda+m|^3)}{[1 
          + \tfrac{1}{2}\tfrac{\psi''(-m)}{\psi'(-m)}(\lambda+m) +
          \tfrac{1}{6}\tfrac{\psi'''(-m)}{\psi'(-m)}(\lambda+m)^2 +
          O(|\lambda+m|^3)]^2}
        -1\right)\biggr]\\
    &= \lim_{\lambda\to-m}\Biggl[\frac{1}{(\lambda+m)^2\psi'(-m)} \\
    & \qquad \times
    \left(\frac{1 + \tfrac{\psi''(-m)}{\psi'(-m)}(\lambda+m)
        + \tfrac{1}{2}\tfrac{\psi'''(-m)}{\psi'(-m)}(\lambda+m)^2 +
        O(|\lambda+m|^3)}{1
        + \tfrac{\psi''(-m)}{\psi'(-m)}(\lambda+m) +
        \bigl[\tfrac{1}{4}(\tfrac{\psi''(-m)}{\psi'(-m)})^2 +
        \tfrac{1}{3}\tfrac{\psi'''(-m)}{\psi'(-m)}\bigr](\lambda+m)^2 +
        O(|\lambda+m|^3)}
      - 1 \right)\Biggr]\\
    &= \lim_{\lambda\to-m}
    \Biggr[\frac{1}{(\lambda+m)^2\psi'(-m)}\\
    & \times \Biggl(\left[1 +
      \tfrac{\psi''(-m)}{\psi'(-m)}(\lambda+m) +
      \tfrac{1}{2}\tfrac{\psi'''(-m)}{\psi'(-m)}(\lambda+m)^2 +
      O(|\lambda+m|^3)\right] \\
    & \times
    \left[ 1 - \frac{\psi''(-m)}{\psi'(-m)}(\lambda+m)
      +\left[\frac{3}{4}\left(\frac{\psi''(-m)}{\psi'(-m)}\right)^2 -
        \frac{1}{3} \frac{\psi'''(-m)}{\psi'(-m)}\right](\lambda+m)^2
       + O(|\lambda+m|^3)\right] - 1 \Biggr) \Biggr]\\
     &= \lim_{\lambda\to-m}\left[\frac{1}{(\lambda+m)^2\psi'(-m)}
       \left(\left(\frac{1}{6}\frac{\psi'''(-m)}{\psi'(-m)} -
           \frac{1}{4}\left(
             \frac{\psi''(-m)}{\psi'(-m)}\right)^2\right)(\lambda+m)^2
         + O(|\lambda+m|^3)\right) \right]\\
     &= \frac{1}{6}\frac{\psi'''(-m)}{[\psi'(-m)]^2} -
     \frac{1}{4}\frac{[\psi''(-m)]^2}{[\psi'(-m)]^3}. 
\end{align*}
We have already computed $\psi'(-m)$ and $\psi''(-m)$
at~\eqref{eq:psipmm} and~\eqref{eq:psippmm}, respectively.  To
calculate \( \psi'''(-m) \) we differentiate~\eqref{eq:B.7} to get,
for \( \lambda \notin \{0,-1,\ldots,-(m-2)\} \),
\begin{align*}
  \psi'''(\lambda) &= \rising{\lambda}{m-1} \Biggl[ 2 \sum_{j=0}^{m-2}
    \frac1{(\lambda+j)^3} - 2 \left( \sum_{j=0}^{m-2} \frac1{\lambda+j}
    \right) \left( \sum_{j=0}^{m-2} \frac1{(\lambda+j)^2} \right)\\
    & \qquad - \left( \sum_{j=0}^{m-2} \frac1{\lambda+j} \right) \left(
      \sum_{j=0}^{m-2} \frac1{(\lambda+j)^2} \right) + \left(
      \sum_{j=0}^{m-2} \frac1{\lambda+j} \right)^3 \Biggr]\\
  &= \rising{\lambda}{m-1} \left[ 2 \sum_{j=0}^{m-2}
    \frac1{(\lambda+j)^3} - 3 \left( \sum_{j=0}^{m-2}
    \frac1{\lambda+j} \right) \left( \sum_{j=0}^{m-2}
    \frac1{(\lambda+j)^2} \right) + \left(
    \sum_{j=0}^{m-2} \frac1{\lambda+j} \right)^3 \right].
\end{align*}
In particular, since \(-m\) is a root of \( \psi \) when \( m \) is
odd,
\begin{equation*}
  \psi'''(-m) = m! [ -2 (H_m^{(3)}-1) + 3(H_m-1)(H_m^{(2)}-1) - (H_m
  -1 )^3].
\end{equation*}
Therefore
\begin{align*}
  \sum_{j=1}^{m-2} \frac1{(\lambda_j+m)^2\psi'(\lambda_j)} &=
  \frac{1}{6}\frac{m![-2(H_m^{(3)}
    - 1) + 3(H_m-1)(H_m^{(2)}-1) - (H_m-1)^3]}{(m!)^2(H_m-1)^2}\\
  & \qquad +
  \frac{1}{4}\frac{[(H_m-1)^2 - (H_m^{(2)}-1)]^2}{m!(H_m-1)^3}\\
  &= \frac{1}{6(m!)}\Bigl[-(H_m-1) + \frac{3(H_m^{(2)}-1)}{H_m-1} -
  \frac{2(H_m^{(3)} -1)}{(H_m-1)^2}\Bigr]\\
  & \qquad + \frac{1}{4 (m!)}\Bigl[ (H_m-1) - \frac{2(H_m^{(2)}-1)}{H_m-1} +
  \frac{(H_m^{(2)} - 1)^2}{(H_m-1)^3}\Bigr]\\
  &= \frac{1}{m!}\Bigl[\frac{1}{12}(H_m-1) -
  \frac{1}{3}\frac{H_m^{(3)}-1}{(H_m-1)^2} +
  \frac{1}{4}\frac{(H_m^{(2)}-1)^2}{(H_m-1)^3}\Bigr],
\end{align*}
as claimed.
\end{proof}


\backmatter

\pagestyle{plain}

\bibliographystyle{habbrv}
\bibliography{msn,leftovers}


\chapter{Vita}

Nevin Kapur was born on May~9, 1974 in Mumbai (formerly Bombay),
India.

He received a Bachelor of Technology in Electrical Engineering, with a
specialization in \textsc{VLSI} design and automation, from the Indian
Institute of Technology, Bombay in 1996.  From 1996 to 1998 he was
part of the Collaborative Benchmarking Laboratory at North Carolina
State University (\textsc{NCSU}).  He received a Master of Science
from the Electrical
and Computer Engineering department, specializing in graph algorithms
and their applications to automatic place-and-route in \textsc{VLSI}
design.  At \textsc{NCSU}, while taking courses in design and analysis
of algorithms, graph theory, and probability, he became interested in
these topics and their intersection.  This motivated his joining the
Department of Mathematical Sciences at The Johns Hopkins University in
1998.

Nevin plans on pursuing an academic career, as a researcher pursuing
problems in discrete mathematics, probability, and their applications
to computer science.

\bigskip\noindent
\textit{Baltimore, 2003}


\end{document}